\documentclass[a4paper,12pt]{article}
\usepackage{amsfonts}
\usepackage{amsmath}
\usepackage{amssymb}
\usepackage{graphicx}
\usepackage{epsf}
\usepackage{epsfig}

\usepackage{color}
\usepackage{afterpage}
\usepackage{verbatim}

\usepackage{times}

\newtheorem{theo}{Theorem}

\begin{document}

\title{An adaptive finite element method  in  reconstruction of coefficients in Maxwell's equations from limited observations}

\author{L. Beilina  \thanks{
Department of Mathematical Sciences, Chalmers University of Technology and
Gothenburg University, SE-42196 Gothenburg, Sweden, e-mail: \texttt{\
larisa@chalmers.se}}
 \and S. Hosseinzadegan \thanks{
Department of Mathematical Sciences, Chalmers University of Technology and
Gothenburg University, SE-42196 Gothenburg, Sweden,e-mail: \texttt{\
samarh{\@@}student.chalmers.se}}}

\date{}

\maketitle


\begin{abstract}
We propose an adaptive finite element method for the solution of a
coefficient inverse problem of simultaneous reconstruction of the
dielectric permittivity and magnetic permeability functions in the
Maxwell's system using limited boundary observations of the electric
field in 3D.

We derive a posteriori error estimates in the Tikhonov functional to be
minimized and in the regularized solution of this functional,  as well as
formulate corresponding adaptive algorithm.  Our numerical experiments
justify the efficiency of our a posteriori estimates and show
significant improvement of the reconstructions obtained on  locally 
adaptively  refined meshes.

\end{abstract}

\section{Introduction}
\label{sec:intro}

This work is a continuation of the recent paper \cite{BCN} and is
focused on the numerical reconstruction of the dielectric permittivity
$\varepsilon(x)$ and the magnetic permeability $\mu(x)$ functions in
the Maxwell's system on  locally refined meshes using an adaptive
finite element method. The reconstruction is performed via
minimization of the corresponding Tikhonov functional from
backscattered single measurement data of the electric field $E(x,t)$.
That means that we use backscattered boundary measurements of the wave
field $E(x,t)$ which are generated by a single direction of a plane
wave. In the minimization procedure we use domain decomposition finite
element/finite difference methods of \cite{BMaxwell} for the numerical
reconstructions of both functions.

Comparing with \cite{BCN}  we present following new points here: we
adopt results of \cite{BOOK, BKK,KBB} to show that the minimizer of the
Tikhonov functional is closer to the exact solution than  guess
of this solution. We  present relaxation property for the mesh
refinements for the case of our inverse problem and we
derive a posteriori error estimates for the error in the minimization
functional and in the reconstructed functions $\varepsilon(x)$ and
$\mu(x)$. Further, we formulate two adaptive algorithms and apply them in the
reconstruction of small inclusions. Moreover, in our numerical
simulations of this work we induce inhomogeneous initial conditions in the
 Maxwell's system.  Non-zero initial conditions involve
uniqueness and stability results of reconstruction of both unknown functions $\varepsilon(x)$ and $\mu(x)$, see details in
\cite{BCN, BCS}.

Using our numerical simulations we can conclude that an adaptive finite
element method can significantly improve reconstructions obtained on a
coarse non-refined mesh in order to accurately obtain shapes,
locations and values of functions $\varepsilon(x)$ and $\mu(x)$.

An outline of this paper is as follows: in Section \ref{sec:model} we
present our mathematical model and in Section \ref{sec:stat} we
formulate forward and inverse problems. In Section \ref{sec:tikhonov}
we present the Tikhonov functional to be minimized and in Section
\ref{sec:spaces} we show different versions of finite element method
used in computations. In Section \ref{sec:relax} we formulate
relaxation property of mesh refinements and in Section
\ref{sec:general} we investigate general framework of a posteriori
error estimates in coefficient inverse problems (CIPs). In Sections
\ref{sec:adaptrelax}, \ref{sec:errorfunc} we present theorems for a
posteriori errors in the regularized solution of the Tikhonov
functional and in the Tikhonov functional, correspondingly. In Sections
\ref{sec:ref}, \ref{sec:alg} we describe mesh refinement
recommendations and formulate adaptive algorithms used in
computations. Finally, in Section \ref{sec:num} we present our
reconstruction results.

\section{The mathematical model}

\label{sec:model}

Let a bounded domain $\Omega \subset \mathbb{R}^d, d=2,3,$  have Lipschitz
boundary $\partial \Omega$ and let us  set $\Omega_T := \Omega \times
(0,T)$, $\partial \Omega_T := \partial \Omega \times (0,T)$,  where
$T >0$.
We consider Maxwell's equations 
 in an inhomogeneous isotropic media in a bounded domain  $\Omega \subset  \mathbb{R}^3$  
\begin{equation}
  \label{eq:maxwell}
  \left \{ \begin{array}{llllll} 
\partial_{t} D - \nabla \times H(x,t) = 0 &&\mbox{ in } \Omega_T\\
 \partial_{t} B + \nabla \times E(x,t) = 0 &&
 \mbox{ in } \Omega_T,\\
 D(x,t)= \varepsilon E(x,t), \quad  B(x,t)= \mu H(x,t),&&\\
E(x,0) = E_0(x), \quad  H(x,0) = H_0(x), &&\\
\nabla \cdot D(x,t) = 0,\quad \nabla \cdot B(x,t)  =0 && \mbox{ in } \Omega_T,\\
      n \times   D(x,t) =0,\quad n \cdot  B(x,t) =0 && \textrm{on }\,\partial \Omega_T, \end{array} \right .
\end{equation}
where $x=(x_1,x_2,x_3)$. Here, $E(x,t)$ is the electric field and
$H(x,t)$ is the magnetic field, $\varepsilon(x) > 0$ and $\mu(x) >0$
are the dielectric permittivity and the magnetic permeability
functions, respectively. $E_0(x)$ and $H_0(x)$ are given initial
conditions. Next, $n = n(x)$ is the unit outward
normal vector to $\partial \Omega$. The electric field $E(x,t)$ is
combined with the electric induction $D(x,t)$ via
\begin{equation*}
D(x,t) = \varepsilon E(x,t) = \varepsilon_{\rm vac} \varepsilon_r E(x,t),
\end{equation*}
where $\varepsilon_{\rm vac} \approx 8.854 \times 10^{-12}$ is the
 vacuum permittivity which is measured in Farads per meter,
and thus $\varepsilon_r$ is the dimensionless relative permittivity.
The magnetic field $H(x,t)$ is combined with the magnetic induction $B(x,t)$
via
\begin{equation*}
B(x,t)  = \mu H(x,t) = \mu_{\rm vac} \mu_r H(x,t),
\end{equation*}
where $\mu_{\rm vac} \approx 1.257 \times 10^{-6}$ is the vacuum permeability
measured in Henries per meter, from what follows that $\mu_r$ is the
dimensionless relative permeability.

By eliminating $B$ and $D$ from (\ref{eq:maxwell}) we obtain
the model problem for the electric field $E$ with the
perfectly conducting boundary conditions
 which is as follows:
\begin{eqnarray}
\varepsilon \frac{\partial^2 E}{\partial t^2} + \nabla \times ( \mu^{-1} \nabla \times E)    &=& 0 ~ \mbox{in}~~ \Omega_T,    \label{model1_1} \\
\nabla \cdot (\varepsilon E) &=& 0 ~ \mbox{in}~~ \Omega_T,  \label{model1_2}  \\
  E(x,0) = f_0(x), ~~~E_t(x,0) &=& f_1(x)~ \mbox{in}~~ \Omega,      \label{model1_3}  \\
E \times n &=& 0 ~ \mbox{on}~~ \partial \Omega_T.  \label{model1_4} 
\end{eqnarray}
Here we assume that
\begin{equation*}
f_{0}\in H^{1}(\Omega), f_{1}\in L^{2}(\Omega). 
\end{equation*}
By this notation we shall mean that every component of the vector
functions $f_0$  and $f_1$  belongs to these spaces. Note that equations similar to  (\ref{model1_1})-(\ref{model1_4}) can be derived also for the magnetic
field $H$.

As in our recent work \cite{BCN}, for the discretization of the
Maxwell's equations we use a stabilized domain decomposition method of
\cite{BMaxwell2}. 
In our numerical simulations we assume
that  the relative permittivity $\varepsilon_r$ and relative
permeability $\mu_r$ does not vary much which is the case of real
applications, see recent experimental work \cite{ BTKB} for
similar observations. We do not impose smoothness
assumptions on the
coefficients $\varepsilon(x), \mu(x)$ and we treat
discontinuities  in a similar way as  in \cite{CWZ14}.
Thus, a discontinuous finite element method should be applied for the
finite element discretization of these functions, see details in
Section \ref{sec:spaces}.  




\section{Statements of forward and inverse problems}

\label{sec:stat}

\begin{figure}[tbp]
 \begin{center}
 \begin{tabular}{cc}
{\includegraphics[width=7.0cm, clip = true, trim = 6.0cm 6.0cm 6.0cm 6.0cm]{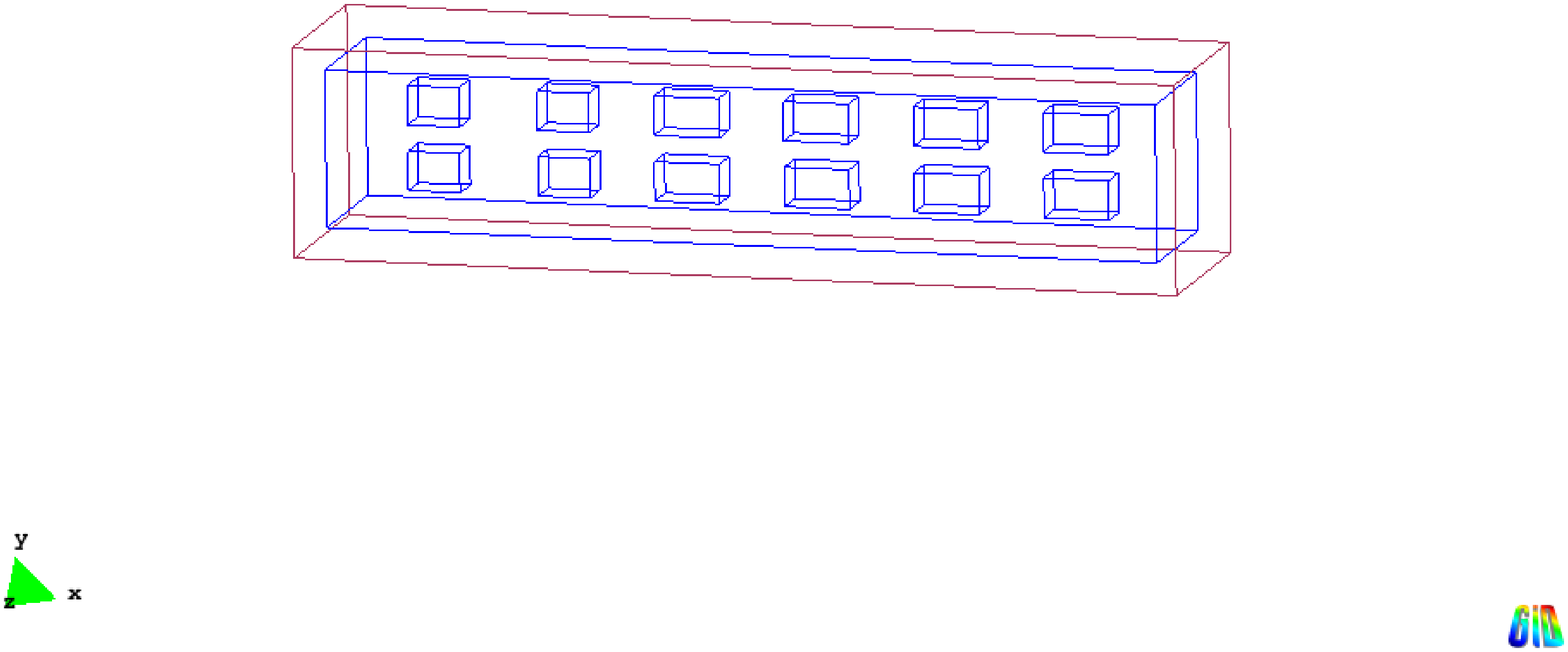}} &
{\includegraphics[width=7.0cm, clip = true, trim = 6.0cm 6.0cm 6.0cm 6.0cm]{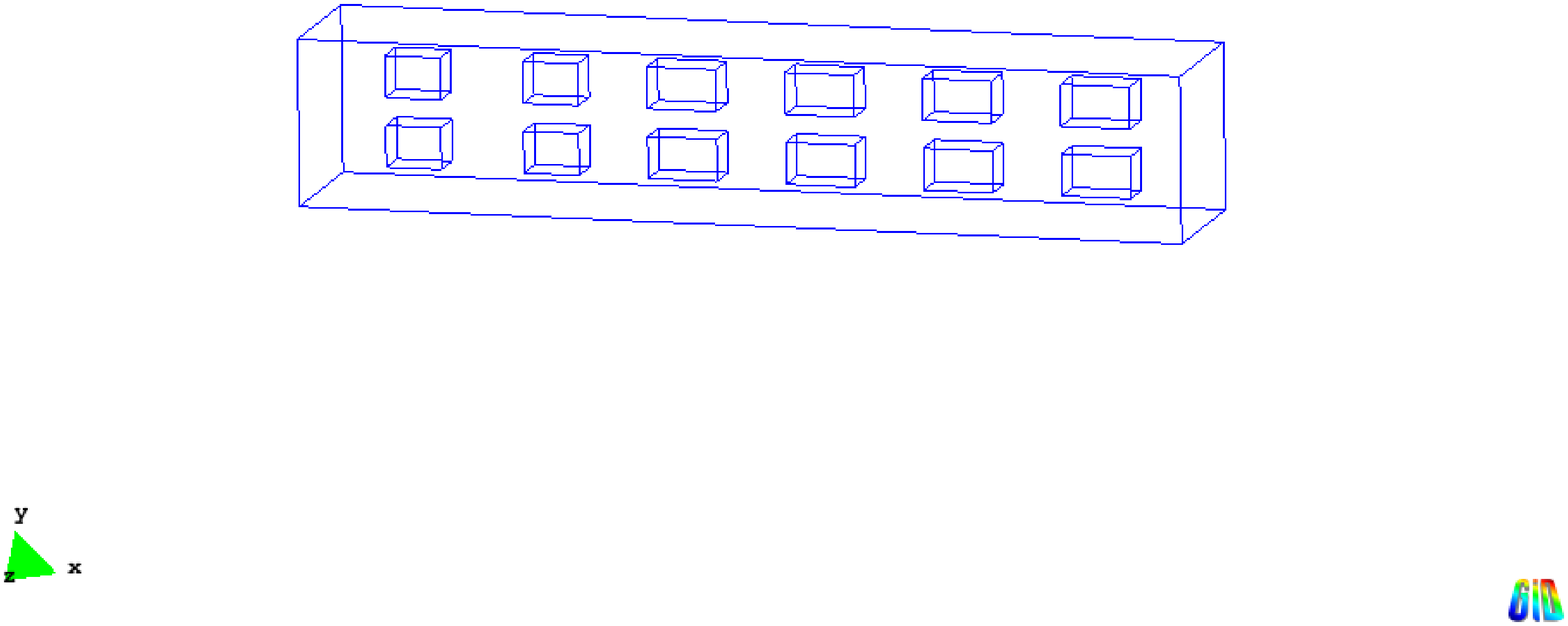}} \\
 a) Test1: $\Omega = \Omega_{FEM} \cup \Omega_{OUT}$ &  b) Test 1: $\Omega_{FEM}$ \\
 {\includegraphics[width=7.0cm, clip = true,  trim = 6.0cm 6.0cm 6.0cm 6.0cm]{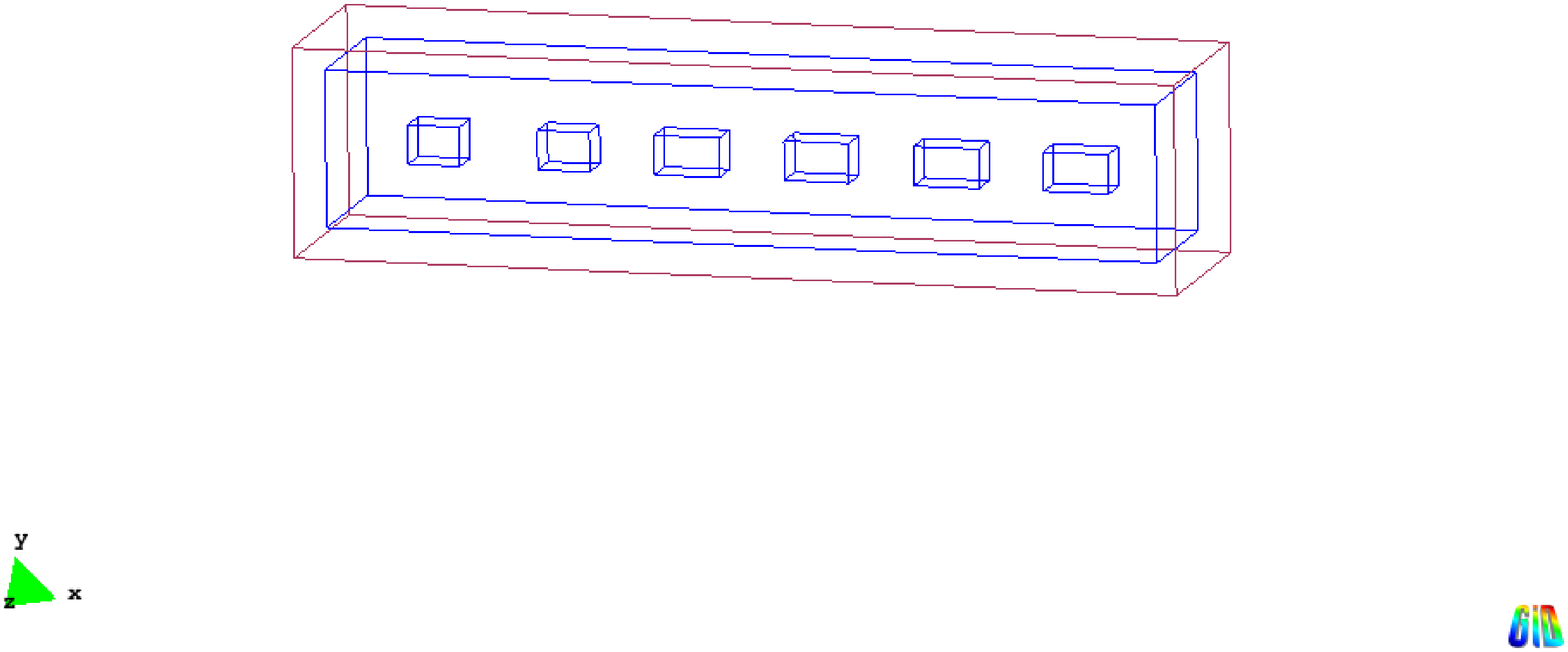}} &
{\includegraphics[width=7.0cm, clip = true, trim = 6.0cm 6.0cm 6.0cm 6.0cm]{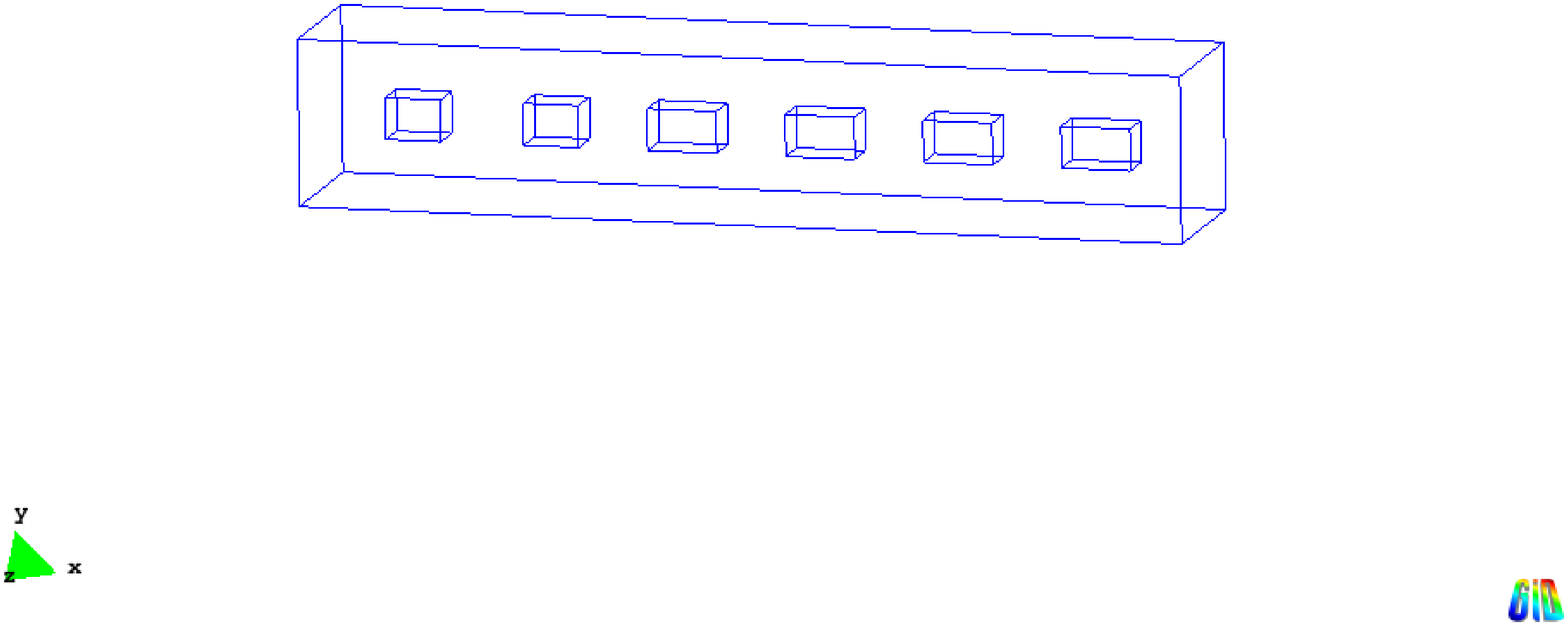}} \\
 c) Test 2: $\Omega = \Omega_{FEM} \cup \Omega_{OUT}$ &  d) Test 2:  $\Omega_{FEM}$ \\
 \end{tabular}
 \end{center}
 \caption{ \protect\small \emph{  Domain decomposition in numerical tests of Section \ref{sec:num}. a), c) The  decomposed domain  $\Omega= \Omega_{FEM} \cup \Omega_{OUT}$. b), d) The finite element  domain $\Omega_{FEM}$. }}
 \label{fig:fig1}
 \end{figure}

We divide
 $\Omega$ into two subregions, $\Omega_{\rm FEM}$ and $\Omega_{\rm OUT}$ such
 that $\overline{\Omega} = \overline{\Omega}_{\rm FEM} \cup \overline{ \Omega}_{\rm OUT}$,  $\Omega_{\rm FEM}
 \cap \Omega_{\rm OUT} = \emptyset$ and $\partial \Omega_{\rm FEM}
 \subset \partial \Omega_{\rm OUT}$. For an illustration of the domain
 decomposition, see Figure \ref{fig:fig1}. 
 The boundary $\partial \Omega$ is  such that $\partial \Omega
=\partial _{1} \Omega \cup \partial _{2} \Omega \cup \partial _{3} \Omega$  where
$\partial _{1} \Omega$ and $\partial _{2} \Omega$ are, respectively, front and
back sides of the domain $\Omega$, and $\partial _{3} \Omega$ is the union
of left, right, top and bottom faces of this domain.
 For numerical solution of
(\ref{model1_1})-(\ref{model1_4}) in $\Omega_{\rm OUT}$ we  can use  either the
finite difference or the finite element method on a structured mesh with constant
coefficients $\varepsilon = 1$ and $\mu =1$. In $\Omega_{\rm FEM}$, we use finite elements on a sequence of
unstructured meshes $K_h = \{K\}$, with elements $K$ consisting of
triangles in $\mathbb{R}^2$ and tetrahedra in $\mathbb{R}^3$
satisfying the maximal angle condition \cite{Brenner}. 
Let $S_T := \partial_1 \Omega \times (0,T)$ where
 $\partial_1 \Omega$ is the backscattering side of the domain $\Omega$ with
 the time domain observations, and define by $S_{1,1} := \partial_1
 \Omega \times (0,t_1]$, $S_{1,2} := \partial_1 \Omega \times
 (t_1,T)$,  $S_2 := \partial_2 \Omega \times (0, T)$,  $S_3 :=
 \partial_3 \Omega \times (0, T)$.

To simplify notations, further we will omit subscript $r$ in
$\varepsilon_r$ and $\mu_r$.
 We add a Coulomb-type gauge condition \cite{Ass,div_cor} to
(\ref{model1_1})-(\ref{model1_4}) for stabilization of the finite
element solution using the standard piecewise continuous functions
with $ 0 \leq s \leq 1$, and our model problem 
(\ref{model1_1})-(\ref{model1_4}) which we use in computations rewrites as
\begin{equation}\label{E_gauge}
\begin{split}
\varepsilon \frac{\partial^2 E}{\partial t^2} + \nabla \times ( \mu^{-1} \nabla \times E)  - s\nabla  ( \nabla \cdot(\varepsilon E))  &= 0~ \mbox{in}~~ \Omega_T, \\
  E(x,0) = f_0(x), ~~~E_t(x,0) &= f_1(x)~ \mbox{in}~~ \Omega,     \\
\partial _{n}E& = (0,f\left( t\right),0) ~\mbox{on}~ S_{1,1},
\\
\partial _{n} E& =-\partial _{t} E ~\mbox{on}~ S_{1,2},
\\
\partial _{n} E& =-\partial _{t} E ~\mbox{on}~ S_2, \\
\partial _{n} E& =0 ~\mbox{on}~ S_3, \\
\mu(x)=\varepsilon \left( x\right) &=1\text{ in }\Omega _{\rm OUT}. 
\end{split}
\end{equation}


In the recent works \cite{BMaxwell, BCN, BTKB} was demonstrated
numerically  that the solution of the problem (\ref{E_gauge})
approximates well the solution of the original Maxwell's system for
the case when $1 \leq \mu(x) \leq 2, 1 \leq \varepsilon(x) \leq 15 $
and $s=1$.


 We assume that our
 coefficients $\varepsilon \left(x\right), \mu(x) $ of equation (\ref{E_gauge})
 are such that
 
\begin{equation} \label{2.3}
\begin{split}
\varepsilon \left(x\right) &\in \left[ 1,d_1\right],~~ d_1 = const.>1,~ \varepsilon(x) =1
\text{ for }x\in  \Omega _{\rm OUT}, \\
\mu(x) &\in \left[ 1,d_2\right],~~ d_2 = const.>1,~ \mu(x) =1   \text{ for }x\in  \Omega _{\rm OUT}, \\
\varepsilon \left(x\right), \mu(x) &\in C^{2}\left( \mathbb{R}^{3}\right) . 
\end{split}
\end{equation}

In our numerical tests the values of constants $d_1, d_2$ in
(\ref{2.3}) are chosen from experimental set-up similarly with
\cite{BTKB, SSMS} and we assume that we know them a priori.

This is in agreement with the availability of a priori information
for an ill-posed problem \cite{BKS, Engl, tikhonov}.
Through the work we use  following  notations:  for any vector function $ u
\in \mathbb{R}^3$  when  we write    $u \in
H^k(\Omega), k=1,2$, we mean that every component of the vector function
$u$ belongs to this space. We consider the following

\textbf{Inverse Problem  (IP) } \emph{Assume that the functions
}$\varepsilon\left(x\right)$ and $\mu(x)$ \emph{\ satisfy conditions (\ref{2.3})  for the known   }$d_1, d_2 >1$\emph{\ and they are unknown in the
  domain }$\Omega \backslash  \Omega_{\rm OUT}$\emph{. Determine the functions }$ \varepsilon\left(x\right), \mu(x) $\emph{\ for }$x\in \Omega \backslash  \Omega_{\rm OUT},$ \emph{\ assuming that the
  following function }$\tilde E\left(x,t\right) $\emph{\ is known}
\begin{equation}
 E \left(x,t\right) = \tilde E \left(x,t\right) ~\forall \left( x,t\right) \in
  S_T.  \label{2.5}
\end{equation}


The function $\tilde E\left(x,t\right)$ in (\ref{2.5}) represents
the time-dependent measurements of the electric wave field $E(x,t)$ at the
backscattering boundary $\partial_1 \Omega$. In real-life experiments,
measurements are performed on a number of detectors, see details in
our recent experimental work \cite{BTKB}.

\section{Tikhonov functional}

\label{sec:tikhonov}

We reformulate our inverse problem as an optimization problem, where
 we seek for two functions,  the permittivity $\varepsilon(x)$ and permeability $\mu(x)$, which result in a solution of
equations (\ref{E_gauge}) with best fit to time and space domain observations
$\tilde E$, measured at a finite number of observation points on
$\partial_1 \Omega$. 
Our goal is to minimize  the Tikhonov functional
\begin{equation}
\begin{split}
J( \varepsilon, \mu) := J(E, \varepsilon, \mu) &= \frac{1}{2} \int_{S_T}(E - \tilde{E})^2 z_{\delta }(t) d \sigma dt \\
&+\frac{1}{2} \gamma_1  \int_{\Omega}( \varepsilon -  \varepsilon_0)^2~~ dx 
+ \frac{1}{2} \gamma_2  \int_{\Omega}( \mu -  \mu_0)^2~~ dx,
\label{functional}
\end{split}
\end{equation}
where $\tilde{E}$ is the observed electric field, $E$ satisfies the
equations (\ref{E_gauge}) and thus depends on $\varepsilon$ and $\mu$,
 $\varepsilon _{0}$ is the initial
guess for $\varepsilon $ and $\mu_{0}$ is the initial guess for $\mu$,
and $\gamma_i, i=1,2$ are the regularization parameters.  Here, $z_{\delta
}(t)$ is a cut-off function, which is introduced to ensure that
the compatibility conditions at $\overline{\Omega}_{T}\cap \left\{
t=T\right\} $ for the adjoint problem (\ref{adjoint}) are satisfied, and $\delta
>0$ is a small number. The function $z_{\delta }$ can be chosen as in \cite{BCN}.

Next, we introduce the following spaces of real valued vector functions
\begin{equation}\label{spaces}
\begin{split}
H_E^1 &:= \{ w \in H^1(\Omega_T):  w( \cdot , 0) = 0 \}, \\
H_{\lambda}^1 &:= \{ w \in  H^1(\Omega_T):  w( \cdot , T) = 0\},\\
U^{1} &=H_{E}^{1}(\Omega_T)\times H_{\lambda }^{1}(\Omega_T)\times C\left(\overline{\Omega}\right)\times C\left(\overline{\Omega}\right),\\
U^{0} &=L_{2}\left(\Omega_{T}\right) \times L_{2}\left(\Omega_{T}\right) \times
L_{2}\left(\Omega \right)\times
L_{2}\left(\Omega \right). 
\end{split}
\end{equation}

We also define the $L_2$  inner product and the norm  over $\Omega_T$  and $\Omega$ as
\begin{equation*}
\begin{split}
((u,v))_{\Omega_T}    &= \int_{\Omega} \int_0^T u v dx dt, \\
||u||^2  &= ((u,u))_{\Omega_T}, \\
(u,v)_{\Omega}    &= \int_{\Omega} u v dx, \\
|u|^2  &= (u,u)_{\Omega}.
\end{split}
\end{equation*}

To solve the minimization problem  we take into account  (\ref{2.3})  and introduce the Lagrangian
\begin{equation}\label{lagrangian}
\begin{split}
L(u) &= J(E, \varepsilon, \mu) 
- \left(\left( \varepsilon  \partial_t \lambda,  \partial_t E \right)\right)_{\Omega_T} 
-(\varepsilon \lambda(x,0), f_1(x))_{\Omega} 
+   \left( \left( \mu^{-1}\nabla \times E,  \nabla \times \lambda \right) \right)_{\Omega_T}   \\
&+ s  \left(\left( \nabla \cdot (\varepsilon E),  \nabla \cdot \lambda \right) \right)_{\Omega_T} - ((\lambda, p(t) ))_{S_{1,1}}   + (( \lambda \partial_t E ))_{S_{1,2}}  + (( \lambda \partial_t E ))_{S_2}, \\
\end{split}
\end{equation}
where $u=(E,\lambda, \varepsilon, \mu) \in U^1$ and $p(t)= (0,f(t), 0)$ and $\partial_t$ define the derivative in time.
We now  search for a stationary
point of the Lagrangian with respect to $u$ satisfying for all $\bar{u}= ( \bar{E}, \bar{\lambda},  \bar{\varepsilon}, \bar{\mu}) \in U^1$
\begin{equation}
 L'(u; \bar{u}) = 0 ,  \label{scalar_lagr2}
\end{equation}
where $ L^\prime (u;\cdot )$ is the Jacobian of $L$ at $u$. Equation
above can be written as
\begin{equation}
 L'(u; \bar{u}) = \frac{\partial L}{\partial \lambda}(u)(\bar{\lambda}) + \frac{\partial L}{\partial E}(u)(\bar{E})  +  \frac{\partial L}{\partial \varepsilon}(u)(\bar{\varepsilon}) + \frac{\partial L}{\partial \mu}(u)(\bar{\mu}) = 0.  \label{scalar_lagr}
\end{equation}
 To find the Frech\'{e}t derivative (\ref{scalar_lagr}) of the
 Lagrangian (\ref{lagrangian}) we consider $L(u + \bar{u}) - L(u)$
 for all $\bar{u} \in U^1$ and single out the linear part of the
 obtained expression with respect to $ \bar{u}$.  In our derivation of
 the Frech\'{e}t derivative we assume that in the Lagrangian
 (\ref{lagrangian}) functions $u=(E,\lambda, \varepsilon, \mu) \in U^1$ can  vary independently of each other. In this approach we
 obtain the same result as by assuming that functions $E$ and
 $\lambda$ are dependent on the coefficients $\varepsilon, \mu$, see
 also Chapter 4 of \cite{BOOK} where similar observations take place.
 Taking into
 account that $E(x,t)$ is the solution of the forward problem (\ref{E_gauge}), assumptions  that $\lambda (x ,T) = \frac{\partial
   \lambda}{\partial t} (x,T) =0$, as well as $\mu=\varepsilon=1$ on
 $\partial \Omega$  and using conditions (\ref{2.3}), we obtain from
 (\ref{scalar_lagr}) that for all $\bar{u}$,
\begin{equation}\label{forward}
\begin{split}
0 = \frac{\partial L}{\partial \lambda}(u)(\bar{\lambda}) =
&- (( \varepsilon \partial_t \bar{\lambda},  \partial_t E ))_{\Omega_T}  -(\varepsilon f_1(x), \bar{\lambda}(x,0))_{\Omega}
+    (( \mu^{-1} \nabla \times E, \nabla \times \bar{\lambda}))_{\Omega_T}   \\
&+ s ((\nabla \cdot(\varepsilon E),  \nabla \cdot \bar{\lambda}))_{\Omega_T} 
- ((\bar{\lambda}, p(t)))_{S_{1,1}}   +  (( \bar{\lambda},  \partial_t E ))_{S_{1,2}}   \\
& + (( \bar{\lambda},  \partial_t E ))_{S_2}~~\forall \bar{\lambda} \in H_{\lambda}^1(\Omega_T),
\end{split}
\end{equation}
\begin{equation} \label{control}
\begin{split}
0 = \frac{\partial L}{\partial E}(u)(\bar{E}) &=
 ((E - \tilde{E},  \bar{E} z_{\delta} ))_{S_T} 
- ((\varepsilon  \partial_t \lambda,  \partial_t \bar{E}))_{\Omega_T}
 +  ((\mu^{-1}  \nabla \times \lambda, \nabla \times \bar{E} ))_{\Omega_T} \\
 &+ s  ((\nabla \cdot \lambda,  \nabla \cdot (\varepsilon \bar{E}))_{\Omega_T}
- (( \partial_t \lambda,  \bar{E} ))_{S_{1,2} \cup S_2}
  -(\varepsilon \bar{E}(x,0), \partial_t \lambda(x,0) )~~\forall \bar{E} \in H_{E}^1(\Omega_T).
\end{split}
\end{equation}
Further, we obtain two equations that express that the gradients
with respect to  $\varepsilon$ and $\mu$ vanish:
\begin{equation} \label{grad1} 
\begin{split}
0 = \frac{\partial L}{\partial \varepsilon}(u)(\bar{\varepsilon})
 =  &- (( \partial_t \lambda, \partial_t E~ \bar{\varepsilon} ))_{\Omega_T} 
- (\lambda(x,0), f_1(x)~\bar{\varepsilon})_{\Omega} \\
&+  s((\nabla \cdot (\bar{\varepsilon} E), \nabla \cdot \lambda ))_{\Omega_T} 
 +\gamma_1 (\varepsilon - \varepsilon_0,  \bar{\varepsilon})_{\Omega}~~\forall  x \in \Omega,  
\end{split}
\end{equation}
\begin{equation} \label{grad2} 
0 = \frac{\partial L}{\partial \mu}(u)(\bar{\mu})
 =  -(( \mu^{-2}~\nabla \times E,  \nabla \times \lambda ~ \bar{\mu}))_{\Omega_T}
 +\gamma_2 (\mu - \mu_0, \bar{\mu})_{\Omega} ~\forall x \in \Omega.
\end{equation}

We observe that the equation (\ref{forward}) is the weak formulation
of the state equation (\ref{E_gauge}) and the equation (\ref{control})
is the weak formulation of the following adjoint problem
\begin{equation}
\begin{split} \label{adjoint}
\varepsilon \frac{\partial^2 \lambda}{\partial t^2} + 
  \nabla \times (\mu^{-1} \nabla \times \lambda) -  s \varepsilon \nabla  (
 \nabla \cdot \lambda) &= -  (E - \tilde{E})|_{S_T} z_{\delta} ~ \mbox{  in } \Omega_T,   \\
\lambda(\cdot, T)& =  \frac{\partial \lambda}{\partial t}(\cdot, T) = 0, \\
\partial _{n} \lambda& = \partial _{t} \lambda ~\mbox{on}~ S_{1,2},
\\
\partial _{n} \lambda& = \partial _{t} \lambda ~\mbox{on}~ S_2, \\
\partial _{n} \lambda& =0 ~\mbox{on}~ S_3,
\end{split}
\end{equation}
which is solved backward in time.

We now define by $E(\varepsilon, \mu), \lambda(\varepsilon, \mu)$ the exact solutions of the
forward and adjoint problems for given $\varepsilon, \mu$, respectively. Then defining by 
\begin{equation*}
u(\varepsilon, \mu) = (E(\varepsilon, \mu), \lambda(\varepsilon, \mu), \varepsilon, \mu) \in U^1,
\end{equation*}
 using the fact that  for exact solutions  $E(\varepsilon, \mu), \lambda(\varepsilon, \mu)$ because of 
 (\ref{lagrangian}) we have
\begin{equation}
J( E(\varepsilon, \mu), \varepsilon, \mu) = L(u(\varepsilon, \mu)),
\end{equation}
~ and assuming that solutions $E(\varepsilon, \mu), \lambda(\varepsilon, \mu) $  are sufficiently stable, see Chapter 5 of book \cite{lad} for details,   we can write that the  Frech\'{e}t derivative of the Tikhonov functional is the function $J'(\varepsilon, \mu, E(\varepsilon, \mu))$ which is defined as
\begin{equation}\label{derfunc}
\begin{split}
J'(\varepsilon, \mu) := J'(\varepsilon, \mu, E(\varepsilon, \mu)) &=  \frac{\partial J}{\partial \varepsilon}(\varepsilon, \mu, E(\varepsilon, \mu) ) + \frac{\partial J}{\partial \mu}(\varepsilon, \mu, E(\varepsilon, \mu))   \\
&=  \frac{\partial L}{\partial \varepsilon}(u(\varepsilon, \mu)) + \frac{\partial L}{\partial \mu}(u(\varepsilon, \mu)).
\end{split}
\end{equation}
Inserting (\ref{grad1}) and (\ref{grad2}) into  (\ref{derfunc}),
we get
\begin{equation} \label{derfunc2}
\begin{split}
J'(\varepsilon, \mu)(x) &:= J'(\varepsilon, \mu, E(\varepsilon, \mu))(x) =
   - \int_0^T \partial_t \lambda~ \partial_t E~ (x,t)~ dt - \lambda(x,0) f_1(x) \\
&+  s  \int_0^T (\nabla \cdot E )  (\nabla \cdot \lambda )  ~ (x,t)~ dt 
 +\gamma_1 (\varepsilon - \varepsilon_0)(x) \\
 &- \int_0^T (\mu^{-2}~\nabla \times E) (  \nabla \times \lambda) ~ (x,t)~ dt  + \gamma_2 (\mu - \mu_0)(x).
\end{split}
\end{equation}

\section{Finite element method}
\label{sec:spaces}

\subsection{Finite element spaces}

For computations we discretize 
$\Omega_{\rm FEM} \times (0,T)$ in space and time. For discretization in space we
denote by $K_h = \{K\}$ a partition of the domain $\Omega_{\rm FEM}$
into tetrahedra $K$ in $ \mathbb{R}^{3}$ or triangles in $
\mathbb{R}^{2}$.  We discretize the time interval $(0,T)$ into
subintervals $J=(t_{k-1},t_k]$ of uniform length $\tau = t_k -
  t_{k-1}$ and denote the time partition by $J_{\tau} = \{J\}$.  In
  our finite element space mesh $K_h$ the elements $K$ are
  such that
\begin{equation*}
 K_h = \cup_{K \in K_h} K=K_1 \cup K_2...\cup K_l,
\end{equation*}
 where $l$ is the total
 number of elements $K$ in $\overline{\Omega}$.

Similarly with \cite{EEJ} we introduce
the mesh function $h=h(x)$ which is a piecewise-constant function such
that
\begin{equation}\label{meshfunction}
h |_K = h_K ~~~ \forall K \in K_h,
\end{equation}
where $h_K$ is the diameter of $K$ which we define as the longest
side of $K$. 
Let $r^{\prime }$ be the radius of the maximal circle/sphere contained
in the element $K$.  For every element $K \in K_h$
we assume the following shape regularity assumption
\begin{equation}
a_{1}\leq h_K \leq r^{\prime }a_{2};\quad
a_{1},a_{2}=const.>0.  \label{2.1}
\end{equation}

To formulate the finite element method  for (\ref{scalar_lagr}),  we
 define the finite element spaces.
First we introduce the finite element trial space $W_h^E$ for every component of the electric field $E$ defined by
\begin{equation}
W_h^E := \{ w \in H_E^1: w|_{K \times J} \in  P_1(K) \times P_1(J),  \forall K \in K_h,  \forall J \in J_{\tau} \}, \nonumber
\end{equation}
where $P_1(K)$ and $P_1(J)$ denote the set of linear functions on $K$
and $J$, respectively.
We also introduce the finite element test space  $W_h^{\lambda}$ defined by
\begin{equation}
W_h^{\lambda} := \{ w \in H_{\lambda}^1: w|_{K \times J} \in  P_1(K) \times P_1(J),  \forall K \in K_h,  \forall J \in J_{\tau} \}. \nonumber
\end{equation}


To approximate the functions 
$\mu(x)$  and $\varepsilon(x)$  we will use the space of piecewise constant functions $V_{h} \subset L_2(\Omega)$, 
\begin{equation*}
V_{h}:=\{u\in L_{2}(\Omega ):u|_{K}\in P_{0}(K),\forall K\in  K_h\}, 
\end{equation*}
where $P_{0}(K)$ is the space of piecewise constant functions on $K$.
In some numerical experiments we will use also  the space of piecewise linear functions $W_{h} \subset H^1(\Omega)$, 
\begin{equation}
W_h = \big\{ v(x) \in H^1(\Omega):~ v|_{K} \in P_1(K)~\forall K \in K_h \big\}.
\end{equation}

  In a general case we allow  the functions $\varepsilon(x),
  \mu(x)$ to be  discontinuous, see \cite{KN}.
Let $F_h$ be the set of all faces of elements in $K_h$ such that $F_h
:= F_h^I \cup F_h^B$, where $F_h^I$ is the set of all interior faces
and $F_h^B$ is the set of all boundary faces of elements in $K_h$.
Let $f \in F_h^I$ be the internal face of the non-empty intersection of the boundaries
of two neighboring elements $K^{+}$ and $K^{-}$.
We denote the jump of the function $v_{h}$ computed from the two
neighboring elements $K^{+}$ and $K^{-}$ sharing the common side $f$ as
\begin{equation}
[ v_{h}]= v_{h}^{+}- v_{h}^{-},
\label{2.4}
\end{equation}
and the jump of the normal component $v_h$ across the side $f$ as
\begin{equation}
[[ v_{h} ]] = v_{h}^{+} \cdot n^+ +  v_{h}^{-} \cdot n^-,
\label{jump_normal}
\end{equation}
where $n^+, n^-$ is the unit outward normal on $f^+, f^-$, respectively.

Let $P_h$
 be  the
 $L_2(\Omega)$ orthogonal projection.  We define by $f_h^I$
the standard nodal interpolant \cite{EEJ} of $f$ into the space of continuous
piecewise-linear functions on the mesh $K_h$.  Then by one of properties of the
orthogonal projection
\begin{equation}
 \left\| f- P_hf\right\| _{L_{2}\left( \Omega \right) }\leq \left\| f-f_{h}^I\right\| _{L_{2}\left( \Omega \right) }. \label{2.5b}
\end{equation}
It follows
from  \cite{SZ} that
\begin{equation}
\left\| f-P_hf\right\| _{L_{2}\left( \Omega \right) }\leq  C_I h \left\|~  f\right\| _{H^1\left( \Omega \right) }~~\forall f\in H^1(\Omega).
\label{2.6}
\end{equation}
where $ C_I = C_I\left( \Omega \right) $ is a positive
constant depending only on the domain $\Omega$.

\subsection{A finite element method for  optimization problem}
\label{sec:fem}

To formulate the finite element method for (\ref{scalar_lagr}) we
define the space $U_h = W_h^E \times W_h^{\lambda} \times V_h \times
V_h$.
The finite
element method reads: Find $u_h \in U_h$ such that
\begin{equation} \label{varlagr}
L'(u_h)(\bar{u})=0 ~\forall
\bar{u} \in U_h . 
\end{equation}

To be more precise, the equation (\ref{varlagr}) expresses that the
finite element method for the forward problem (\ref{E_gauge}) in
$\Omega_{FEM}$ for continuous $(\varepsilon, \mu)$ will be: find $E_h =({E_1}_h,{E_2}_h, {E_3}_h ) \in
W_h^E$, such that for all $\bar{\lambda} \in W_h^\lambda$ and for the
known $(\varepsilon_h, \mu_h) \in V_h \times V_h$
\begin{equation}\label{varforward}
\begin{split}
&- ((\varepsilon_h \frac{\partial \bar{\lambda}}{\partial t} \frac{\partial E_h}{\partial t})) 
-   (\varepsilon_h f_1,\bar{\lambda}(x,0) )_{\Omega}  +   (( \mu_h^{-1} \nabla \times E_h, \nabla \times \bar{\lambda}))_{\Omega_T}   
+ s ((\nabla \cdot(\varepsilon_h E_h),  \nabla \cdot \bar{\lambda}))_{\Omega_T} \\
&- ((\bar{\lambda}, p(t)))_{S_{1,1}}   + (( \bar{\lambda},  \partial_t E_h ))_{S_{1,2}}   
 + (( \bar{\lambda},  \partial_t E_h ))_{S_2}=0~~\forall \bar{\lambda} \in H_{\lambda}^1(\Omega_T),
\end{split}
\end{equation}
and the finite element method for the adjoint problem (\ref{adjoint})
in $\Omega_{FEM}$ for continuous $(\varepsilon, \mu)$ reads: find $\lambda_h = ({\lambda_h}_1,{\lambda_h}_2, {\lambda_h}_3)  \in W_h^\lambda$, such that
for the computed approximation $E_h =({E_1}_h,{E_2}_h, {E_3}_h ) \in
W_h^E$ of (\ref{varforward}) and for all $ \bar{E} \in W_h^E$ and for the
known $(\varepsilon_h, \mu_h) \in V_h \times V_h$
\begin{equation} \label{varadjoint}
\begin{split}
& (((E_h - \tilde{E})|_{S_T} z_{\delta},  \bar{E} )) 
- ((\varepsilon_h  \partial_t \lambda_h,  \partial_t \bar{E}))_{\Omega_T}
 +  ((\mu_h^{-1}  \nabla \times \lambda_h, \nabla \times \bar{E} ))_{\Omega_T} \\
 &+ s  ((\nabla \cdot \lambda_h,  \nabla \cdot (\varepsilon_h \bar{E}))_{\Omega_T} - (( \partial_t \lambda_h,  \bar{E} ))_{S_{1,2} \cup S_2}  -( \varepsilon_h \bar{E}(x,0), \partial_t \lambda_h(x,0) ) = 0~~\forall \bar{E} \in H_{E}^1(\Omega_T).\\
\end{split}
\end{equation}

Similar finite element method for the forward and adjoint problems can
be written for discontinuous functions $\varepsilon, \mu$ which will
include additional terms with jumps for computation of coefficients.
In our work similarly with \cite{CWZ14} we compute the discontinuities
of coefficients $\epsilon$ and $\mu$ by computing the jumps from the
two neighboring elements, see (\ref{2.4}) and (\ref{jump_normal}) for
definitions of jumps.

Since we are usually working in  finite dimensional spaces $U_{h}$
and $U_{h} \subset U^{1}$ as a set, then $U_{h}$ is a discrete
analogue of the space $U^{1}.$
It is well known that in finite dimensional spaces all norms are
equivalent, and in our computations we compute approximations of
smooth functions $\varepsilon(x), \mu(x)$ in the space $V_h$.

\subsection{Fully discrete scheme}
\label{sec:discrete}

To write fully discrete schemes for (\ref{varforward}) and (\ref{varadjoint}) 
we expand $E_h$ and $\lambda_h$ in terms of the standard continuous piecewise
linear functions $\{\varphi_i(x)\}_{i=1}^M$ in space and
$\{\psi_k(t)\}_{k=1}^N$ in time, respectively, as 
\begin{equation*}
\begin{split}
E_h (x,t) &=\sum_{k=1}^N \sum_{i=1}^M \mathbf{E_{h}}  \varphi_i(x)\psi_k(t), \\
\lambda_h(x,t) &=\sum_{k=1}^N \sum_{i=1}^M \mathbf{\lambda_{h}} \varphi_i(x)\psi_k(t),
\end{split}
\end{equation*}
 where $\mathbf{E_h} := E_{h_{i,k}}$ and $\mathbf{\lambda_h} := \lambda_{h_{i,k}}$ denote unknown
coefficients at the point $x_i \in K_h$ and time level $t_k \in J_{\tau}$,
substitute them into (\ref{varforward}) and  (\ref{varadjoint}) to obtain the
following system of linear equations:
\begin{equation} \label{femod1}
\begin{split}
 M (\mathbf{E}^{k+1} - 2 \mathbf{E}^k  + \mathbf{E}^{k-1})  &= 
  - \tau^2  K \mathbf{E}^k - s \tau^2  C \mathbf{E}^k + \tau^2 F^k + \tau^2 P^k - \frac{1}{2}\tau (MD)\cdot(\mathbf{E}^{k+1} - \mathbf{E}^{k-1}),   \\
M (\boldsymbol{\lambda}^{k+1} - 2 \boldsymbol{\lambda}^k + \boldsymbol{
\lambda}^{k-1}) &=  -\tau^2  S^k - \tau^2  K \boldsymbol{\lambda}^k - s \tau^2  C \boldsymbol{\lambda}^k +
  \frac{1}{2} \tau (MD)\cdot(\boldsymbol{\lambda}^{k+1} - \boldsymbol{\lambda}^{k-1}) + \tau^2 (D\lambda)^k.\\
\end{split}
\end{equation}
  Here, $M$ is the block mass matrix in space and $MD$ is the block
  mass matrix in space corresponding to the elements at the boundaries
  $\partial_1 \Omega, \partial_2 \Omega$, $K$ is the block stiffness
  matrix corresponding to the rotation term, $C$ is the stiffness
  matrix corresponding to the divergence term, $ F^k, P^k, D\lambda^k,
  S^k$ are load vectors at time level $t_k$, $\mathbf{E}^k$ and $
  \boldsymbol{\lambda}^k$ denote the nodal values of $\mathbf{E_h}$
  and $\mathbf{\lambda_h}$, respectively, at time level $t_k$, and
  $\tau$ is the time step. We refer to \cite{BMaxwell} for details of derivation of
  these matrices.

Let us define the mapping $F_K$ for the reference element $\hat{K}$
such that $F_K(\hat{K})=K$ and let $\hat{\varphi}$ be the piecewise
linear local basis function on the reference element $\hat{K}$ such
that $\varphi \circ F_K = \hat{\varphi}$.  Then the explicit formulas
for the entries in system (\ref{femod1}) at each element $K$ can be
given as:
\begin{equation*}
\begin{split}
  M_{i,j}^{K} & =    (\varepsilon_h ~\varphi_i \circ F_K, \varphi_j \circ F_K)_K, \\
  K_{i,j}^{K} & =   ( \mu_h^{-1} \nabla \times \varphi_i \circ F_K, \nabla \times \varphi_j \circ F_K)_K,\\
  C_{i,j}^{K} & =   ( \nabla\cdot (\varepsilon_h \varphi_i) \circ F_K, \nabla \cdot \varphi_j \circ F_K)_K,\\
  S_{j}^{K}&= ((E_h -\tilde{E})_{S_T} z_{\sigma}, \varphi_j \circ F_K )_{K}, \\
 F_{j}^{K}&= (\varepsilon_h f_1, \varphi_j \circ F_K )_{K}, \\
 P_{j}^{K}&= ( p, \varphi_j \circ F_K )_{\partial_1 \Omega_K}, \\
 MD_{j}^{K}&= (~\varphi_i \circ F_K, \varphi_j \circ F_K )_{\partial_1 \Omega_K \cup \partial_2 \Omega_K}, \\
 D\lambda_{j}^{K}&= ( \varepsilon_h \partial_t \lambda_h(x,0), \varphi_j \circ F_K )_K, \\
\end{split}
 \end{equation*}
where $(\cdot,\cdot)_K$ denotes the $L_2(K)$ scalar product, and $\partial_1 \Omega_K, \partial_2 \Omega_K$
are  boundaries $\partial K$ of elements $K$, which belong to  $\partial_1 \Omega, \partial_2 \Omega$, respectively.

To obtain an explicit scheme, we approximate $M$ with the lumped mass matrix
$M^{L}$ (for further details, see \cite{Cohen}). Next, we multiply (\ref{femod1}) with
$(M^{L})^{-1}$ and get the following  explicit method:
\begin{equation}   \label{fem_maxwell_full}
\begin{split}
(I + \frac{1}{2} \tau (M^{L})^{-1} MD) \mathbf{E}^{k+1} = &2\mathbf{E}^k
  - \tau^2  (M^{L})^{-1} K\mathbf{E}^k
 +\tau^2 (M^{L})^{-1} F^k + \tau^2 (M^{L})^{-1} P^k\\
&+ \frac{1}{2} \tau (M^{L})^{-1} (MD) \mathbf{E}^{k-1}  
- s \tau^2  (M^{L})^{-1} C \mathbf{E}^k   -\mathbf{E}^{k-1},  \\
(I + \frac{1}{2} \tau (M^L)^{-1} MD) \boldsymbol{\lambda}^{k-1} = &-\tau^2 (M^{L})^{-1} S^k  
  + 2\boldsymbol{\lambda}^k
  - \tau^2  (M^{L})^{-1} K \boldsymbol{\lambda}^k - s \tau^2  (M^{L})^{-1} C \boldsymbol{\lambda}^k  \\
&+ \tau^2 (M^{L})^{-1} (D\lambda)^k  -\boldsymbol{\lambda}^{k+1}
+ \frac{1}{2} \tau (M^L)^{-1} (MD) \boldsymbol{\lambda}^{k+1}. 
\end{split}
\end{equation} 
In the case of the domain decomposition FEM/FDM method   when  the schemes above are used only in $\Omega_{FEM}$   we have
\begin{equation}   \label{fem_maxwell}
\begin{split}
 \mathbf{E}^{k+1} = &2\mathbf{E}^k
  - \tau^2  (M^{L})^{-1} K\mathbf{E}^k
 +\tau^2 (M^{L})^{-1} F^k + \tau^2 (M^{L})^{-1} P^k
- s \tau^2  (M^{L})^{-1} C \mathbf{E}^k   -\mathbf{E}^{k-1},  \\
 \boldsymbol{\lambda}^{k-1} = &-\tau^2 (M^{L})^{-1} S^k  
  + 2\boldsymbol{\lambda}^k
  - \tau^2  (M^{L})^{-1} K \boldsymbol{\lambda}^k - s \tau^2  (M^{L})^{-1} C \boldsymbol{\lambda}^k  
+ \tau^2 (M^{L})^{-1} D\lambda  -\boldsymbol{\lambda}^{k+1}. 
\end{split}
\end{equation}


\section{Relaxation property of mesh refinements}

\label{sec:relax}

In this section we reformulate results of \cite{BKK} for the case  of our \textbf{IP}.
For simplicity, we shall sometimes write 
 $||\cdot||$  for the $L_2$ norm.

 We use the theory of ill-posed problems \cite{tikhonov, TGSK}
 and introduce noise level $\delta $ in the function $\tilde{E}(x,t)$ in
 the Tikhonov functional (\ref{functional}). This means that
\begin{equation}
\tilde{E}(x,t)= \tilde{E}^{\ast }(x,t) + \tilde{E}_{\delta
}(x,t);\text{ } \tilde{E}^{\ast }, \tilde{E}_{\delta }\in L_{2}\left(
S_{T}\right) =H_{2}, \label{4.247}
\end{equation}
where $\tilde{E}^{\ast }(x,t)$ is the exact data corresponding to the exact function $z^*=(\varepsilon^*, \mu^*)$, and the function $\tilde{E}_{\delta }(x,t)$
represents the error in these data. In other words, we can write that
\begin{equation}
\left\Vert \tilde{E}_{\delta }\right\Vert _{L_{2}\left( S_{T}\right) }\leq \delta .
\label{4.248}
\end{equation}

The question of stability and uniqueness of our \textbf{IP} is
addressed in \cite{BCN, BCS} which is needed in the local strong
convexity theorem formulated below.  Let $H_1$ be the finite dimensional
linear space.
 Let $Y$ be  the set of admissible functions $(\varepsilon, \mu)$ which we defined in (\ref{2.3}), and let
  $Y_1 := Y \cap H_1$ with $G := \bar{Y}_1$.
We introduce now the operator $F:
G \to H_2$ corresponding to the Tikhonov functional (\ref{functional})
such that
\begin{equation}\label{opF}
F(z)(x,t) := F(\varepsilon, \mu)(x,t) =  (E(x,t,\varepsilon, \mu) - \tilde{E})^2 z_{\delta }(t)~~ \forall (x,t) \in S_T,
\end{equation} 
where $E(x,t,\varepsilon, \mu) := E(x,t)$ is the weak solution of the
forward problem (\ref{E_gauge}) and thus, depends on $\varepsilon$ and
$\mu$. Here, $z=(\varepsilon, \mu)$ and $z_{\delta }(t)$ is a cut-off
function chosen as in \cite{BCN}.

 We now assume that the operator $F(z)(x,t) $  which we defined in (\ref{opF}) is one-to-one. 
Let us denote by 
\begin{equation}\label{neigh}
V_{d}\left( z\right) =\left\{z^{\prime }\in H_{1}:\left\| z^{\prime } - z \right\|   < d ~~\forall z=(\varepsilon, \mu)  \in H_{1}\right\} 
\end{equation}
the neighborhood of $z$ of the diameter $d$. 
 We also
assume  that the operator $F$ is  Lipschitz continuous what means that for  $N_{1},N_{2} > 0$
\begin{equation}
\left\| F^{\prime }(z) \right\| \leq N_{1},\left\| F^{\prime
}(z_1)  - F^{\prime }(z_2) \right\| \leq N_{2}\left\|
 z_1 - z_2 \right\|~~\forall z_1, z_2 \in V_{1}\left( z^{\ast }\right) .  \label{2.7}
\end{equation}

Let the constant $D= D\left( N_{1},N_{2}\right) =const.>0$ is such that 
\begin{equation}
\left\| J^{\prime }\left( z_1\right) -J^{\prime }\left(z_2\right) \right\| \leq D\left\| z_1 - z_2\right\|~~\forall z_1, z_2 \in V_{1}(z^*),
\label{2.10}
\end{equation}
where $(\varepsilon^*, \mu^*)$ is the exact solution of the equation $F(\varepsilon^*, \mu^*)=0$.
Similarly with   \cite{BKK},  we assume that
\begin{equation}\label{2.11a}
\begin{split}
\left\| \varepsilon_0 - \varepsilon^{\ast }\right\| &\leq\delta ^{\nu _{1}}, ~\nu_1 =const.\in \left(0,1\right),   \\
\left\| \mu_0 - \mu^{\ast }\right\| &\leq \delta ^{\nu_2}, ~\nu_2 =const.\in \left(0,1\right), \\
\gamma_1 &= \delta ^{\zeta_1}, \zeta_1= const.\in ( 0,\min (\nu_1,2(1- \nu_1)), \\
\gamma_2 &= \delta ^{\zeta_2}, \zeta_2= const.\in ( 0,\min (\nu_2,2(1- \nu_2)),
\end{split}
\end{equation}
which in closed form can be written as
\begin{eqnarray}
\left\| z_0 - z^{\ast }\right\| &\leq &\delta ^{(\nu _{1}, \nu_2)},
~z_0=(\varepsilon_0, \mu_0), ~(\nu_1, \nu_2) =const.\in \left( 0,1\right) ,  \label{2.11} \\
(\gamma_1, \gamma_2) &=&\delta ^{(\zeta_1, \zeta_2)}, 
(\zeta_1, \zeta_2) =const.\in \left( 0,\min \left((\nu_1, \nu_2), 2\left( 1- (\nu_1, \nu_2)\right) \right) \right), \label{2.12}
\end{eqnarray}
where $(\gamma_1, \gamma_2)$ are regularization parameters in
(\ref{functional}).  Equation (\ref{2.11}) means that we assume that
all initial guesses $z_0=(\varepsilon_0, \mu_0)$ are located in a
sufficiently small neighborhood $V_{\delta ^{\mu _{1}}}(z^*) $ of the
exact solution $z^*=(\varepsilon^*, \mu^*)$.  Conditions (\ref{2.12})
imply that $(z^{\ast }, z_0)$ belong to an appropriate neighborhood of
the regularized solution of the functional (\ref{functional}), see
proofs in Lemmata 2.1 and 3.2 of \cite{BKK}.

Below we reformulate Theorem 1.9.1.2 of \cite{BOOK} for the Tikhonov
functional (\ref{functional}). Different proofs of it can be found in
\cite{BOOK} and in \cite{BKK} and are straightly applied to our  \textbf{IP}.
We note here that if functions
$(\varepsilon, \mu) \in H_{1}$ and satisfy conditions (\ref{2.3}) then
$(\varepsilon, \mu) \in \rm{Int}\left( G\right) .$

\textbf{Theorem 1} \emph{Let }$\Omega \subset \mathbb{R}^{3}$
\emph{\ be a convex bounded domain with the boundary }$\partial \Omega
\in C^{3}.$ \emph{  Suppose that
conditions (\ref{4.247}) and (\ref{4.248}) hold.  Let the function
}$E(x,t) \in H^{2}(\Omega_T)$ \emph{\ in
 the Tikhonov functional  (\ref{functional}) be the solution of the forward problem (\ref{E_gauge}) for
  the functions }$(\varepsilon, \mu) \in G$. \emph{ Assume that there exists
    the exact solutions }$(\varepsilon^*, \mu^*) \in G$\emph{\ of the
        equation }$F(\varepsilon^*, \mu^*) =0$ \emph{\ for the case of
        the exact data }$ \tilde{E}^{\ast }$\emph{\ in (\ref{4.247}).
      Let regularization parameters } $(\gamma_1, \gamma_2)$ \emph{ in (\ref{functional}) are such
      that }
\begin{equation*}
(\gamma_1, \gamma_2) = (\gamma_1, \gamma_2) \left( \delta \right) =\delta ^{2(\nu_1, \nu_2) },~~(\nu_1, \nu_2)  =const.\in \left( 0,\frac{1}{4}\right)~~\quad \forall \delta
\in \left( 0,1\right).
\end{equation*}
\emph{Let } $z_0=(\varepsilon_0, \mu_0)$ \emph{ satisfy
  (\ref{2.11}). Then the Tikhonov functional  (\ref{functional})
    is strongly convex in the neighborhood }$V_{(\gamma_1, \gamma_2) \left( \delta
    \right) }\left(\varepsilon^*, \mu^* \right) $ \emph{\ with the strong convexity constants }$(\alpha_1, \alpha_2)=(\gamma_1, \gamma_2)
  /2.$ \emph{The strong convexity property can be also written as} 
\begin{equation}
\left\Vert z_{1} - z _{2}\right\Vert ^{2}\leq \frac{2}{\delta ^{2(\nu_1, \nu_2) }}\left(
J'(z_1) - J'(z_2), z_{1} - z_{2}\right)~~\forall z_1 =(\varepsilon_1, \mu_1), z_2 =(\varepsilon_2, \mu_2) \in H_{1}.  \label{4.249}
\end{equation}
\emph{ Alternatively,  using the expression for the Fr\'{e}chet
  derivative given in (\ref{derfunc}) we can write (\ref{4.249}) as}
\begin{equation}\label{convex}
\begin{split}
\left\Vert \varepsilon_{1} - \varepsilon_{2} \right \Vert ^{2}  &\leq \frac{2}{\delta ^{2\nu_1 }}\left(
J_{\varepsilon}'(\varepsilon_1, \mu_1) - J_{\epsilon}'(\varepsilon_2, \mu_2), \varepsilon_{1} - \varepsilon_{2}\right)~~\forall (\varepsilon_1, \mu_1),  (\varepsilon_2, \mu_2)  \in H_{1},  \\
\left\Vert \mu_{1} - \mu_{2} \right \Vert ^{2}  &\leq \frac{2}{\delta ^{2\nu_2 }}\left(
J_{\mu}'(\varepsilon_1, \mu_1) - J_{\mu}'(\varepsilon_2, \mu_2), \mu_{1} - \mu_{2}\right)~~\forall (\varepsilon_1, \mu_1), (\varepsilon_2, \mu_2)  \in H_{1},
\end{split}
\end{equation}
\emph{where }$\left(\cdot , \cdot \right) $\emph{\ is the} $ L_2(\Omega)$
\emph{ inner product.} \emph{Next, there exists
  the unique regularized solution } $(\varepsilon_{\gamma_1}, \mu_{\gamma_2})$ \emph{ of the
  functional (\ref{functional}) in } $ (\varepsilon_{\gamma_1}, \mu_{\gamma_2}) \in
V_{\delta ^{3(\nu_1, \nu_2) }/3}(\varepsilon^*, \mu^*).$\emph{\ The
  gradient method of the minimization of the functional
  (\ref{functional}) which starts at }$(\varepsilon_0,
\mu_0)$\emph{\ converges to the regularized solution of this
  functional. Furthermore,}

\begin{equation}\label{accur}
\begin{split}
\left\Vert \varepsilon_{\gamma_1}  -  \varepsilon^* \right\Vert &\leq \Theta_1
\left\Vert \varepsilon_0 - \varepsilon^* \right\Vert, ~~\Theta_1 \in (0,1), \\
\left\Vert  \mu_{\gamma_2}  - \mu^* \right\Vert &\leq \Theta_2
\left\Vert \mu_0  -  \mu^* \right\Vert, ~~\Theta_2 \in (0,1).
\end{split}
\end{equation}

The property(\ref{accur}) means that the regularized solution of the Tikhonov
functional (\ref{functional}) provides a better accuracy than the
initial guess $(\varepsilon_0, \mu_0)$  if it satisfies condition
(\ref{2.11}).  The next theorem presents the estimate of the norm
$\left\Vert (\varepsilon, \mu) - (\varepsilon_{\gamma_1}, \mu_{\gamma_2}) \right\Vert $
via the norm of the Fr\'{e}chet derivative of the Tikhonov functional (\ref{functional}).

\textbf{Theorem 2} \emph{Assume that conditions of Theorem 1
  hold. Then for any functions }$(\varepsilon, \mu) \in V_{(\gamma_1, \gamma_2)(\delta)}(\varepsilon^*, \mu^*)$ \emph{the following error estimate holds}
\begin{equation}
\left\Vert (\varepsilon, \mu)  - (\varepsilon_{\gamma_1(\delta)}, \mu_{\gamma_2(\delta)}) \right\Vert 
\leq \frac{2}{\delta ^{2(\nu_1, \nu_2) }}
\left\Vert P_h J^{\prime }(\varepsilon, \mu) \right\Vert 
\leq \frac{2}{\delta ^{2(\nu_1, \nu_2) }}
\left\Vert J^{\prime }(\varepsilon, \mu )
\right\Vert,  \label{4.250}
\end{equation}
\emph{which explicitly can be written as}
\begin{equation} \label{error_theorem2}
\begin{split}
\left\Vert \varepsilon  - \varepsilon_{\gamma_1(\delta)} \right\Vert 
&\leq \frac{2}{\delta ^{2 \nu_1}}
\left\Vert P_h J_{\varepsilon}^{\prime }(\varepsilon, \mu) \right\Vert 
\leq \frac{2}{\delta ^{2 \nu_1 }}
\left\Vert J_{\epsilon}^{\prime }(\varepsilon, \mu )\right\Vert
=\frac{2}{\delta ^{2 \nu_1 }}
\left\Vert L_{\varepsilon}^{\prime }(u(\varepsilon, \mu) )\right\Vert, \\
\left\Vert \mu  - \mu_{\gamma_2(\delta)} \right\Vert 
&\leq \frac{2}{\delta ^{2 \nu_2 }}
\left\Vert P_h J_{\mu}^{\prime }(\varepsilon, \mu) \right\Vert 
\leq \frac{2}{\delta ^{2 \nu_2 }} \left\Vert J_{\mu}^{\prime }(\varepsilon, \mu )\right\Vert
 = \frac{2}{\delta ^{2 \nu_2 }} \left\Vert L_{\mu}^{\prime }(u(\varepsilon, \mu))\right\Vert,
\end{split}
\end{equation}
\emph{where } $(\varepsilon_{\gamma_1(\delta)}, \mu_{\gamma_2(\delta)})$  \emph{ are minimizers of the
 Tikhonov functional (\ref{functional})  computed with regularization parameters  } $(\gamma_1(\delta), \gamma_2(\delta))$ and $
P_h: L_{2}\left( \Omega \right) \rightarrow H_{1}$\emph{\ is the
operator of orthogonal projection  of the space }$L_{2}\left( \Omega \right) $
\emph{\ on its subspace }$H_{1}$\emph{.}

\textbf{Proof.}


Since $ (\varepsilon_{\gamma_1}, \mu_{\gamma_2}) :=(\varepsilon_{\gamma_1(\delta)}, \mu_{\gamma_2(\delta)})$ is the
minimizer of the functional (\ref{functional}) on the set $G$ and
$(\varepsilon_{\gamma_1}, \mu_{\gamma_2}) \in {\rm Int}\left( G\right),$ then $P_hJ^{\prime }(\varepsilon_{\gamma_1}, \mu_{\gamma_2}) = 0$, or
\begin{equation}\label{4.2511}
\begin{split}
P_hJ^{\prime }_{\varepsilon}(\varepsilon_{\gamma_1}, \mu_{\gamma_2}) &= 0, \\
 P_hJ^{\prime }_{\mu}(\varepsilon_{\gamma_1}, \mu_{\gamma_2})  &=0. 
\end{split}
\end{equation}
Similarly with Theorem 4.11.2 of \cite{BOOK} since $ (\varepsilon, \mu)  -  (\varepsilon_{\gamma_1}, \mu_{\gamma_2}) \in H_{1},$ then
\begin{equation*}
\begin{split}
(J^{\prime }(\varepsilon, \mu) &- J^{\prime } (\varepsilon_{\gamma_1}, \mu_{\gamma_2}),  (\varepsilon, \mu) - (\varepsilon_{\gamma_1}, \mu_{\gamma_2}))  =
(P_h J^{\prime }(\varepsilon, \mu) - P_h J^{\prime }(\varepsilon_{\gamma_1}, \mu_{\gamma_2}),  (\varepsilon, \mu) - 
(\varepsilon_{\gamma_1}, \mu_{\gamma_2})).
\end{split}
\end{equation*}
Hence, using (\ref{4.249}) and (\ref{4.2511}) we can write that

\begin{equation*}
\begin{split}
\left\Vert (\varepsilon, \mu) - (\varepsilon_{\gamma_1}, \mu_{\gamma_2})  \right \Vert ^{2} 
&\leq \frac{2}{\delta ^{2(\nu_1, \nu_2) }}
\left( J^{\prime }(\varepsilon, \mu) - J^{\prime }(\varepsilon_{\gamma_1}, \mu_{\gamma_2}),
(\varepsilon, \mu) - (\varepsilon_{\gamma_1}, \mu_{\gamma_2})\right) \\
&=\frac{2}{\delta ^{2(\nu_1, \nu_2) }}\left( P_h J^{\prime }(\varepsilon, \mu)  - 
P_h J^{\prime }(\varepsilon_{\gamma_1}, \mu_{\gamma_2}) ),(\varepsilon, \mu) - (\varepsilon_{\gamma_1}, \mu_{\gamma_2}) \right) \\
&=\frac{2}{\delta ^{2(\nu_1, \nu_2) }}( P_h J^{\prime }(\varepsilon, \mu),  
(\varepsilon, \mu) - (\varepsilon_{\gamma_1}, \mu_{\gamma_2})) \\ 
&\leq
\frac{2}{\delta ^{2(\nu_1, \nu_2) }}\left\Vert P_h J ^{\prime }(\varepsilon, \mu)
 \right\Vert  \cdot \left\Vert (\varepsilon, \mu) - (\varepsilon_{\gamma_1}, \mu_{\gamma_2})   \right\Vert .
\end{split}
\end{equation*}

Thus,  from the expression above we can get
\begin{equation*}
\left\Vert (\varepsilon, \mu) - (\varepsilon_{\gamma_1}, \mu_{\gamma_2}) \right\Vert ^{2} \leq \frac{2}{
\delta ^{2(\nu_1, \nu_2) }} \left \Vert
 P_h J^{\prime }(\varepsilon, \mu ) \right\Vert \cdot \left\Vert
(\varepsilon, \mu) - (\varepsilon_{\gamma_1}, \mu_{\gamma_2}) \right\Vert.
\end{equation*}
We now divide the expression above  by $\left\Vert (\varepsilon, \mu) - (\varepsilon_{\gamma_1}, \mu_{\gamma_2})  \right\Vert $.
Using the fact that
\begin{equation*}
\left\Vert P_h J^{\prime}(\varepsilon, \mu ) \right\Vert \leq 
\left\Vert J^{\prime }(\varepsilon, \mu ) \right\Vert,
\end{equation*}
we obtain (\ref{4.250}), and using definition of the derivative of the
Tikhonov functional (\ref{derfunc}) we get (\ref{error_theorem2}),
where explicit entries of $L_{\varepsilon}^{\prime }(u(\varepsilon,
\mu)), L_{\mu}^{\prime }(u(\varepsilon, \mu))$ are given by
(\ref{grad1}), (\ref{grad2}), respectively.

$\square $

Below we reformulate Lemmas 2.1 and 3.2 of \cite{BKK} for the case of
Tikhonov functional (\ref{functional}).

\textbf{Theorem 3} \emph{Let the assumptions  of Theorems 1,2  hold.
  Let }$\left\Vert (\varepsilon^*, \mu^*) \right\Vert \leq
C,$\emph{\ with a given constant }$C$. \emph{\ We define by
}$M_{n}\subset H_{1}$\emph{\ the subspace which is obtained after
}$n$\emph{\ mesh refinements of the mesh} $K_h$.  \emph{Let } $h_n$ be
the mesh function on $M_n$ \emph{as defined in Section
    \ref{sec:spaces}.}  \emph{Then there exists the unique minimizer
  }$(\varepsilon_n, \mu_n) \in G\cap M_{n}$\emph{\ of the Tikhonov
    functional (\ref{functional}) such that the following inequalities hold}
\begin{equation}
\begin{split}
\left\Vert \varepsilon_n - \varepsilon_{\gamma_1(\delta)}\right\Vert 
&\leq
 \frac{2}{\delta ^{2 \nu_1}}   \left\Vert J^{\prime }_{\varepsilon}(\varepsilon, \mu )  \right\Vert, \\
\left\Vert  \mu_n - \mu_{\gamma_2(\delta)})  \right\Vert 
&\leq \frac{2}{
\delta ^{2 \nu_2}}\left\Vert J_{\mu}^{\prime }(\varepsilon, \mu )  \right\Vert.  \label{4.253}
\end{split}
\end{equation}

Now we present relaxation property of mesh refinements for the Tikhonov
functional (\ref{functional})  which follows from the
Theorem 4.1  of \cite{BKK}.
  
\textbf{Theorem 4}
\textbf{. }\emph{Let the assumptions of Theorems 2, 3 hold.  Let
}$(\varepsilon_n, \mu_n) \in V_{\delta ^{3\mu }}\left(\varepsilon^*, \mu^*\right) \cap
M_{n}$\emph{\ be the minimizer of the Tikhonov functional
  (\ref{functional}) on the set }$G\cap M_{n}.$ \emph{ The existence of
  the minimizer is guaranteed by Theorem 3.
 Assume that the
  regularized solution }$  (\varepsilon, \mu) \neq
(\varepsilon_n, \mu_n)$ \emph{\ which means that }$  (\varepsilon, \mu)  \notin M_{n}.$
\emph{Then the following relaxation properties hold}
\begin{equation*}
\begin{split}
\left\Vert \varepsilon_{n+1} -  \varepsilon \right\Vert &\leq \eta_{1,n}
\left\Vert  \varepsilon_{n} -  \varepsilon \right\Vert, \\
\left\Vert  \mu_{n+1} -  \mu \right\Vert &\leq \eta_{2,n}
\left\Vert   \mu_{n} - \mu \right\Vert
\end{split}
\end{equation*}
for $\eta_{1,n}, \eta_{2,n} \in (0,1)$.

\section{General framework of a posteriori error estimate}

\label{sec:general}

In this section we briefly present a posteriori error estimates for three kinds of errors:
\begin{itemize}

\item for the error $|L(u) - L(u_h)|$ in the Lagrangian 
(\ref{lagrangian});

\item for the error $|J(\varepsilon, \mu) - J(\varepsilon_h, \mu_h)|$ in the Tikhonov functional 
(\ref{functional});

\item for the errors $|\varepsilon - \varepsilon_h|$ and $|\mu -
  \mu_h|$ in the regularized solutions $\varepsilon, \mu$ of this functional.

\end{itemize}

Here, $u_h, \varepsilon_h, \mu_h$ are finite element approximations of the functions $u, \varepsilon, \mu$, respectively.
A posteriori error estimate in the Lagrangian was already derived in
\cite{BMaxwell2} for the case when only the function $\varepsilon(x)$
in system (\ref{E_gauge}) is unknown. In \cite{Bondestam1, Bondestam2}
were derived a posteriori error estimates in the Lagrangian which
corresponds to modified system (\ref{E_gauge}) for $\mu=1$.  A posteriori error in
the Lagrangian (\ref{lagrangian}) can be derived straightforwardly from a
posteriori error estimate presented in \cite{BMaxwell2} and thus, all
details of this derivation are not presented here.

However, to make clear how   a
posteriori errors in the Lagrangian and in the Tikhonov functional can be obtained,
we present general framework for them.
First we note that
\begin{equation} \label{femerrors}
\begin{split}
J(\varepsilon, \mu) - J(\varepsilon_h, \mu_h) &= J_{\varepsilon}^{\prime}(\varepsilon_h, \mu_h)(\varepsilon - \varepsilon_h) +
J_{\mu}^{\prime}(\varepsilon_h, \mu_h)(\mu - \mu_h) +  R( \varepsilon, \varepsilon_h) +  R(\mu, \mu_h), \\
L(u) -  L(u_h) &= L^{\prime }(u_h)(u- u_h) + R(u,u_h), \\
\end{split}
\end{equation}
where $R( \varepsilon, \varepsilon_h), R(\mu, \mu_h),  R(u,u_h),$ are remainders
of the second order.  We assume that $(\varepsilon_h, \mu_h)$ are
located in the small neighborhood of the regularized solutions
$(\varepsilon, \mu)$, correspondingly. Thus, since the terms $ R(u,u_h), R(
\varepsilon, \varepsilon_h), R(\mu, \mu_h)$ are of the second order then they will be small and we can
neglect them in (\ref{femerrors}).

We now use the splitting
\begin{equation}\label{splitting}
\begin{split}
u - u_h &= (u - u_h^I) + (u_h^I - u_h), \\
\varepsilon - \varepsilon_h &= (\varepsilon - \varepsilon_h^I) + (\varepsilon_h^I -\varepsilon_h), \\
\mu - \mu_h &= (\mu - \mu_h^I) + (\mu_h^I - \mu_h),
\end{split}
\end{equation}
 together with
the Galerkin orthogonality principle
\begin{equation}
\begin{split}
L^{\prime }(u_h)(\bar{u}) &= 0~~ \forall \bar{u} \in U_h,\\
 J^{\prime}(z_h)(b) &= 0~~ \forall b \in V_h,
\end{split}
\end{equation}
insert (\ref{splitting})  into (\ref{femerrors}) and  get the following error representations:
\begin{equation}\label{errorfunc} 
\begin{split}
 L(u) - L(u_h) &\approx L^{\prime}(u_h)( u - u_h^I), \\ J(\varepsilon,
 \mu) - J(\varepsilon_h, \mu_h) & \approx
 J_{\varepsilon}^{\prime}(\varepsilon_h, \mu_h)(\varepsilon -
 \varepsilon_h^I) + J_{\mu}^{\prime}(\varepsilon_h, \mu_h)(\mu - \mu_h^I).
\end{split}
\end{equation}
In (\ref{splitting}), (\ref{errorfunc}) functions  $u_h^I \in U_h$ and $\varepsilon_h^I, \mu_h^I \in V_h$  denote the interpolants of $u, \varepsilon, \mu$, respectively.

Using (\ref{errorfunc}) we conclude that a posteriori error estimate
in the Lagrangian involves the derivative of the Lagrangian
$L^{\prime}(u_h)$ which  we define as a  residual, multiplied by weights
$u - u_h^I$.  Similarly, a posteriori error
estimate in the Tikhonov functional involves the derivatives of the
Tikhonov functional $J_{\varepsilon}^{\prime}(\varepsilon_h, \mu_h)$ and
$J_{\mu}^{\prime}(\varepsilon_h, \mu_h)$ which represents residuals, multiplied by weights
$\varepsilon - \varepsilon_h^I$ and $\mu - \mu_h^I$, correspondingly.

To derive the errors $|\varepsilon - \varepsilon_h|$ and $|\mu -
\mu_h|$ in the regularized solutions $\varepsilon, \mu$ of the
functional (\ref{functional}) we will use the convexity property of
the Tikhonov functional together with the interpolation property
(\ref{2.6}).  We now make both error estimates more explicit.

\section{A posteriori error estimate in the regularized solution}

\label{sec:adaptrelax}

In this section we formulate theorem for a posteriori error estimates
$|\varepsilon - \varepsilon_h|$ and $|\mu - \mu_h|$ in the regularized
solution $\varepsilon, \mu$ of the functional (\ref{functional}).
During the proof we reduce notations and denote the scalar product
$(\cdot, \cdot)_{L_2}$ as $(\cdot, \cdot)$, as well as we denote the norm
$\left \Vert\cdot, \cdot \right \Vert_{L_2}$ as $\left \Vert\cdot, \cdot\right \Vert $. However, if norm
should be specified,  we will write it explicitly.

\textbf{Theorem 5}

 \emph{ Let the assumptions of Theorems 1,2 hold. Let $z_h
   =(\varepsilon_h, \mu_h) \in W_h$ be a finite element approximations
   of the regularized solution $z=(\varepsilon, \mu)$ on the finite
   element mesh $K_h$. Then there exists a constant $D$ defined in
   (\ref{2.10}) such that the following a posteriori error estimates
   hold}
\begin{equation}\label{theorem1}
\begin{split}
 \left \Vert  \varepsilon - \varepsilon_h \right \Vert  &\leq \frac{D}{\alpha_1} C_I \left (h ||  \varepsilon_h ||  + \left \Vert [\varepsilon_h] \right \Vert  \right ) =\frac{2 D}{  \delta^{2 \nu_1}} C_I  \left ( h ||  \varepsilon_h ||  + \left \Vert [\varepsilon_h] \right \Vert  \right ) ~ \forall \varepsilon_h \in V_h, \\
 \left \Vert \mu - \mu_h \right \Vert  &\leq \frac{D}{\alpha_2} C_I  \left ( h \left \Vert  \mu_h  \right \Vert  + \left \Vert [\mu_h] \right \Vert \right ) =  \frac{2 D}{\ \delta^{2 \nu_2}  } C_I  \left (h  \left \Vert   \mu_h  \right \Vert  + \left \Vert [\mu_h] \right \Vert \right ) 
~ \forall \mu_h \in V_h.
\end{split}
\end{equation}

\textbf{Proof.}       

 Let $z_h=(\varepsilon_h, \mu_h)$ be the minimizer of the Tikhonov functional
(\ref{functional}).
 The existence and uniqueness of this minimizer  is
guaranteed by Theorem 2. By the Theorem 1, the functional
(\ref{functional}) is strongly convex on the space $L_2$ with the strong
convexity constants  $(\alpha_1, \alpha_2) = (\gamma_1/2, \gamma_2/2)$.
 This  fact implies, see (\ref{4.249}),  that
\begin{equation} \label{4.222}
(\alpha_1, \alpha_2) \left\Vert z - z_h \right \Vert_{L_{2}(\Omega)} ^{2}
\leq  \left(J^{\prime} (z) - J^{\prime}(z_h), z - z_h \right), 
\end{equation}
where $J^{\prime}\left(z_h\right), J^{\prime}(z)$ are the Fr\'{e}chet
derivatives of the functional (\ref{functional}).

Using (\ref{4.222})  with the  splitting
\begin{equation*}
z - z_h =\left( z -  z_h^I \right) + \left(z_h^I - z_h \right) ,
\end{equation*}
where $z_h^I$ is the standard interpolant of $z$,
  and combining it with  the Galerkin orthogonality principle  
\begin{equation}
\left(J^{\prime }(z_h) - J^{\prime}(z), z_h^I -  z_h \right) =0
\label{4.223}
\end{equation}
such that
$(z_h, z_h^I) \in  W_h$,  we will obtain
\begin{equation} \label{4.224}
(\alpha_1, \alpha_2) \left \Vert z -z_h \right \Vert_{L_2}^{2}
 \leq  (J^{\prime}(z) - J^{\prime}(z_h), z - z_h^I). 
\end{equation}

The
right-hand side of (\ref{4.224}) can be estimated 
 using (\ref{2.10}) as 
\begin{equation*}
\left(J^{\prime}\left( z \right) - J^{\prime}(z_h), z - z_h^I \right) \leq D || z - z_h || \cdot ||  z - z_h^I  ||.
\end{equation*}
Substituting above equation into (\ref{4.224})  we obtain 
\begin{equation}
|| z - z_h || \leq \frac{D}{(\alpha_1, \alpha_2) } || z - z_h^I ||.  \label{theorem1_1}
\end{equation}
Using the interpolation property (\ref{2.6})  
\begin{equation*}
 || z-z_h^I||_{L^2(\Omega)} \leq C_I h || z||_{H^1(\Omega)}
\end{equation*}
 we get a posteriori error estimate for the regularized solution $z$ with the interpolation constant $C_I$:
\begin{equation}\label{theorem1_1_1}
|| z - z_h || \leq \frac{D}{(\alpha_1, \alpha_2) } || z -z_h^I || \leq   \frac{D}{(\alpha_1, \alpha_2) } C_I h ||z||_{H^1(\Omega)}.
\end{equation}

 We can estimate  $ h ||z||_{H^1(\Omega)}$  as
\begin{equation}\label{theorem1_2}
\begin{split}
h ||z||_{H^1(\Omega)}  &\leq \sum_K h_K || z||_{H^1(K)} = \sum_K    || (z  +  \nabla z) ||_{L_2(K)} h_K \\
&\leq  \sum_K \left(h_K || z_h||_{L_2(K)} + \left \Vert \frac{|[z_h]|}{h_K} h_K \right \Vert_{L_2(K)}  \right) \\
&\leq  h || z_h ||_{L_2(\Omega)} +   \sum_K (\left \Vert [z_h] \right \Vert_{L_2(K)}).
\end{split}
\end{equation}
 We denote in (\ref{theorem1_2}) by $[z_h]$  the jump of the function $z_h$ over the element $K$, $h_K$
 is the diameter of the element $K$.
In (\ref{theorem1_2}) we also used the fact 
that \cite{JS}
\begin{equation} \label{jumps}
|\nabla z | \leq \frac{|[z_h]|}{h_K}.
\end{equation}

Substituting the above estimates into the right-hand side of (\ref{theorem1_1_1})
we get 
\begin{equation*}
|| z  - z_h || \leq   \frac{D}{(\alpha_1, \alpha_2) } C_I  h ||  z_h ||  +  
 \frac{D}{(\alpha_1, \alpha_2) } C_I \left \Vert [z_h] \right \Vert ~ \forall z_h \in W_h.
\end{equation*}

Now taking into account $z_h=(\varepsilon_h, \mu_h)$ we get estimate
(\ref{theorem1}) for $|\varepsilon - \varepsilon_h|$ and $|\mu -
\mu_h|$, correspondingly.

 $\square $

\section{A posteriori error estimates for the Tikhonov functional}

\label{sec:errorfunc}

In Theorem 2 we derive a posteriori error estimates for the error in
the Tikhonov functional (\ref{functional})  obtained on the finite element mesh  $K_h$.

\textbf{Theorem 6}

\emph{\ Suppose that there exists minimizer $(\varepsilon, \mu) \in
  H^1(\Omega)$ of the Tikhonov functional (\ref{functional}) on the
  mesh $K_h$. Suppose also that there exists a finite element
  approximation $z_h =(\varepsilon_h, \mu_h)$ of $z=(\varepsilon, \mu)$ of
  $J(\varepsilon, \mu)$ on the set $W_h$ and mesh $K_h$ with the mesh
  function $h$.  Then the following approximate a posteriori error
  estimate for the error $ e=| J(\varepsilon, \mu) - J(\varepsilon_h, \mu_h) |$ in the
  Tikhonov functional (\ref{functional}) holds}
\begin{equation}\label{theorem2}
\begin{split}
e= | J(\varepsilon, \mu) - J(\varepsilon_h, \mu_h) | &\leq C_I   \left ( \left\| J_{\varepsilon}^{\prime}(\varepsilon_h, \mu_h)\right\| (
h ||  \varepsilon_h || + \left \Vert [\varepsilon_h] \right \Vert  \right) \\
&+ \left\| J_{\mu}^{\prime}(\varepsilon_h, \mu_h)\right\| \left( h || \mu_h || + \left \Vert [\mu_h] \right \Vert  \right) ) \\
&=
 C_I   (\left\| L_{\varepsilon}^{\prime}(u(\varepsilon_h, \mu_h))\right\| \left(
 h ||\varepsilon_h || + \left \Vert [\varepsilon_h] \right \Vert  \right) \\
&+ \left\| L_{\mu}^{\prime}(u(\varepsilon_h, \mu_h))\right\| \left(
h ||  \mu_h || + \left \Vert [\mu_h] \right \Vert  \right) ).
\end{split}
\end{equation}

\textbf{Proof}

By the definition of the Frech\'{e}t derivative of the Tikhonov
functional (\ref{functional}) with $z=(\varepsilon, \mu), z_h=(\varepsilon_h, \mu_h)$ we can write that on the mesh $K_h$
\begin{equation}\label{theorem2_1}
J(z) - J(z_h) = J'(z_h)(z - z_h) + R(z, z_h),
\end{equation}
where remainder $R(z, z_h) =O((z - z_h)^2),~~
(z - z_h) \to 0 ~~\forall z, z_h \in W_h$ and $
J^{\prime }(z_h)$ is the Fr\'{e}chet derivative of the
functional (\ref{functional}).  We can neglect the term
$R(z, z_h)$ in the estimate (\ref{theorem2_1}) since it is small.
This is because we assume that $z_h$ is the minimizer
of the Tikhonov functional on the mesh $K_h$ and this minimizer is located
in a small neighborhood of the regularized solution $z$.
For similar results  for the case of a general nonlinear operator equation  we refer to \cite{BKS, BKK}.
We again use the splitting
\begin{equation}
z - z_h =  z - z_h^I +  z_h^I -  z_h
\end{equation}
and the Galerkin orthogonality \cite{EEJ} 
\begin{equation}
 J'(z_h)( z_h^I -  z_h) = 0 \text{   } \forall z_h^I, z_h \in  W_h
\end{equation}
to get
\begin{equation}\label{dertikh}
J(z) - J(z_h) \leq J'(z_h)(z - z_h^I),
\end{equation}
where $z_h^I$ is a standard interpolant of $z$ on the mesh $K_h$ \cite{EEJ}.
Using (\ref{dertikh}) we can also write
\begin{equation}\label{theorem2_2}
|J(z) - J(z_h) | \leq
 || J'(z_h)|| \cdot ||z - z_h^I||,
\end{equation}
where the term $||z - z_h^I||$  can be estimated through the
interpolation estimate
\begin{equation*}
||z - z_h^I||_{L_2(\Omega)} \leq C_I || h~z||_{H^1(\Omega)}.
\end{equation*}
Substituting above estimate into (\ref{theorem2_2}) we get
\begin{equation}\label{theorem2_3}
 | J(z) - J(z_h) |  \leq C_I  \left\| J^{\prime }(z_h)\right\| h ||z||_{H^1\Omega)}.
\end{equation}
Using (\ref{jumps})
 we can estimate  $ h || z||_{H^1(\Omega)}$  similarly with (\ref{theorem1_2}) to get
\begin{equation}\label{theorem2_4}
 | J(z) - J(z_h) | \leq C_I  \left\| J^{\prime}(z_h)\right\|  \left(
h || z_h || + \left \Vert [z_h] \right \Vert  \right)~ \forall z_h \in W_h.
\end{equation}

Now taking into account $z_h=(\varepsilon_h, \mu_h)$ and using (\ref{derfunc}) we get estimate
(\ref{theorem2}) for $|J(\varepsilon, \mu) - J(\varepsilon_h, \mu_h)|$.

$\square$

\section{Mesh refinement  recommendations}

\label{sec:ref}

In this section we will show how to use
 Theorems 5 and 6 for the local mesh refinement recommendation.
This recommendation will allow  improve
accuracy of the reconstruction of the regularized solution
$(\varepsilon, \mu)$ of our problem \textbf{IP}.

Using the estimate (\ref{theorem1}) we observe that the main
contributions of the norms of the reconstructed functions
$(\varepsilon_h, \mu_h)$ are given by neighborhoods of thus points in
the finite element mesh $K_h$ where computed values of $|h
\varepsilon_h|$ and $|h \mu_h|$
achieve its maximal values.

 We also note that terms with jumps in the estimate (\ref{theorem1})
 disappear in the case of the conforming finite element meshes and
 with $(\varepsilon_h, \mu_h) \in V_h$.  Our idea of the local
 finite element mesh refinement is that it should be refined all
 neighborhoods of all points in the mesh $K_h$ where the functions $|h
 \varepsilon_h| $ and $|h \mu_h|$ achieves its maximum values.

Similarly, the estimate (\ref{theorem2}) of  Theorem 6 gives us the
idea where locally refine the finite element mesh $K_h$ to improve the
accuracy in the Tikhonov functional (\ref{functional}).  Using the
estimate (\ref{theorem2}) we observe that the main contributions of
the norms in the right-hand side of (\ref{theorem2}) are given by
neighborhoods of thus points in the finite element mesh $K_h$ 
 where computed values of $|h \varepsilon_h|$, $|h \mu_h|$, as well as 
computed values of 
 $|J_{\varepsilon}'(\varepsilon_h, \mu_h)|,  |J_{\mu}'(\varepsilon_h, \mu_h)|$ 
achieve its maximal values.

Recalling (\ref{derfunc}) and (\ref{grad1}), (\ref{grad2})  we have
\begin{equation} \label{dertikhonov1} 
\begin{split}
J_{\varepsilon}'(\varepsilon_h, \mu_h)(x)
 =  &- \int_0^T (\partial_t \lambda~ \partial_t E)~ (x,t)~dt  
+ s \int_0^T (\nabla \cdot  E)~ (\nabla \cdot \lambda) (x,t)~dt \\
 &-\lambda(x,0)f_1(x) + \gamma_1 (\varepsilon_h - \varepsilon_0)(x), ~ x \in \Omega,  
\end{split}
\end{equation}
\begin{equation} \label{dertikhonov1a} 
J_{\mu}'(\varepsilon_h, \mu_h)(x)
 =  - \int_0^T (\mu_h^{-2}~\nabla \times E ~\nabla \times \lambda) (x,t)~dt 
 +\gamma_2 (\mu_h - \mu_0)(x),~ x \in \Omega.
\end{equation}

 Thus, the second idea where to refine the finite element mesh $K_h$
 is that the neighborhoods of all points in $K_h$ where
 $|J_{\varepsilon}'(\varepsilon_h, \mu_h)|+ |J_{\mu}'(\varepsilon_h, \mu_h)|$ achieve its maximum, or both functions $|h
 \varepsilon_h| +|h \mu_h| $ and  $|J_{\varepsilon}'(\varepsilon_h, \mu_h)|+ |J_{\mu}'(\varepsilon_h, \mu_h)|$ 
 achieve their maximum, should be refined.
We include  the term $|h
 \varepsilon_h| +|h \mu_h| $  in the first mesh
refinement recommendation, and the term  $|J_{\varepsilon}'(\varepsilon_h, \mu_h)|+ |J_{\mu}'(\varepsilon_h, \mu_h)|$ 
in the second mesh refinement recommendation.
 In our computations of Section \ref{sec:num} we use 
 the first  mesh refinement recommendation and check
 performance of this mesh refinement criteria.

 \textbf{The First Mesh Refinement Recommendation  for IP.} \emph{Applying 
  Theorem 5  we conclude that we should refine the mesh in
  neighborhoods of those points in }$\Omega_{FEM}$\emph{\ where the
  function }$|h \varepsilon_h| +|h \mu_h| $\emph{\ attains its maximal
  values. More precisely, we refine the mesh in such subdomains of }$\Omega_{FEM}
$\emph{\ where}
\begin{equation*}
|h \varepsilon_h| +|h \mu_h|  \geq \widetilde{\beta} \max \limits_{{\Omega_{FEM}}} (|h \varepsilon_h| +|h \mu_h|), 
\end{equation*}
\emph{where $ \widetilde{\beta} \in (0,1)$ is the number which should
be chosen computationally and $h$ is the mesh function (\ref{meshfunction}) of the finite element mesh $K_h$}.
\\
\textbf{The Second Mesh Refinement Recommendation for IP.} \emph{Using
  Theorem 6 we  conclude that we should refine the mesh in
  neighborhoods of those points in }$ \Omega_{FEM} $\emph{\ where the
  function } $|J_{\varepsilon}'(\varepsilon_h, \mu_h)|+ |J_{\mu}'(\varepsilon_h, \mu_h)|$ 
\emph{\ attains its maximal
  values. More precisely, let }$\beta \in (0,1)
$\emph{\ be the tolerance number which should be chosen in
  computational experiments. Refine the mesh $K_h$ in such subdomains of }$
\Omega_{FEM} $\emph{\ where}
\begin{equation*}
 |J_{\varepsilon}'(\varepsilon_h, \mu_h) + J_{\mu}'(\varepsilon_h, \mu_h)| 
  \geq \beta \max_{{\Omega_{FEM}}}(|J_{\varepsilon}'(\varepsilon_h, \mu_h) + J_{\mu}'(\varepsilon_h, \mu_h)|).
\end{equation*}

\textbf{Remarks}

\begin{itemize}

\item
1. We note that in (\ref{dertikhonov1}), (\ref{dertikhonov1a}) we have
exact values of $E(x,t), \lambda(x,t)$ obtained with the computed
functions $(\varepsilon_h, \mu_h)$. However, in our algorithms of
Section \ref{sec:alg} and in computations of Section \ref{sec:num} we
approximate exact values of $E(x,t), \lambda(x,t)$ by the computed
ones $E_h(x,t), \lambda_h(x,t)$.

\item 2. In both mesh refinement recommendations we used the fact that
  functions $\varepsilon, \mu$ are unknown only in $\Omega_{FEM}$.

\end{itemize}

\section{Algorithms for solution IP}

\label{sec:alg}

In this section we will present three different algorithms which can
be used for solution of our \textbf{IP}: usual conjugate gradient
algorithm and two different adaptive finite element
algorithms. Conjugate gradient algorithm is applied on every finite
element mesh $K_h$ which we use in computations. We note that in our
adaptive algorithms we refine not only the space mesh $K_h$ but also
the time mesh $J_{\tau}$ accordingly to the CFL condition of
\cite{CFL67}.  However, the time mesh $J_{\tau}$ is refined globally
and not locally.  It can be thought as a new research task to check how
will  adaptive finite element method  work when both space and time
meshes are refined locally.

Taking into account remark of Section  \ref{sec:ref} we
denote by
\begin{equation} \label{dertikhonov2} 
\begin{split}
g_{\varepsilon}^n(x)
 =  &- \int_0^T (\partial_t \lambda_h~ \partial_t E_h)~ (x,t,\varepsilon_h^n, \mu_h^n)~dt  
+ s \int_0^T (\nabla \cdot E_h)~ ( \nabla \cdot \lambda_h) (x,t,\varepsilon_h^n, \mu_h^n )~dt \\
&-\lambda_h(x,0)f_1(x) +\gamma_1 (\varepsilon_h^n - \varepsilon_0)(x), ~ x \in \Omega,  
\end{split}
\end{equation}
\begin{equation} \label{dertikhonov3} 
g_{\mu}^n(x)
 =  - \int_0^T ((\mu_h^n)^{-2}~\nabla \times E_h ~\nabla \times \lambda_h) (x,t, \varepsilon_h^n, \mu_h^n)~dt 
 +\gamma_2 (\mu_h^n - \mu_0)(x),~ x \in \Omega,
\end{equation}
where functions $\lambda_h,  E_h$  are approximated finite element solutions of state and adjoint problems computed with $\varepsilon :=\varepsilon_h^n$ and $\mu := \mu_h^n$, respectively, and $n$ is the
number of iteration in the conjugate gradient algorithm.

 \subsection{Conjugate Gradient Algorithm}
\label{sec:cg}

\begin{itemize}
\item[Step 0.]  Discretize the computational space-time domain $\Omega
  \times [0,T]$ using partitions $K_{h}$ and $J_{\tau}$, respectively,
  see Section \ref{sec:spaces}.  Start with the initial approximations
  $\varepsilon_{h}^{0}= \varepsilon_0$ and $\mu_{h}^{0}= \mu_0$ and
  compute the sequences of $\varepsilon_{h}^{n}, \mu_{h}^{n}$ as:

\item[Step 1.]  Compute solutions $E_{h}\left(x,t,\varepsilon_{h}^{n},
  \mu_h^n\right) $ and $\lambda _{h}\left(x,t,\varepsilon_{h}^{n},
  \mu_h^n\right) $ of state (\ref{E_gauge}) and adjoint
  (\ref{adjoint}) problems, respectively, using explicit schemes
  (\ref{fem_maxwell}).

\item[Step 2.]  Update the coefficient $\varepsilon_h:=\varepsilon_{h}^{n+1}$ and $\mu_h:=\mu_{h}^{n+1}$
  on $K_{h}$ and $J_{\tau}$ via the conjugate gradient method
\begin{equation*}
\begin{split}
\varepsilon_h^{n+1} &=  \varepsilon_h^{n}  + \alpha_{\varepsilon} d_{\varepsilon}^n(x),\\
\mu_h^{n+1} &=  \mu_h^{n}  + \alpha_{\mu} d_{\mu}^n(x),
\end{split}
\end{equation*}
where
\begin{equation*}
\begin{split}
 d_{\varepsilon}^n(x)&=  -g_{\varepsilon}^n(x)  + \beta_{\varepsilon}^n  d_{\varepsilon}^{n-1}(x),\\
 d_{\mu}^n(x)&=  -g_{\mu}^n(x)  + \beta_{\mu}^n  d_{\mu}^{n-1}(x),
\end{split}
\end{equation*}
with
\begin{equation*}
\begin{split}
 \beta_{\varepsilon}^n &= \frac{|| g_{\varepsilon}^n(x)||^2}{|| g_{\varepsilon}^{n-1}(x)||^2},\\
 \beta_{\mu}^n &=   \frac{|| g_{\mu}^n(x)||^2}{|| g_{\mu}^{n-1}(x)||^2}.
\end{split}
\end{equation*}
Here, $d_{\varepsilon}^0(x)= -g_{\varepsilon}^0(x), d_{\mu}^0(x)=
-g_{\mu}^0(x)$ and $\alpha_{\varepsilon}, \alpha_{\mu} $  are
step-sizes in the gradient update which can be computed as in
\cite{Peron}.

\item[Step 3.]  Stop computing $\varepsilon_{h}^{n}$ at the iteration
  $M:=n$ and obtain the function $\varepsilon_h^M :=\varepsilon_h^n$ if
  either $||g_1^{n}||_{L_{2}( \Omega)}\leq \theta$ or norms
  $||\varepsilon_{h}^{n}||_{L_{2}(\Omega)}$ are stabilized. Here,
  $\theta$ is the tolerance in $n$ updates of the gradient method.

\item[Step 4.]  Stop computing $\mu_{h}^{n}$ at the iteration $N:=n$
  and obtain the function $\mu_h^N := \mu_h^n$ if either
  $||g_2^{n}||_{L_{2}( \Omega)}\leq \theta$ or norms
  $||\mu_{h}^{n}||_{L_{2}(\Omega)}$ are stabilized. Otherwise set
  $n:=n+1$ and go to step 1.

\end{itemize}

\subsection{Adaptive algorithms}
\label{sec:adaptalg}

In this section we present two adaptive algorithms for the solution of
our \textbf{IP}.
  In Adaptive algorithm 1 we apply first mesh refinement
recommendation of Section \ref{sec:ref}, while in Adaptive algorithm 2 we use
second mesh refinement recommendation of Section \ref{sec:ref}.

We define the minimizer of the Tikhonov functional (\ref{functional})
and its approximated finite element solution on $k$ times adaptively
refined mesh $K_{h_k}$ by $(\varepsilon, \mu)$ and $(\varepsilon_k,
\mu_k)$, correspondingly.  In our both mesh refinement recommendations
of Section \ref{sec:ref} we need compute the functions $\varepsilon_k,
\mu_k$ on the mesh $K_{h_k}$. To do that we apply conjugate gradient
algorithm of Section \ref{sec:cg}. We will define by $\varepsilon_k :=
\varepsilon_h^M, \mu_k := \mu_h^N$ values obtained at steps 3 and 4 of
the conjugate gradient algorithm.

\vspace{0.5cm}

\textbf{Adaptive Algorithm 1}

\vspace{0.5cm}

\begin{itemize}

\item[Step 0.]  Choose an initial space-time mesh ${K_h}_{0} \times
  J_{\tau_0}$ in $\Omega_{FEM} \times [0,T]$.
 Compute the sequences of
  $\varepsilon_k, \mu_k, k >0$, via following steps:

\item[Step 1.]  Obtain numerical solutions  $\varepsilon_k, \mu_k$   on $K_{h_k}$ using the  Conjugate Gradient Method of  Section \ref{sec:cg}.

\item[Step 2.] 
Refine such elements in the mesh $K_{h_k}$  where the expression
\begin{equation} \label{alg2_2}
|h \varepsilon_k| +|h \mu_k|  \geq \widetilde{\beta}_k \max_{{\Omega_{FEM}}} (|h \varepsilon_k| +|h \mu_k|)
\end{equation}
is satisfied. Here, the tolerance numbers $ \widetilde{\beta}_k  \in \left(
0,1\right) $ are chosen by the user.

\item[Step 3.]  Define a new refined mesh as $K_{h_{k+1}}$ and
  construct a new time partition $J_{\tau_{k+1}}$ such that the CFL
  condition of \cite{CFL67} for explicit schemes (\ref{fem_maxwell})
  is satisfied.  Interpolate $\varepsilon_k, \mu_k$ on a new mesh
  $K_{h_{k+1}}$ and perform steps 1-3 on the space-time mesh
    $K_{h_{k+1}} \times J_{\tau_{k+1}}$. Stop mesh refinements when
      $||\varepsilon_k - \varepsilon_{k-1}|| < tol_1$ and $||\mu_k -
      \mu_{k-1}|| < tol_2$ or $|| g_{\varepsilon}^k(x)|| < tol_3$ and
      $|| g_{\mu}^k(x)|| < tol_4$, where $tol_i, i=1,...,4$ are
      tolerances chosen by the user.

\end{itemize}

\vspace{0.5cm}

\textbf{Adaptive Algorithm 2}

\vspace{0.5cm}
\begin{itemize}

\item[Step 0.]  Choose an initial space-time mesh $K_{h_{0}}\times
  J_{\tau_0}$ in $\Omega_{FEM}$.
Compute the sequence $\varepsilon_k, \mu_k, k >0$, on a refined meshes $K_{h_k}$
  via following steps:

\item[Step 1.] Obtain numerical solutions  $\varepsilon_k, \mu_k$   on $K_{h_k} \times
  J_{\tau_k}$ using
the  Conjugate Gradient Method of
  Section \ref{sec:cg}.

\item[Step 2.] 
Refine the mesh $K_{h_k}$ at all points where 
\begin{equation}
| g_{\varepsilon}^k(x)| + |g_{\mu}^k(x) | \geq \beta_k  \max_{\Omega}( | g_{\varepsilon}^k(x)| + |g_{\mu}^k(x)|),  \label{62}
\end{equation}
where a posteriori error indicators  $g_{\varepsilon}^k, g_{\mu}^k$  are defined in
(\ref{dertikhonov1}), (\ref{dertikhonov2}).  We choose the tolerance
number $\beta_k \in \left( 0,1\right) $ in numerical examples.

\item[Step 3.] Define a new refined mesh as $K_{h_{k+1}}$ and
  construct a new time partition $J_{\tau_{k+1}}$ such that the CFL
  condition of \cite{CFL67} for explicit schemes (\ref{fem_maxwell})
  is satisfied. Interpolate $\varepsilon_k, \mu_k$ on a new mesh
  $K_{h_{k+1}}$ and perform steps 1-3 on the space-time mesh $K_{h_{k+1}}
  \times J_{\tau_{k+1}}$. Stop mesh refinements when $||\varepsilon_k
  - \varepsilon_{k-1}|| < tol_1$ and $||\mu_k - \mu_{k-1}|| < tol_2$
  or $|| g_{\varepsilon}^k(x)|| < tol_3$ and $|| g_{\mu}^k(x)|| <
  tol_4$, where $tol_i, i=1,...,4$ are tolerances chosen by the user.

\end{itemize}

\vspace{0.5cm}

\textbf{Remarks}

\begin{itemize}

\item 1. First we make comments how to choose the tolerance numbers
  $\beta_k, \widetilde{\beta_k }$ in (\ref{62}), (\ref{alg2_2}). Their
  values depend on the concrete values of $ \max \limits_{\Omega_{FEM}}( |
  g_{\varepsilon}^k(x)| + |g_{\mu}^k(x)|)$ and $\max \limits_{\Omega_{FEM}} (|h
  \varepsilon_k| +|h \mu_k|)$, correspondingly.  If we will take
  values of $\beta_k, \widetilde{\beta_k }$ which are very close to $1$ then we
  will refine the mesh in very narrow region of the $\Omega_{FEM}$,
  and if we will choose $\beta_k, \widetilde{\beta_k} \approx 0$ then
  almost all elements in the finite element mesh  will be refined, and thus, we
  will get global and not local mesh refinement.  Our numerical tests
  of Section \ref{sec:num} show that the choice of $\beta_k,
  \widetilde{\beta}_k =0.7$ is almost optimal one since with these
  values of the parameters $\beta_k, \widetilde{\beta_k }$ the finite
  element mesh $K_h$ is refined exactly at the places where we have computed the
  functions $(\varepsilon_h, \mu_h)$.

\item 2.  To compute $L_2$ norms $||\varepsilon_k - \varepsilon_{k-1}||$,
  $||\mu_k - \mu_{k-1}||$ in step 3 of adaptive algorithms the
  solutions $\varepsilon_{k-1}, \mu_{k-1}$ are interpolated from the
  mesh $K_{h_{k-1}}$ to the mesh $K_{h_k}$.


\end{itemize}

\section{Numerical studies of the adaptivity technique}

\label{sec:num}

In this section we present numerical tests for solution of our
\textbf{IP} using adaptive algorithm 1 of Section \ref{sec:adaptalg}.
Goal of our simulations is to show performance of the adaptivity
 technique in order to improve
reconstruction which was obtained on a coarse non-refined mesh.

 In our tests we reconstruct two symmetric structures of Figure \ref{fig:fig1} which
 represents model of a waveguide with small magnetic metallic
 inclusions with the relative permittivity $\varepsilon_r =12$ and
 the relative magnetic permeability $\mu_r = 2.0$. We note that we
 choose in metallic targets $\varepsilon_r =12$ similarly with our
 recent work \cite{BCN} and experimental works \cite{ BTKB,BTKF,NBKF}
 where  metallic targets were treated  as dielectrics with large dielectric
 constants and they were called  \emph{effective} dielectric constants. Values
 of them we choose similarly with \cite{BCN, BTKB, BTKF, NBKF} in the interval
\begin{equation}
\varepsilon_r \in
\left(10,30\right).  \label{2.51}
\end{equation}
In our tests we choose 
 $\mu_r = 2.0$ because the relative magnetic
permeability belongs to the interval $\mu_r \in [1,3]$, see
\cite{SSMS} and \cite{BCN} for a similar choice.

As in \cite{BCN} we initialize
  only one component $E_2$ of the electrical field $E=(E_1,E_2,E_3)$
  on $S_T$ as a plane wave $f(t)$ such that (see boundary condition in
  (\ref{E_gauge})) 
 \begin{equation}\label{f}
 \begin{split}
 f(t) =\left\{ 
 \begin{array}{ll}
 \sin \left(\omega t \right) ,\qquad &\text{ if }t\in \left(0,\frac{2\pi }{\omega }
 \right) , \\ 
 0,&\text{ if } t>\frac{2\pi }{\omega }.
 \end{array}
 \right. 
 \end{split}
 \end{equation}

 Compared with \cite{BCN} where in computations only zero initial conditions in (\ref{E_gauge})
 were used, in Test 2 of our study  we use non-zero initial
 condition for the second component $E_2$  given by the function 
\begin{equation}\label{initcond}
\begin{split}
f_0 (x) = E_2(x,0) &= \exp^{-(x_1^2 + x_2^2 + x_3^3)}  \cdot \cos  t|_{t=0} = \exp^{-(x_1^2 + x_2^2 + x_3^3)}  , \\
f_1 (x) =\frac{ \partial E_2}{\partial t} (x,0) &= -\exp^{-(x_1^2 + x_2^2 + x_3^3)} \cdot  \sin t|_{t=0} \equiv 0.
\end{split}
\end{equation}

\begin{figure}[tbp]
 \begin{center}
 \begin{tabular}{cc}
\includegraphics[width = 2.40cm, clip = true, trim = 9.0cm 3.0cm 4.0cm 2.0cm, angle = -90.0]{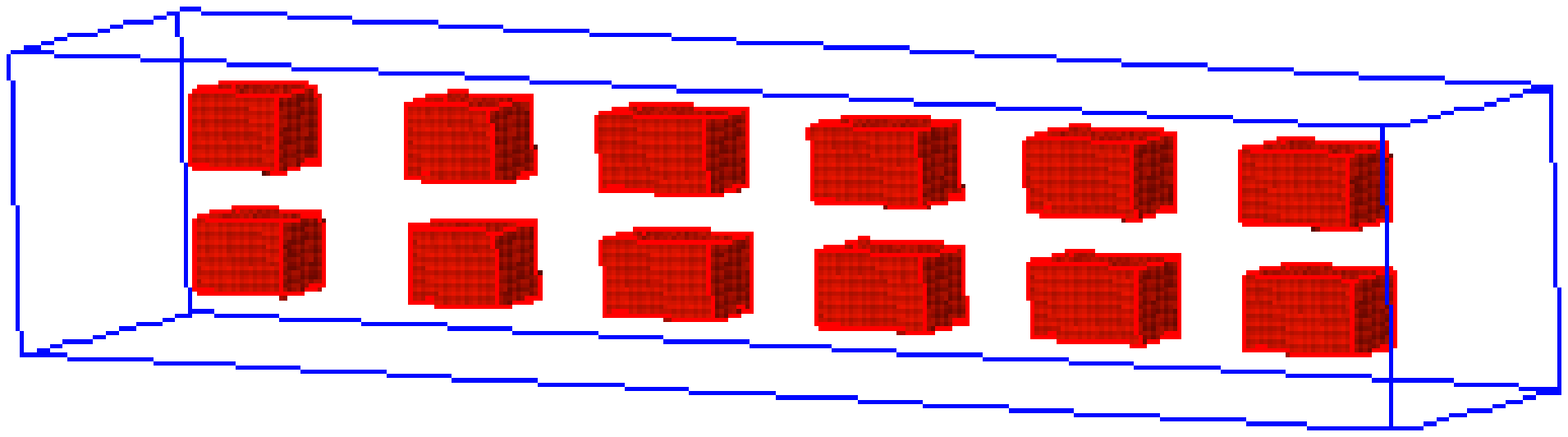} &
\includegraphics[width = 2.7cm, clip = true, trim = 8.0cm 3.0cm 4.0cm 2.0cm, angle = -90.0]{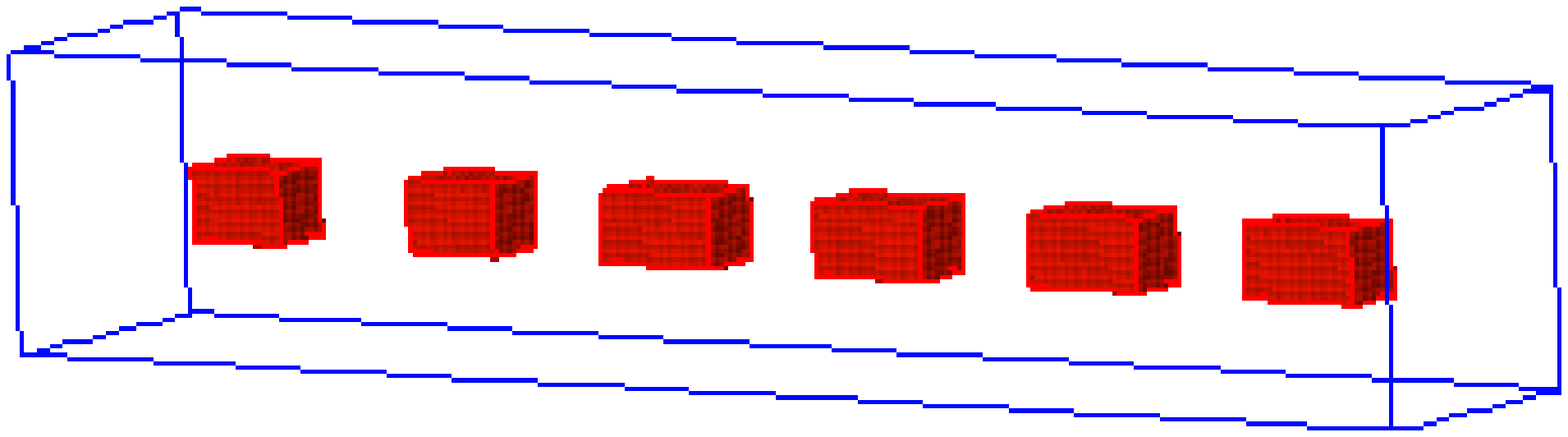} \\
 a) Test 1 &  b) Test 2 \\
\end{tabular}
 \end{center}
 \caption{The exact
 values of functions $\varepsilon(x)$ and $\mu(x)$ are:
 $\varepsilon(x)=12.0, \mu(x)=2$ inside the small scatterers, and
 $\varepsilon(x)=\mu(x)=1.0$  everywhere else in $\Omega_{\rm FEM}$.}
 \label{fig:fig2}
 \end{figure}

We perform two different tests with different inclusions to be reconstructed:
\begin{itemize}

\item
  Test 1. Reconstruction of two layers of scatterers of figure
  \ref{fig:fig2} -a) with additive noise  $\sigma=7 \%$ and
  $\sigma=17\%$ in backscattered data on the frequency interval for $\omega
  \in [45, 60]$ with zero initial conditions in (\ref{E_gauge}).

\item Test 2. Reconstruction of one layer of scatterers of figure
  \ref{fig:fig2}-b) with additive noise $\sigma=7 \%$ and
  $\sigma=17\%$ in backscattered data on the frequency interval
for $\omega   \in [45, 60]$  with  one non-zero initial condition  (\ref{initcond}) in (\ref{E_gauge}).

\end{itemize}

\subsection{Computational domains}

 For simulations of forward and adjoint problems we use the domain
 decomposition method of \cite{BMaxwell}.  This method is convenient
 for our computations since it is efficiently implemented in the software package WavES \cite{waves} using PETSc \cite{petsc}.
To apply method of \cite{BMaxwell} we divide our computational domain
$\Omega$ into two subregions as described in Section \ref{sec:stat},
and we define $\Omega_{FDM} := \Omega_{\rm OUT}$ such that $\Omega =
\Omega_{\rm FEM} \cup \Omega_{\rm FDM}$, see Figure \ref{fig:fig1}.
 In
$\Omega_{\rm FEM}$ we use finite elements and in $\Omega_{\rm FDM}$ we
will use finite difference method.  We set functions $\varepsilon(x) =
\mu(x)=1$ in $\Omega_{FDM}$ and assume that they are unknown only in
$\Omega_{FEM}$.  We choose the dimensionless domain $\Omega_{FEM}$
such that
 \begin{equation*}
 \Omega_{ FEM} = \left\{ x' = (x_1,x_2,x_3) \in (
 -3.2,3.2) \times (-0.6,0.6) \times (-0.6,0.6) \right\} .
 \end{equation*}
and the dimensionless  domain
 $\Omega$  is set to be
 \begin{equation*}
 \Omega = \left\{ x' = (x_1,x_2,x_3) \in (
 -3.4,3.4) \times (-0.8,0.8) \times (-0.8,0.8) \right\}.
 \end{equation*}
Here, the dimensionless spatial variable 
 $x^{\prime}= x/\left(1m\right)$.
  In the domain decomposition between n $\Omega_{ FEM}$ and $\Omega_{
    FDM}$ we choose the mesh size $h=0.1$.  We use also this mesh size
  for the coarse mesh ${K_h}_0$ in both adaptive algorithms of Section
  \ref{sec:adaptalg}.  As in \cite{BMaxwell, BMaxwell2, BCN} in all our tests we set $s=1$
  in (\ref{model2_1}) in $\Omega_{FEM}$.

\begin{figure}[tbp]
  \begin{center}
  \begin{tabular}{cc}
 
  {\includegraphics[scale=0.33,clip = true, trim = 6cm 6cm 6cm 3cm, angle = 0]{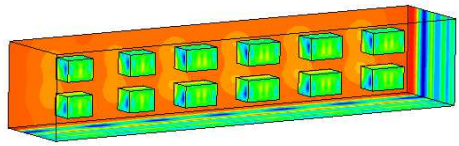}}  &
  {\includegraphics[scale=0.33, clip = true, trim = 6cm 6cm 6cm 3cm, angle = 0]{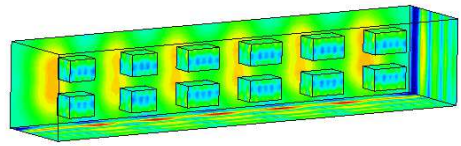}} \\
   a) FEM solution at $t= 1.5$ &   b)  FEM solution at $t= 1.8$ \\
   
  {\includegraphics[scale=0.28, clip = true, trim = 6cm 6cm 6cm 3cm, angle = -3]{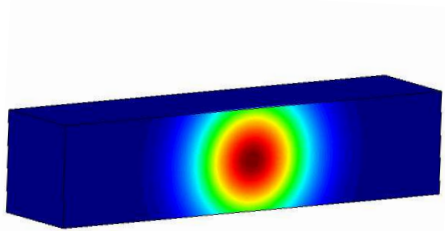}} &
  {\includegraphics[scale=0.28, clip = true, trim = 6cm 6cm 6cm 3cm, angle = -3]{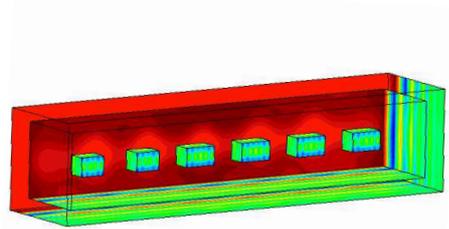}}
 \\
  c) FEM/FDM solution at $t=0$ &
  d) FEM/FDM solution at $t=1.8$ 
  \end{tabular}
  \end{center}
  \caption{Isosurfaces of the simulated  FEM/FDM 
    solution of the model problem at different
    times: a), b) Test 1;  c), d) Test 2.}
  \label{fig:Isosurfaces}
  \end{figure}

\begin{figure}
 \begin{center}
  \begin{tabular}{cccc}
 
 {\includegraphics[angle=0,width=6.0cm, clip = true, trim = 0cm 0cm 0cm 0cm, angle = 00]{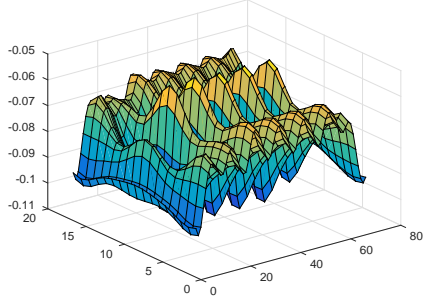}}
  &{\includegraphics[width=6.0cm, clip = true, trim = 0cm 0cm 0cm 0cm, angle = 00]{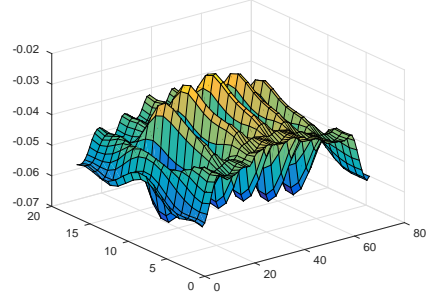}} \\
   \hspace{-1cm} (a) Test 1. $\omega = 50$, $\sigma = 7\%$, $t = 3$ &
    (b) Test 2. $\omega = 60$, $\sigma =17\%$, $t = 3$ \\
   {\includegraphics[width=6.0cm, clip = true, trim = 0cm 0cm 0cm 0cm, angle = 00]{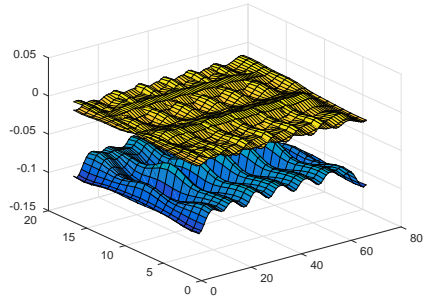}}
  &{\includegraphics[width=6.0cm, clip = true, trim = 0cm 0cm 0cm 0cm, angle = 00]{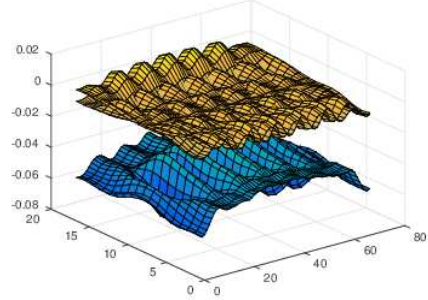}} \\
    \hspace{-1cm} (c) Test 1. $\omega = 50$, $\sigma =7\%$, $t = 3$ &
      (d) Test 2. $\omega = 60$, $\sigma =17\%$, $t = 3$ \\
   
  \end{tabular}
  \end{center}
       \caption{ a), b) Backscattered data of the one component, $E_2(x,t)$, of the electric field $E(\mathbf{x},t)$. c), d) Computed components  $E_2$ (below) and $E_1 \ \text{and}\ E_3$ (on top)
   of the backscattered electric field   $E(\mathbf{x},t)$.}
   \label{fig:backscatdata}
\end{figure}

Because of the domain decomposition
 the Maxwell's system (\ref{E_gauge}) transforms to the wave equation in
 $\Omega_{FDM}$ such that
 \begin{equation}\label{waveeq}
 \begin{split}
  \frac{\partial^2 E}{\partial t^2} - \triangle E  &=0,~ \mbox{in}~~ \Omega_{FDM} \times [0,T],    \\
   E_2(x,0) &= f_0(x), E_1(x,0) = E_3(x,0) = 0~ \mbox{ for}~~ x \in \Omega,     \\
   E_t(x,0) &= 0~ \mbox{for}~~ x \in \Omega,     \\
  E(x,t)& = (0, f\left(t\right),0) ,~ \mbox{on}
 ~\partial \Omega_{1}\times (0,t_{1}], \\
 \partial _{n}E(x,t)& =-\partial _{t} E(x,t),~\mbox{on}
 ~\partial \Omega_{1}\times (t_{1},T), \\
 \partial _{n}E(x,t)& =-\partial _{t} E(x,t),~ \mbox{on}
 ~\partial \Omega_{2}\times (0,T), \\
 \partial _{n} E(x,t)& =0,~ \mbox{on}~\partial \Omega_{3}\times (0,T). 
 \end{split}
 \end{equation}
In  $\Omega_{ FEM}$ we   solve 
 \begin{equation}\label{maxweq}
 \begin{split}
 \varepsilon \frac{\partial^2 E}{\partial t^2} + \nabla \times ( \mu^{-1} \nabla \times E)  - s\nabla  ( \nabla \cdot(\varepsilon E))  &= 0,~ \mbox{in}~~ 
 \Omega_{{\rm FEM}},    \\
   E(x,0) = 0, ~~~E_t(x,0) &= 0~ \mbox{in}~~ \Omega_{\rm FEM},    \\
 E(x,t)|_{\partial \Omega_{\rm FEM}} &= E(x,t)|_{\partial \Omega_{{ FDM}_I}}. 
 \end{split}
 \end{equation}
In (\ref{maxweq}), $\partial \Omega_{{\rm FDM}_I}$ denotes the
internal boundary of the domain $\Omega_{FDM}$, and $\partial
\Omega_{FEM}$ denotes the boundary of the domain $\Omega_{FEM}$.  In a
similar way transforms also the adjoint problem (\ref{adjoint}) into
two problems in $\Omega_{FDM}$ and in $\Omega_{FEM}$, which will be
the same as in \cite{BCN}.  We solve the forward and adjoint problems
in time $[0,T]=[0,3]$ in both adaptive algorithms and choose the time
step $\tau=0.006$ which satisfies the CFL condition \cite{CFL67}.  To
be able test adaptive algorithms we first generate backscattered data
at $S_T$ by solving the forward problem (\ref{E_gauge}) with the plane
wave $f(t)$ given by (\ref{f}) in the time interval $t=[0,3]$ with
$\tau=0.006$ and with known values of $\varepsilon_ r =12.0, \mu_r =2$
inside scatterers of Figure \ref{fig:fig2} and $\varepsilon_r = \mu_r
=1.0$ everywhere else in $\Omega$. Figure \ref{fig:Isosurfaces}
presents isosurfaces of the exact simulated solution at different
times.  Particularly, in Figure \ref{fig:Isosurfaces}-c) we observe
behaviour of non-zero initial condition (\ref{initcond}). Our data
were generated on a specially constructed mesh for the solution of the
forward problem: this mesh was several times refined in the places
where inclusions of Figure \ref{fig:fig2} are located. This mesh is
completely different than meshes used in computations in Tests 1,
2. Thus, the variational crime in our computations is avoided.
Figures \ref{fig:backscatdata}-a), b) illustrate typical behavior of
noisy backscattered data in Test 1 running it with $\omega =50$ in
(\ref{f}).  Figure \ref{fig:backscatdata}-b) shows result of
computations of the forward problem
in Test 2 when we take $\omega =60$ in (\ref{f}).  Figure
\ref{fig:backscatdata}-c),d) show the difference in backscattered data
for all components of the electrical field at final time of
computations $t=3$.

\subsection{Reconstructions}

\begin{table}[h] \label{tab:1}
{\footnotesize Table 1. \emph{Results of reconstruction on a coarse
    meshes of Tables 5,6 for $\sigma =7\%$ together with computational errors between $\max\limits_{\Omega_{FEM} } \varepsilon_{\overline{N}}$ and exact $\varepsilon^{*}$ in
    procents. Here, $\overline{N}$ is the final iteration number 
    in the conjugate gradient method for  computation of
    $\varepsilon_r$, and $\overline{M}$ is the final iteration number
    for computation of $\mu_r$.}}  \par
\vspace{2mm}
\centerline{
\begin{tabular}{|c|}
 \hline
   $\sigma= 7\%$ 
 \\
 \hline
\begin{tabular}{c|c|c|c|c|c|c}
Test 1 & $ \max\limits_{\Omega_{FEM} } \varepsilon_{\overline{N}}$ & error, $\%$  & $\overline{N}$  &
$ \max\limits_{\Omega_{FEM} } \mu_{\overline{M}}$ & error, $\%$  & $\overline{M}$  
\\ \hline
$\omega=45$ & $15$ & $25$ & $10$ &$2.58$ &$29$ &$10$   \\
$\omega=50$  & $15$ & $25$ &$10 $  &$2.38$ &$19$ &$10$   \\
$\omega=60$  & $15$ &$25$    &$10$ &$2.46$ &$23$ &$10$    \\
\end{tabular}
\\
\hline
\begin{tabular} {c|c|c|c|c|c|c}
Test 2 & $\max\limits_{\Omega_{FEM} } \varepsilon_{\overline{N}}$ & error, $\%$  & $\overline{N}$  &
$ \max\limits_{\Omega_{FEM} } \mu_{\overline{M}}$ & error, $\%$  & $\overline{M}$  
\\ \hline
$\omega=45$ & $13.32$ &$11$& $10$ & $3.07$ &$53.5$ & $10$    \\
$\omega=50$  & $15 $&$25$  &$10$  &$2.62$ &$31$ &$10$   \\
$\omega=60$  & $9.3$ &$22.4$ &$10$ &$2.88$ & $44 $&$10$   \\
\end{tabular}
\\
\hline
\end{tabular}
}

\end{table}

\begin{table}[h] \label{tab:2}
{\footnotesize Table 2. \emph{Results of reconstruction 
    on a coarse meshes of Tables 5,6 for $\sigma =17\%$ together with
    computational errors  between $\max\limits_{\Omega_{FEM} } \varepsilon_{\overline{N}}$ and exact $\varepsilon^{*}$ in procents. Here, $\overline{N}$ is the
    final iteration number in the conjugate gradient method for
    computation of $\varepsilon_r$, and $\overline{M}$ is the final
    iteration number for  computation of $\mu_r$.}}
\par
\vspace{2mm}
\centerline{
\begin{tabular}{|c|}
 \hline
   $\sigma= 17\%$ 
 \\
 \hline
\begin{tabular}{c|c|c|c|c|c|c}
Test 1 & $ \max\limits_{\Omega_{FEM} } \varepsilon_{\overline{N}}$ & error, $\%$  & $\overline{N}$  &
$ \max\limits_{\Omega_{FEM} } \mu_{\overline{M}}$ & error, $\%$  & $\overline{M}$  
\\ \hline
$\omega=45$ & $15$ & $25$ & $10$ &$2.35$ &$17.5$ &$10$   \\
$\omega=50$  & $15$ & $25$ &$10 $  &$2.89$ &$44.5$ &$10$   \\
$\omega=60$  & $15$ &$25$    &$8$ &$3.09$ &$53.6$ &$8$    \\
\end{tabular}
\\
\hline
\begin{tabular} {c|c|c|c|c|c|c}
Test 2 & $\max\limits_{\Omega_{FEM} } \varepsilon_{\overline{N}}$ & error, $\%$  & $\overline{N}$  &
$ \max\limits_{\Omega_{FEM} } \mu_{\overline{M}}$ & error, $\%$  & $\overline{M}$  
\\ \hline
$\omega=45$ & $15$ &$25$& $10$ & $2.39$ &$19.5$ & $10$    \\
$\omega=50$  & $15 $&$25$  &$10$  &$2.24$ &$12$ &$10$   \\
$\omega=60$  & $8.46$ &$29.5$ &$10$ &$2.50$ & $25 $&$10$   \\
\end{tabular}
\\
\hline
\end{tabular}
}

\end{table}

\begin{table}[h] \label{tab:3}
{\footnotesize Table 3. \emph{Results of reconstruction on a 5 times
    adaptively refined meshes of Tables 5,6 for $\sigma =7\%$ together with
    computational errors  between $\max\limits_{\Omega_{FEM} } \varepsilon_{\overline{N}}$ and exact $\varepsilon^{*}$ in procents. Here, $\overline{N}$ is the
    final iteration number in the conjugate gradient method for
    computation of $\varepsilon_r$, and $\overline{M}$ is the final
    iteration number for computation of $\mu_r$.}}  
\par
\vspace{2mm}
\centerline{
\begin{tabular}{|c|}
 \hline
   $\sigma= 7\%$ 
 \\
 \hline
\begin{tabular}{c|c|c|c|c|c|c}
Test 1 & $ \max\limits_{\Omega_{FEM} } \varepsilon_{\overline{N}}$ & error, $\%$  & $\overline{N}$  &
$ \max\limits_{\Omega_{FEM} } \mu_{\overline{M}}$ & error, $\%$  & $\overline{M}$  
\\ \hline
$\omega=45$ & $14.96$ & $24.6$ & $3$ &$1.82$ &$9$ &$3$   \\
$\omega=50$  & $14.96$ & $24.6$ &$3 $  &$1.73$ &$13.5$ &$3$   \\
$\omega=60$  & $14.95$ &$24.5$    &$3$ &$1.76$ &$12$ &$3$    \\
\end{tabular}
\\
\hline
\begin{tabular} {c|c|c|c|c|c|c}
Test 2 & $\max\limits_{\Omega_{FEM} } \varepsilon_{\overline{N}}$ & error, $\%$  & $\overline{N}$  &
$ \max\limits_{\Omega_{FEM} } \mu_{\overline{M}}$ & error, $\%$  & $\overline{M}$  
\\ \hline
$\omega=45$ & $12.97$ &$8$& $3$ & $1.99$ &$0.5$ & $3$    \\
$\omega=50$  & $14.57$ &$21.4$  &$3$  &$1.79$ &$10.5$ &$3$   \\
$\omega=60$  & $9.3$ &$22.5$ &$3$ &$1.91$ & $4.5$& $3$    \\
\end{tabular}
\\
\hline
\end{tabular}
}
\end{table}

\begin{table}[h] \label{tab:4}
{\footnotesize Table 4. \emph{Results of reconstruction on a 5 times
    adaptively refined meshes of Tables 5,6 for $\sigma =17\%$ together with
    computational errors between $\max\limits_{\Omega_{FEM} } \varepsilon_{\overline{N}}$ and exact $\varepsilon^{*}$ in procents. Here, $\overline{N}$ is the
    final iteration number in the conjugate gradient method for
    computation of $\varepsilon_r$, and $\overline{M}$ is the final
    iteration number for computation of $\mu_r$. }}  \par
\vspace{2mm}
\centerline{
\begin{tabular}{|c|}
 \hline
   $\sigma= 17\%$ 
 \\
 \hline
\begin{tabular}{c|c|c|c|c|c|c}
Test 1 & $ \max\limits_{\Omega_{FEM} } \varepsilon_{\overline{N}}$ & error, $\%$  & $\overline{N}$  &
$ \max\limits_{\Omega_{FEM} } \mu_{\overline{M}}$ & error, $\%$  & $\overline{M}$  
\\ \hline
$\omega=45$ & $14.96$ &$24.6$& $3$ & $1.65$ &$17.5$ & $3$    \\
$\omega=50$  & $14.96$ &$24.6$  &$3$  &$1.97$ &$1.5$ &$3$ \\
$\omega=60$  & $14.95$ &$24.5$    &$3$ &$2.04$ &$20$ &$3$    \\
\end{tabular}
\\
\hline
\begin{tabular}{c|c|c|c|c|c|c}
Test 2 & $ \max\limits_{\Omega_{FEM} } \varepsilon_{\overline{N}}$ & error, $\%$  & $\overline{N}$  &
$ \max\limits_{\Omega_{FEM} } \mu_{\overline{M}}$ & error, $\%$  & $\overline{M}$  
\\
 \hline
$\omega=45$ & $14.69$ & $22.4$ & $3$& $1.71$ &$14.5$ & $3$ \\
$\omega=50$  & $14.47$ & $20.5$ &$3$ &$1.63$ &$18.5$ & $3$  \\
$\omega=60$  & $8.44$ &$29.7$ &$3$ &$1.74$ &$13$ &$3$    \\
\end{tabular}
\\
\hline
\end{tabular}
}

\end{table}

\begin{table}[tbp] \label{tab:5}
{\footnotesize Table 5. \emph{Test 1. Computed values of $\varepsilon_\mathrm{r}^{\mathrm{comp}} := \max \limits_{\Omega_{FEM}} \varepsilon_r \ \text{and} \ \mu_\mathrm{r}^{\mathrm{comp}} := \max \limits_{\Omega_{FEM}}\mu_r$ on the adaptively refined meshes. Computations are done with the noise $\sigma = 7\%.$}}
  \begin{center}
   {\footnotesize
    \begin{tabular}{|c|l|r|r|r|r|r|r|} \hline
    $\omega$ & & coarse mesh & 1 ref. mesh & 2 ref. mesh & 3 ref. mesh & 4 ref. mesh & 5 ref. mesh \\ \hline
      45 & \# nodes & 10958 & 11028 & 11241 & 11939 & 14123 & 18750   \\ \hline
        & \# elements & 55296 & 55554 & 56624   & 60396 & 73010 &  96934  \\ \hline
        & $\varepsilon_\mathrm{r}^{\mathrm{comp}}$ & 15 & 15 &  15 & 15 & 15 & 14.96\\ \hline
         & $\mu_\mathrm{r}^{\mathrm{comp}}$ & 2.58 & 2.58 &  2.58 & 2.58 & 2.58 & 1.82 \\ \hline
      50 & \# nodes & 10958 & 11031 & 11212 & 11887 & 13761 & 17892 \\ \hline
        & \# elements & 55296 & 55572 & 56462  & 60146 & 71010 & 92056 \\ \hline
        & $\varepsilon_\mathrm{r}^{\mathrm{comp}}$ & 15 & 15 & 15 & 15 & 15 & 14.96 \\ \hline
         & $\mu_\mathrm{r}^{\mathrm{comp}}$ & 2.38 & 2.38 &  2.38 & 2.38 & 2.38 & 1.73\\ \hline
      60 & \# nodes & 10958 & 11050 & 11255 & 11963 & 13904 &18079 \\ \hline
        & \# elements & 55296 & 56666 & 60564 & 71892 & 61794 & 92926 \\ \hline
        &$\varepsilon_\mathrm{r}^{\mathrm{comp}}$ & 15 & 15 & 15  & 15 & 15 & 14.96  \\ \hline
         & $\mu_\mathrm{r}^{\mathrm{comp}}$ & 2.46 & 2.46 &  2.46 & 2.46 & 2.46 & 1.76 \\ \hline

    \end{tabular}
    }
    \end{center}
    
\end{table}
\begin{table}[tbp] \label{tab:6}
{\footnotesize Table 6.  \emph{Test 2. Computed values of $\varepsilon_\mathrm{r}^{\mathrm{comp}} := \max \limits_{\Omega_{FEM}} \varepsilon_r \ \text{and} \ \mu_\mathrm{r}^{\mathrm{comp}} := \max \limits_{\Omega_{FEM}}\mu_r$ on the adaptively refined meshes. Computations are done with the noise $\sigma = 17\%.$}}
  \begin{center}
   {\footnotesize
  \begin{tabular}{|c|l|r|r|r|r|r|r|} \hline
 
     $\omega$ & & coarse mesh & 1 ref. mesh & 2 ref. mesh & 3 ref. mesh & 4 ref. mesh & 5 ref. mesh \\ \hline
      45 & \# nodes & 10958 & 11007 & 11129 & 11598 & 12468 & 14614   \\ \hline
        & \# elements & 55428 & 55428 &56024  & 58628 &63708 &  74558  \\ \hline
        & $\varepsilon_\mathrm{r}^{\mathrm{comp}}$ & 15 & 15 &  15 & 15 & 15 & 14.96\\ \hline
         & $\mu_\mathrm{r}^{\mathrm{comp}}$ & 2.39 & 2.39 &  2.39 & 2.39 & 2.39 & 1.71 \\ \hline
      50 & \# nodes & 10958 & 11002 & 11106 & 11527 & 12433 & 14494 \\ \hline
        & \# elements & 55296 & 55398 & 55908  & 58240 & 63540 & 73900 \\ \hline
        & $\varepsilon_\mathrm{r}^{\mathrm{comp}}$ & 15 & 15 & 15 & 15 & 15 & 14.47 \\ \hline
         & $\mu_\mathrm{r}^{\mathrm{comp}}$ & 2.24 & 2.24 &  2.24 & 2.24 & 2.24 & 1.63\\ \hline
      60 & \# nodes & 10958 & 11002 & 11104 & 11560 & 12459 &14888 \\ \hline
        & \# elements & 55296 &55398 & 55904 & 58402 & 63628&76068  \\ \hline
        &$\varepsilon_\mathrm{r}^{\mathrm{comp}}$ & 8.46 & 8.46 &8.46& 8.46  & 8.46  & 8.44  \\ \hline
         & $\mu_\mathrm{r}^{\mathrm{comp}}$ & 2.50& 2.50 &  2.50 & 2.50 & 2.50 & 1.74 \\ \hline

    \end{tabular}
    
    }
    
  \end{center}

\end{table}

We start to run adaptive algorithms with guess values of
$\varepsilon_r =1.0, \mu_r = 1.0$ at all points in $\Omega$.  In our
recent work \cite{BCN} was shown that such choice of the initial guess
gives a good reconstruction for both functions $\varepsilon_r$ and
$\mu_r$, see also \cite{BKS, BMaxwell} for a similar choice of initial
guess for other coefficient inverse problems (CIPs).  Taking into
account (\ref{2.51}) we choose following sets of admissible parameters
for $\varepsilon_r$ and $\mu_r$
 \begin{equation}\label{admpar}
 \begin{split}
  M_{\varepsilon} \in \{\varepsilon\in C(\overline{\Omega })|1\leq \varepsilon(x)\leq 15\},\\
  M_{\mu} \in \{\mu\in C(\overline{\Omega })|1\leq \mu(x)\leq 3\}.
 \end{split}
 \end{equation}

In our simulations we choose two 
constant regularization parameters $\gamma_1 =0.01, \gamma_2=0.7$ in the 
 Tikhonov functional (\ref{functional}).
These parameters
satisfy conditions (\ref{2.11})
and  were chosen because of our
computational experience: such choices for the regularization
parameters were optimal since they gave the smallest relative
errors $e_{\varepsilon} = \frac{||\varepsilon - \varepsilon_h ||}{||\varepsilon_h ||}$ 
and $ e_{\mu} =        \frac{||\mu - \mu_h||}{||\mu_h ||}$ in the reconstruction,
       see \cite{BCN} for details.
Iteratively regularized adaptive finite element method
for our \textbf{IP}  when zero initial conditions $f_0=f_1=0$  in
(\ref{E_gauge}) are initialized, is recently presented in
\cite{samar}.  Currently we perform numerical experiments with iteratively
regularized adaptive finite element method for the case when we initialize one
non-zero initial condition (\ref{initcond})  in (\ref{E_gauge}). This work will be
described in the forthcoming paper. In the above mentioned works iterative
regularization is performed via algorithms of \cite{BKS}.
We also refer to \cite{Engl, IJT11}  for
different techniques for the choice of regularization parameters. 

 To get our reconstructions of Figures
 \ref{fig:Recos_omega45noise7_zero}~-~\ref{fig:Recos_omega50noise17_nonzero},
 we use image post-processing procedure described in \cite{BCN}.
Tables 1-6 present computed results of reconstructions for
$\varepsilon_r$ and $\mu_r$ on different adaptively refined meshes
after applying adaptive algorithm 1. Similar results are obtained for
adaptive algorithm 2, and thus they are not presented here.  

\subsubsection{Test 1}

In this example we performed simulations with two additive noise
levels in data: $\sigma= 7 \%$ and $\sigma= 17\%$, see Tables 1-6 for
results. Using these tables we observe that the best reconstruction
results for both noise levels are obtained for $\omega=45$ in
(\ref{f}). Below we describe  reconstructions
obtained with $\omega = 45$ in (\ref{f}) and $\sigma=7\%$.

  The reconstructions of $\varepsilon_r$ and $\mu_r$ on initial coarse
  mesh are presented
  in Figure \ref{fig:Recos_omega45noise7_zero_coarse}.  Using Table 1
  we observe that we achieve good values of contrast for both
  functions already on a coarse mesh. However, Figures
  \ref{fig:yz_reconstruction_zero}-a), b) show us that the locations
  of all inclusions in $x_3$ direction should be improved.  The
  reconstructions of $\varepsilon_r$ and $\mu_r$ on a final adaptively
  refined mesh are presented in Figure \ref{fig:Recos_omega45noise7_zero}.
  We observe significant improvement of reconstructions of
  $\varepsilon_r$ and $\mu_r$ in $x_3$ direction on the final
  adaptively refined mesh compared with reconstructions obtained on a
  coarse mesh, see Figure \ref{fig:yz_reconstruction_zero}.  Figures
  \ref{fig:Recos_omega45noise7_zero_mesh}-a), c), e) show different
  projections of final adaptively refined mesh which was used for
  computations of images of Figures
  \ref{fig:Recos_omega45noise7_zero},
  \ref{fig:yz_reconstruction_zero}-c), d).


\subsubsection{Test 2}

In this test we again used two additive noise levels in data, $\sigma=
7 \%$ and $\sigma= 17\%$, as well as non-zero initial condition
(\ref{initcond}) in (\ref{E_gauge}). Results of computations are
presented in Tables 1-6. Using these tables we see that the best
reconstruction results in  this test for both noise levels are obtained for $\omega=50$ 
 in (\ref{f}). We now describe
reconstructions obtained for 
 $\omega = 50$ in (\ref{f})  and $\sigma= 17\%$.

The reconstructions of $\varepsilon_r$ and $\mu_r$ on a coarse mesh
are shown in Figure \ref{fig:Recos_omega50noise17_nonzero_coarse}.
The reconstructions of $\varepsilon_r$ and $\mu_r$ on a final
adaptively refined mesh are given in Figure
\ref{fig:Recos_omega50noise17_nonzero}.  We again observe significant
improvement of reconstructions of $\varepsilon_r$ and $\mu_r$ in $x_3$
direction on the final adaptively refined mesh in comparison to the
reconstruction obtained on a coarse mesh, see Figure
\ref{fig:yz_reconstruction}. Figures
\ref{fig:Recos_omega45noise7_zero_mesh}-b), d), f) show different
projections of final adaptively refined mesh which was used for
computations of images of Figures
\ref{fig:Recos_omega50noise17_nonzero}, \ref{fig:yz_reconstruction}-c),d).


\section{Conclusion}

This work is a continuation of our previous study in \cite{BCN} and is
focused on the solution of coefficient inverse problem for simultaneously reconstruction of
functions $\varepsilon$ and $\mu$ from time-dependent backscattered
data in the Maxwell's equations.  To do that we have used optimization
approach of \cite{BCN} applied on adaptively refined meshes.

We derived a posteriori error estimates in the reconstructed
coefficients $\varepsilon$ and $\mu$ and in the Tikhonov functional to
be minimized. We then formulated two adaptive algorithms which
allow   reconstruction of
$\varepsilon$ and $\mu$ on the locally adaptively refined meshes 
 using these estimates.

\begin{figure}
\begin{center}

    \begin{tabular}{c c c}
          
 {\includegraphics[width = 3.5cm, clip = true, trim = 5.0cm 0.0cm 1.0cm 0.0cm, angle = -90.0]{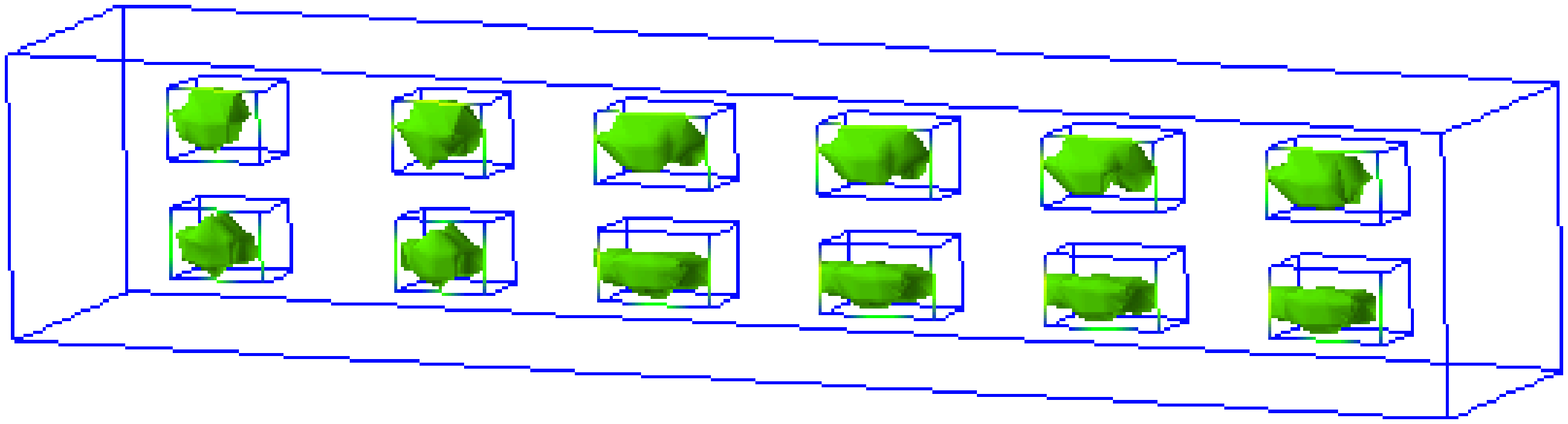}}
 & {\includegraphics[width = 3.5cm, clip = true, trim = 5.0cm 0.0cm 1.0cm 0.0cm, angle = -90.0]{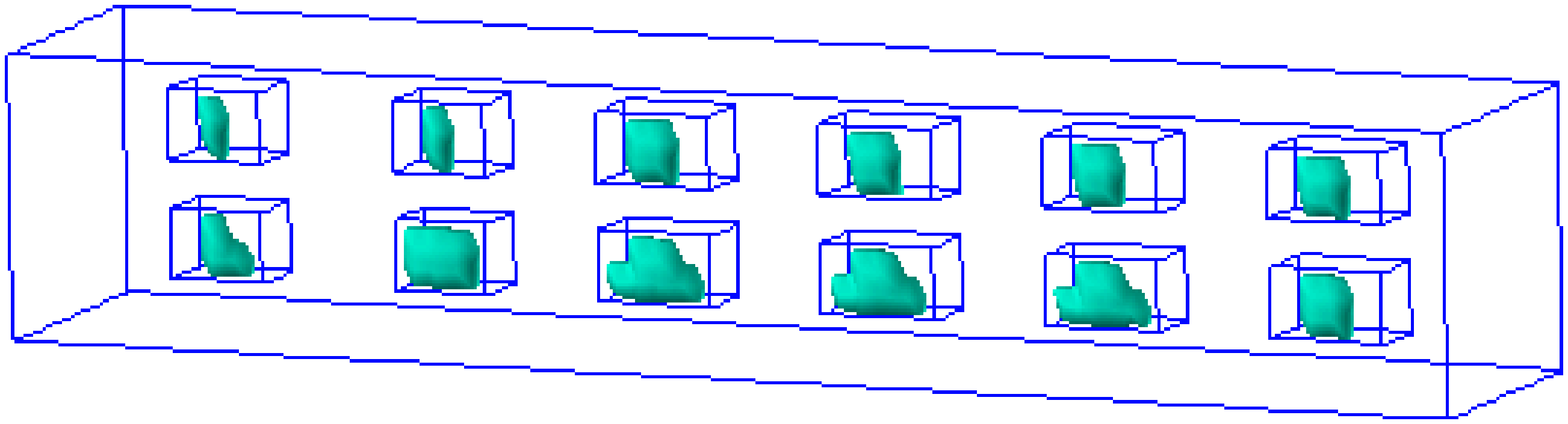}} & \\
 (a) $\max\limits_{\Omega_{FEM} }\varepsilon_r \approx 15$  &
 (b) $\max\limits_{\Omega_{FEM} }\mu_r \approx 2.5$
  \end{tabular}
 \caption{Test 1. Computed images of reconstructed functions $\varepsilon_r(\mathbf{x}) \ \text{and} \ \mu_r(\mathbf{x})$ on a coarse mesh for $\omega = 45$, $\sigma =7\%$.}
 \label{fig:Recos_omega45noise7_zero_coarse}      

\end{center}     
 \end{figure}

\begin{figure}
\begin{center}
      \begin{minipage}{1\textwidth}
 \begin{tabular}{c c c c}
 {\includegraphics[width = 3.5cm, clip = true, trim = 5.0cm 0.0cm 1.0cm 0.0cm, angle = -90.0]{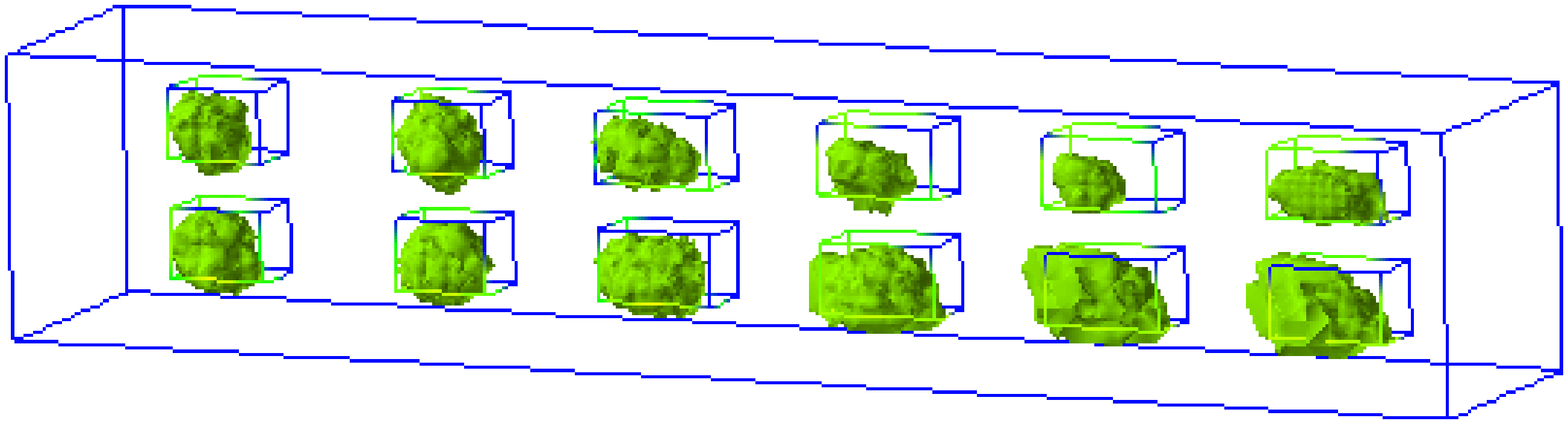}} 
 & {\includegraphics[width = 3.5cm, clip = true, trim = 5.0cm 0.0cm 1.0cm 0.0cm, angle = -90.0]{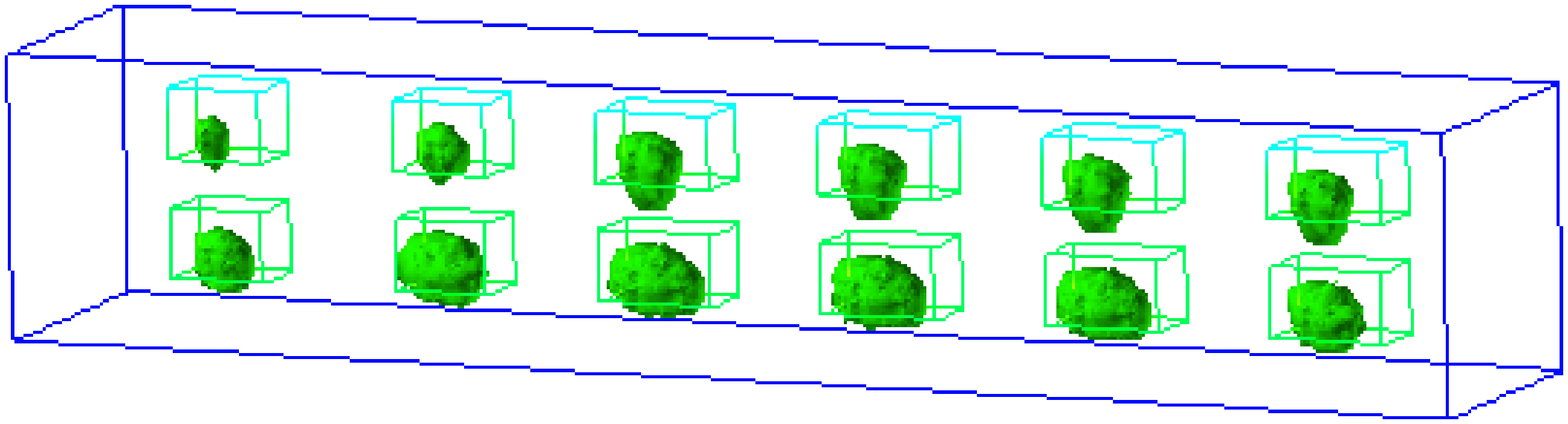}} \\
 (a) $\max\limits_{\Omega_{FEM} } \varepsilon_r \approx 14.9$  &
 (b) $\max\limits_{\Omega_{FEM} } \mu_r \approx 1.8$
  \end{tabular}
 \caption{Test 1. Computed images of reconstructed functions $\varepsilon_r(\mathbf{x}) \ \text{and} \ \mu_r(\mathbf{x})$  on a 5 times adaptively refined mesh presented in Figure \ref{fig:Recos_omega45noise7_zero_mesh}. Computations are done  for $\omega = 45$, $\sigma = 7\%$.}
 \label{fig:Recos_omega45noise7_zero}
 \end{minipage}
 \end{center}
 \end{figure}

\begin{figure}
 \begin{center}
 
    \begin{tabular} {c c c } 
      \includegraphics[width = 0.3\textwidth, clip = true, trim = 5.50cm 3.0cm 2.0cm 1.50cm,, angle = -90]{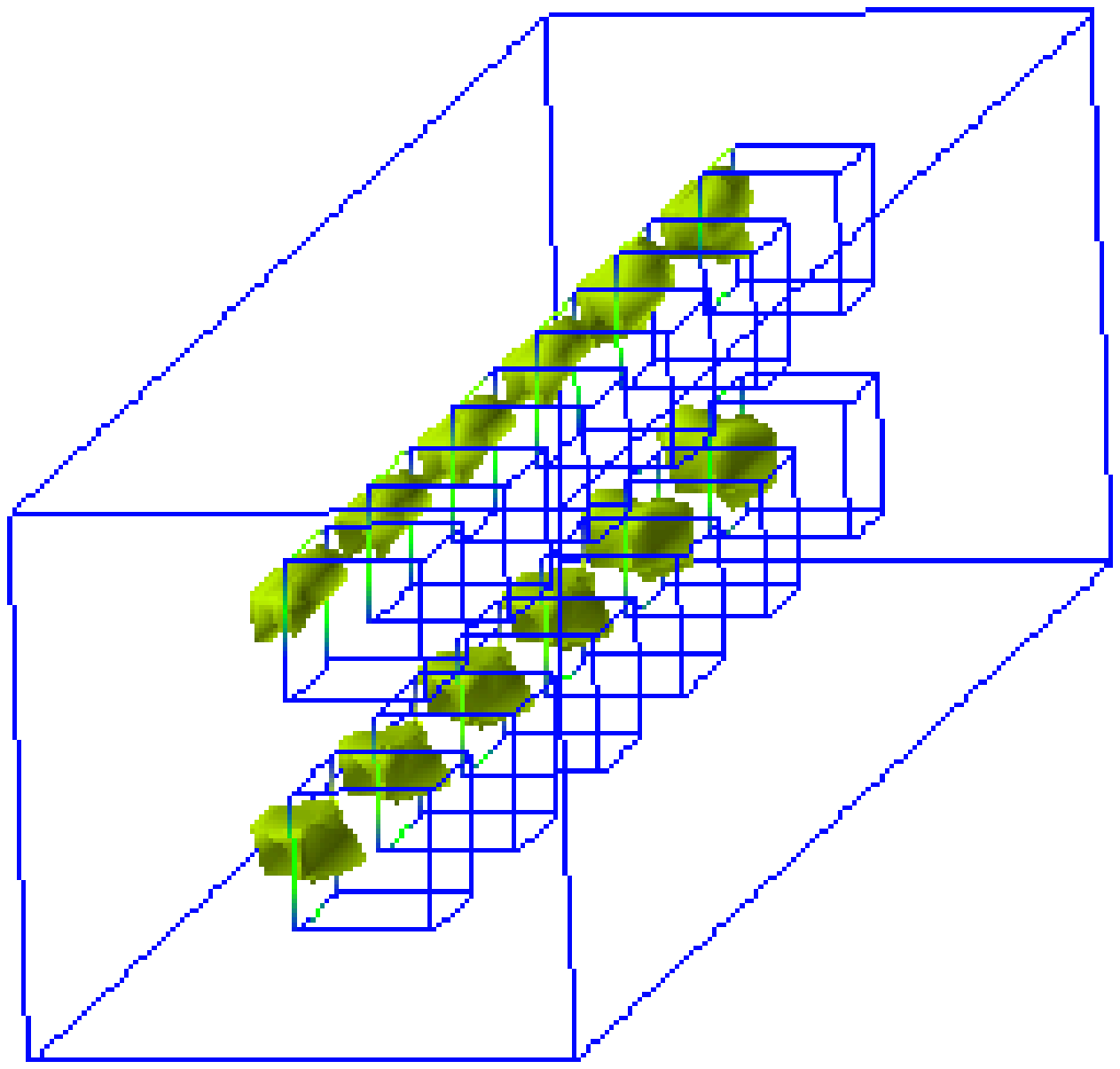} &  
      \includegraphics[width =0.3\textwidth, clip = true, trim = 5.50cm 6.0cm 2.0cm 1.50cm, angle = -90]{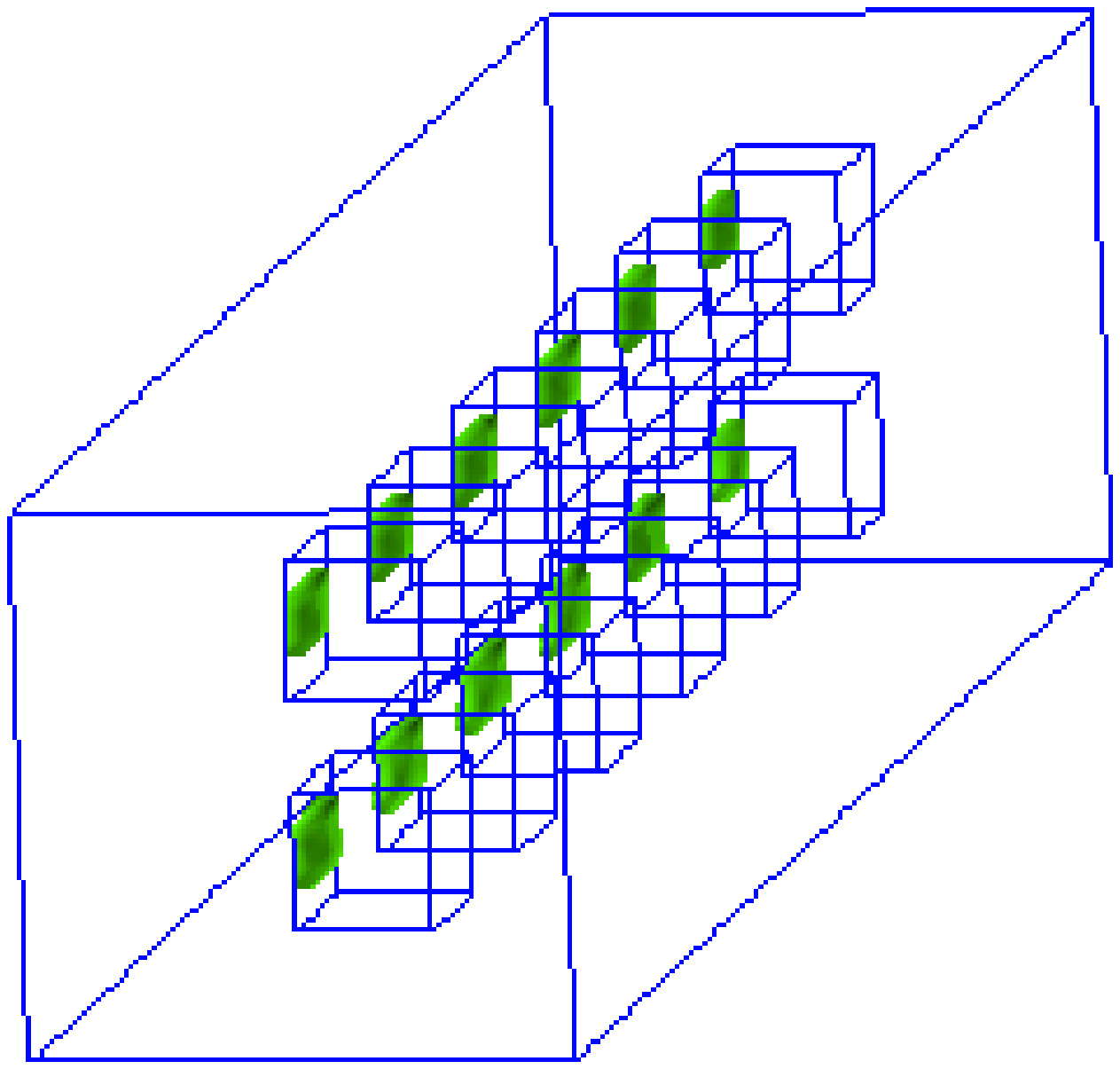}  &\\
   \hspace{-1.5cm} (a)  $\max\limits_{\Omega_{FEM} } \varepsilon_r \approx 15$  &  \hspace{-1.5cm}(b) $\max\limits_{\Omega_{FEM} } \mu_r \approx 2.5$ \\
 \includegraphics[width = 0.3\textwidth, clip = true, trim = 5.50cm 3.0cm 2.0cm 1.50cm,, angle = -90]{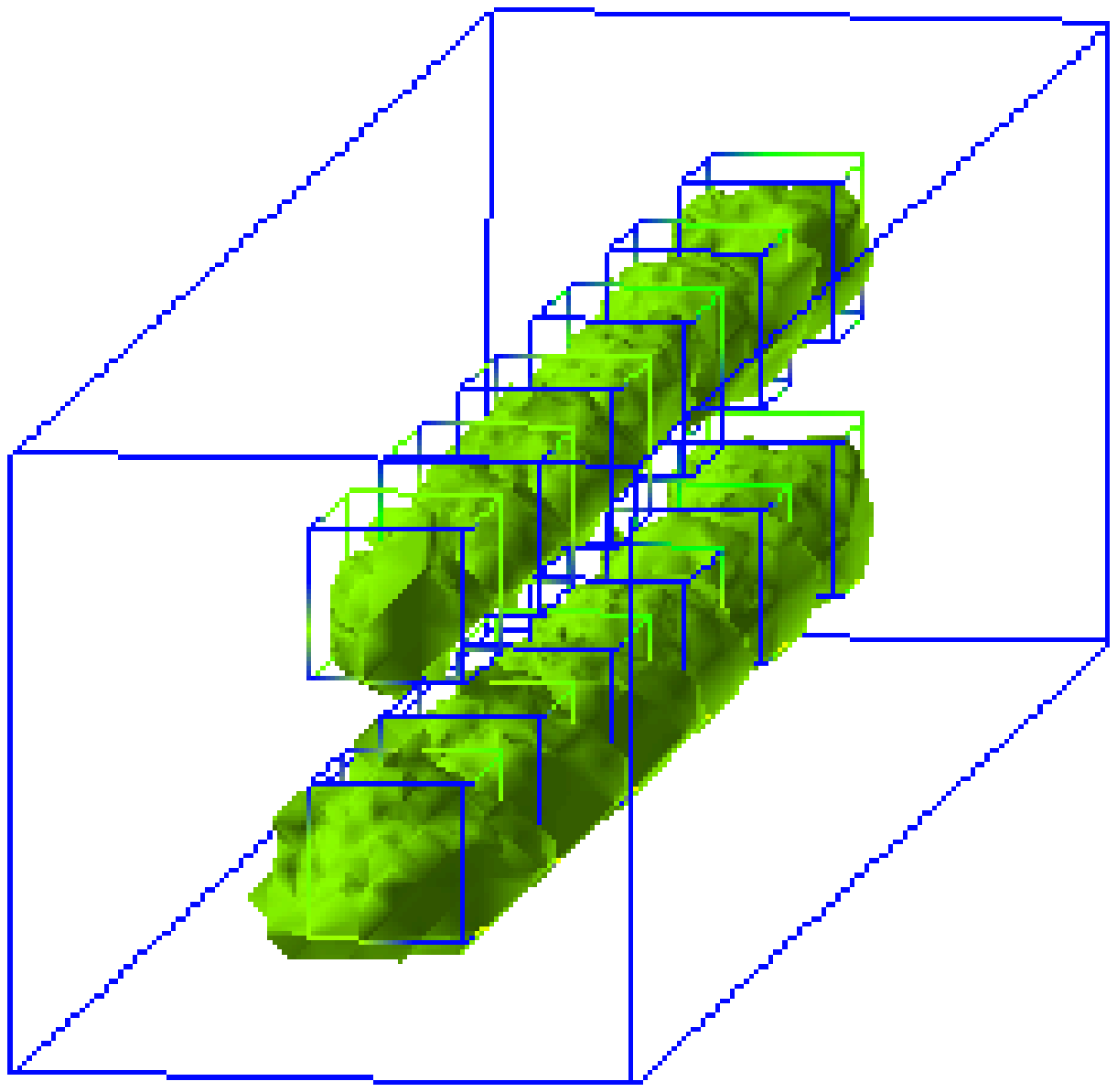} &  
      \includegraphics[width =0.3\textwidth, clip = true, trim = 5.50cm 6.0cm 2.0cm 1.50cm, angle = -90]{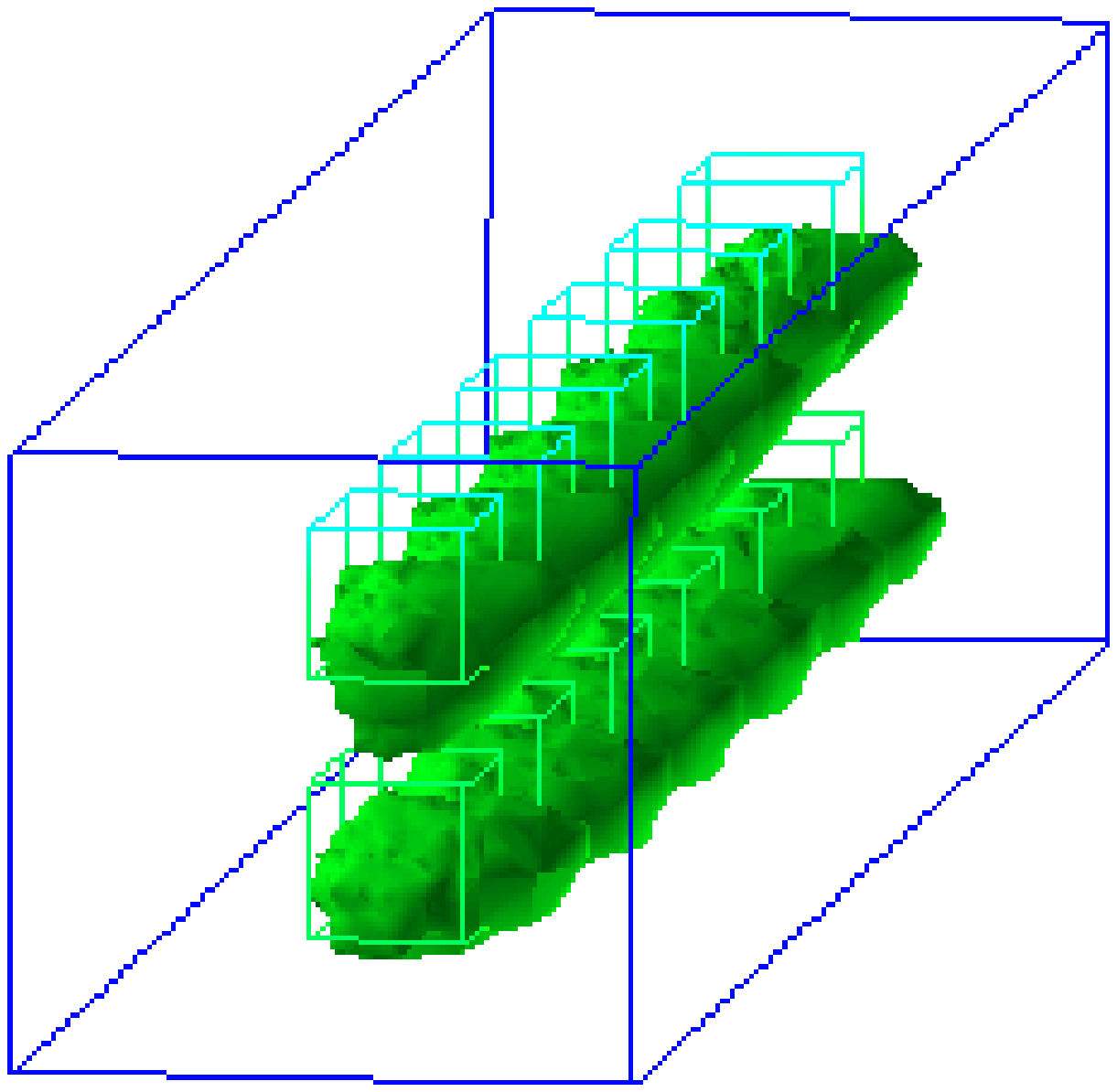}  &\\
   \hspace{-1.5cm} (c)  $\max\limits_{\Omega_{FEM} } \varepsilon_r  \approx 14.9$ &  \hspace{-1.5cm}(d) $\max\limits_{\Omega_{FEM} } \mu_r \approx 1.8$ \\   
   
\end{tabular}   
   
    \caption{Test 1.  Computed images of reconstructed functions $\varepsilon_r(\mathbf{x}) \ \text{and} \ \mu_r(\mathbf{x})$  in $x_2 x_3$ view: a), b) on a coarse mesh, c), d) on a 5 times adaptively refined mesh.  Computations are done for $\omega = 45$, $\sigma =7\%$.}
\label{fig:yz_reconstruction_zero} 

\end{center}

  \end{figure}

\begin{figure}
  \centering
    \begin{tabular}{c c c c }
     
      \includegraphics[width = 0.4\textwidth, clip = true, trim = 0cm 6cm 0cm 2cm, angle = 00]{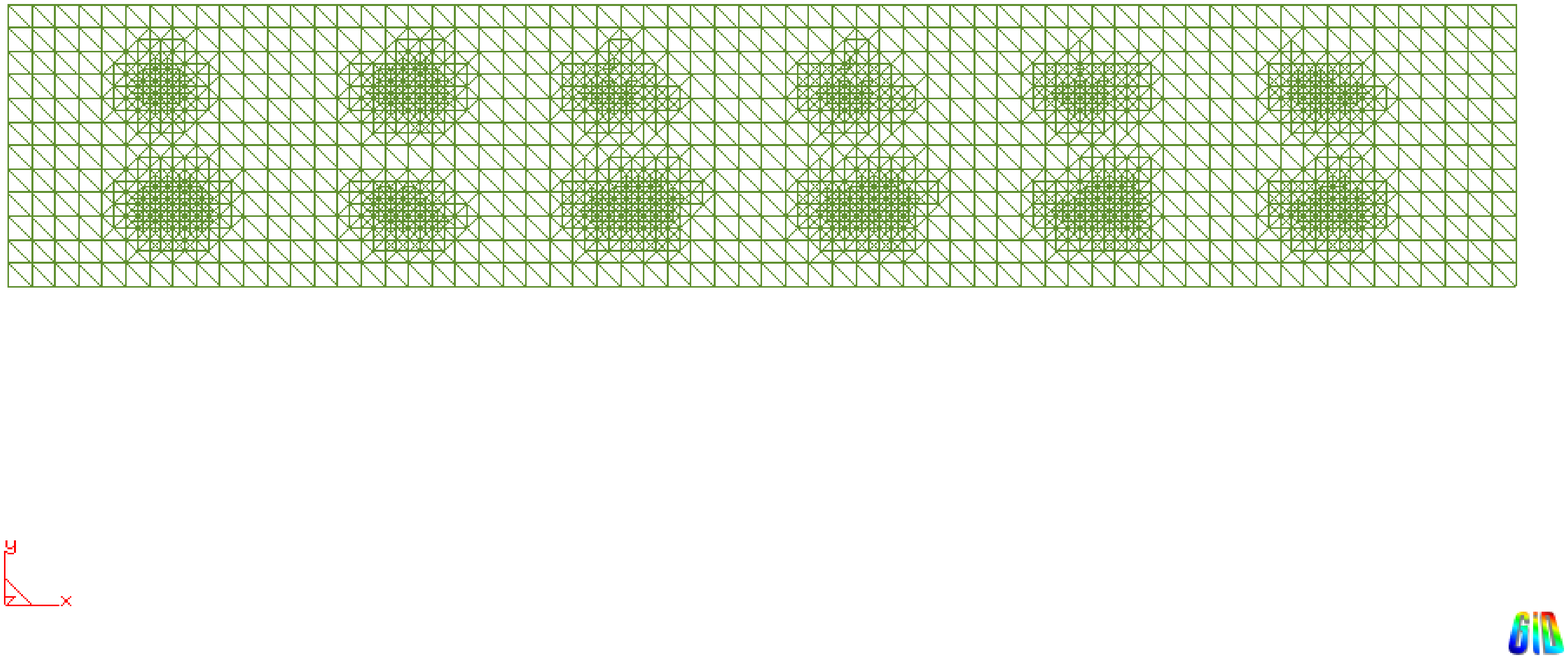} &
       \includegraphics[width = 0.4\textwidth, clip = true, trim = 1cm 6.5cm 1.0cm 3.0cm, angle = 0]{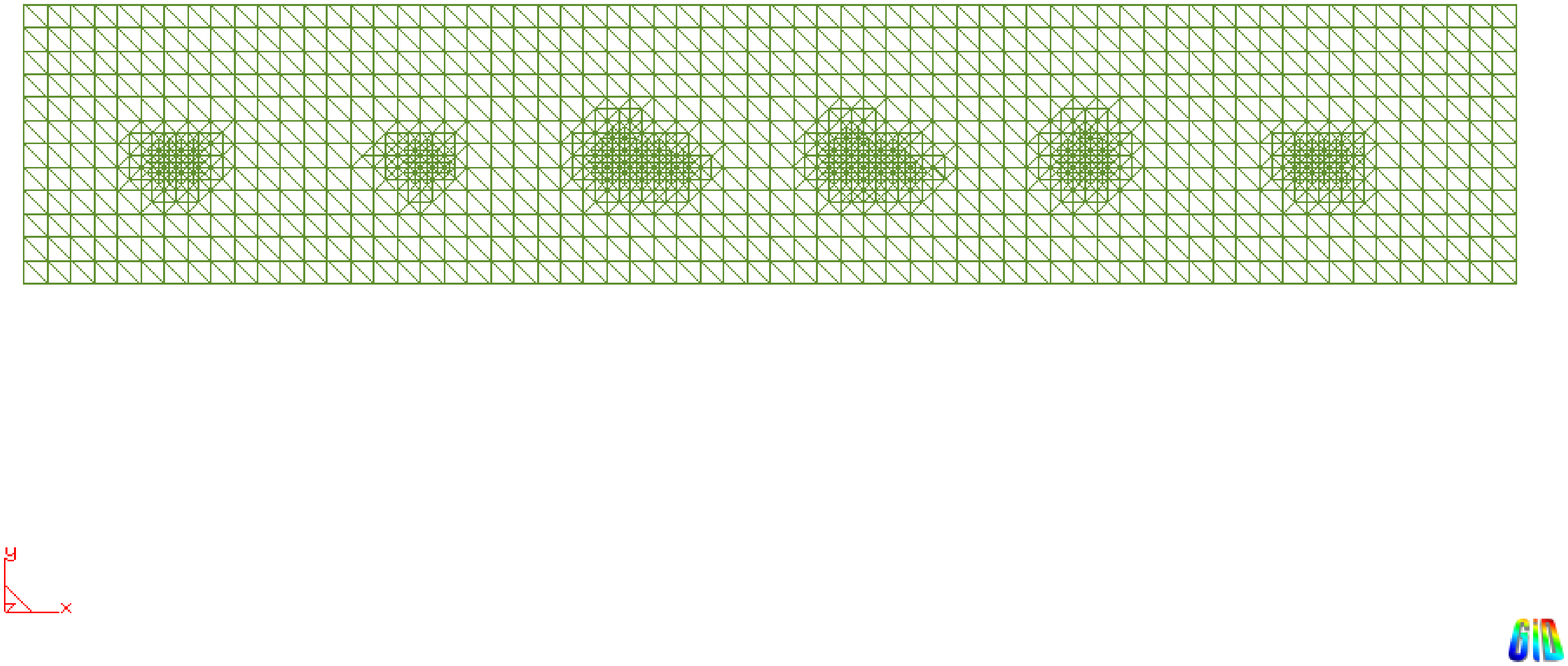} \\
      (a) $x_1 x_2$-view  & (b)  $x_1 x_2$-view \\ 
       \includegraphics[width = 0.4\textwidth, clip = true,  trim = 0cm 6cm 0cm 2cm, angle = 0]{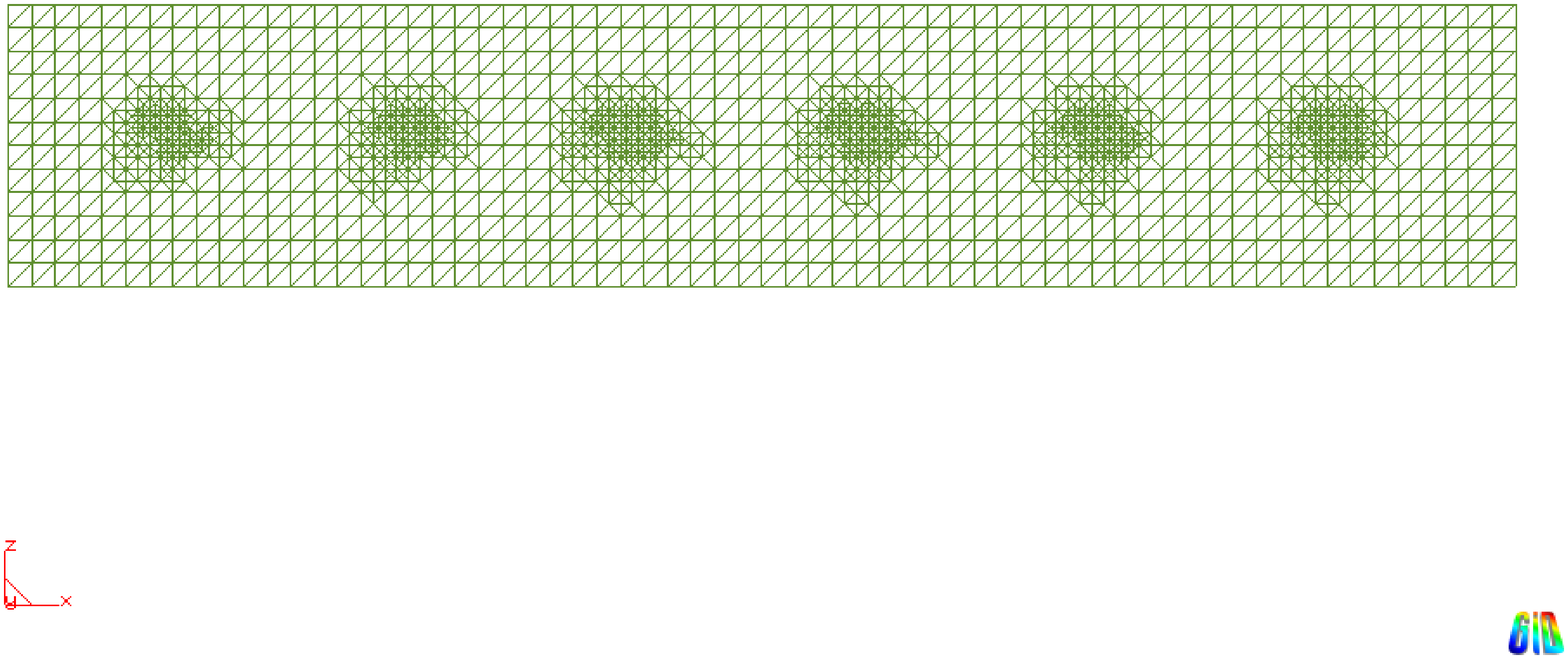} &
       \includegraphics[width = 0.4\textwidth, clip = true, trim = 1cm 6.5cm 1.0cm 3.0cm, angle = 0]{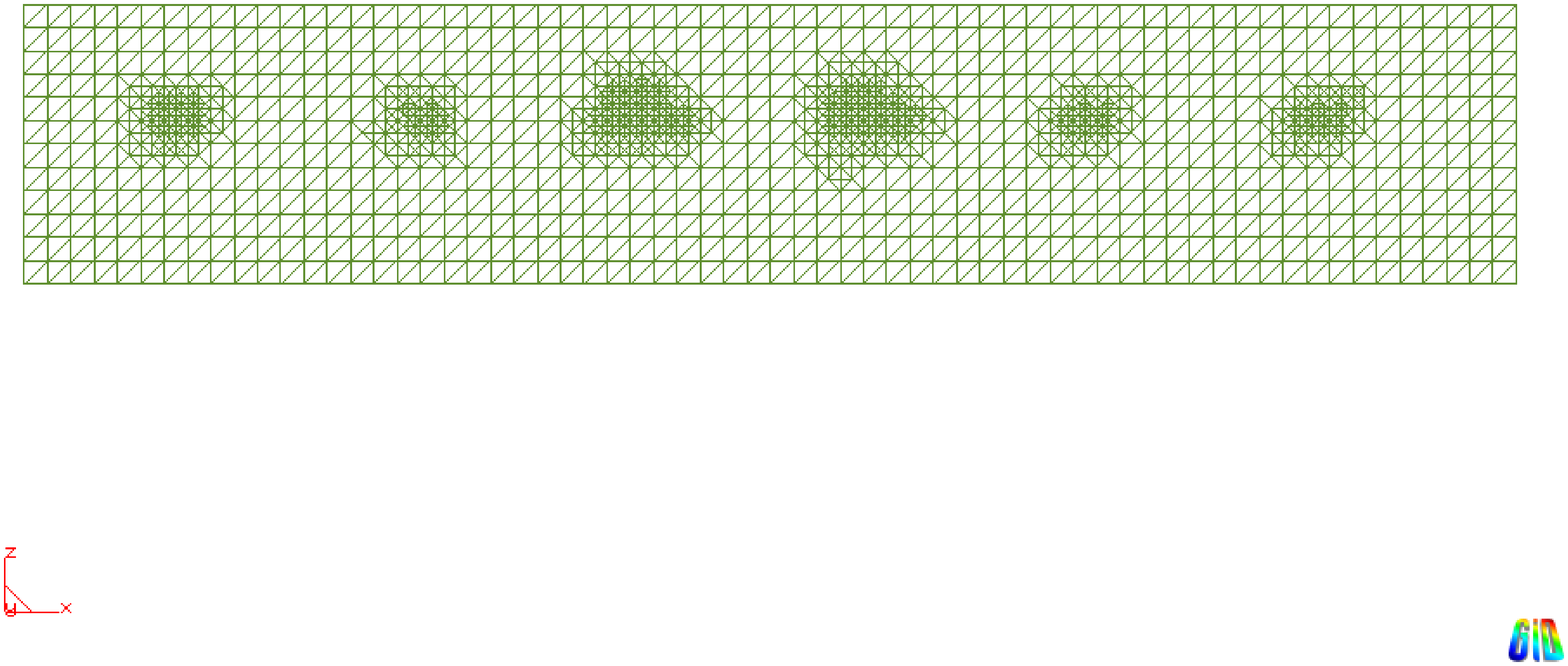} \\
         (c)  $x_1 x_3$-view  & (d)  $x_1 x_3$-view \\
       \includegraphics[width = 0.4\textwidth, clip = true, trim = 1cm 6.5cm 1.0cm 3.0cm, angle = 0]{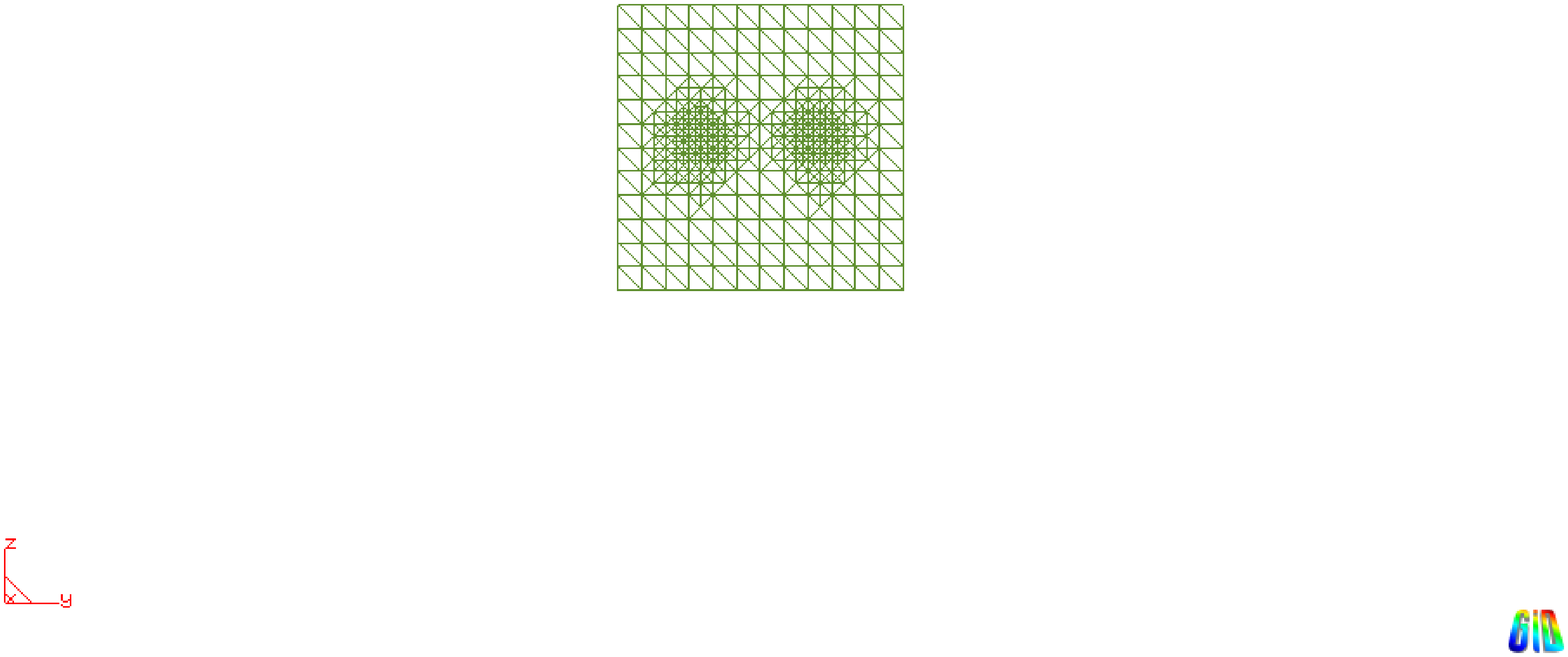} &
       \includegraphics[width = 0.4\textwidth, clip = true, trim = 1cm 6.5cm 1.0cm 3.0cm, angle = 0]{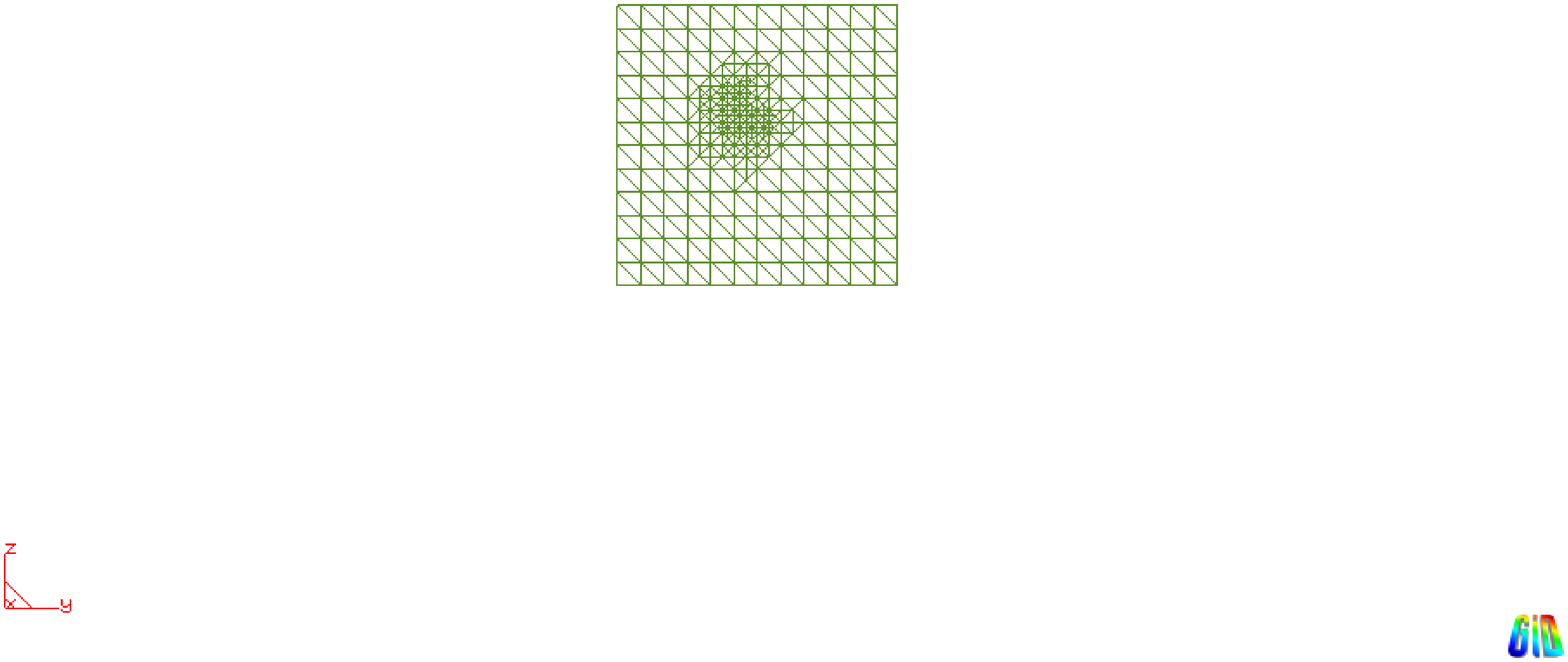} \\
     
       (e) $x_2 x_3$-view  & (f)  $x_2 x_3$-view \\
      \end{tabular}
\caption{Different projections of 5 times adaptively refined meshes for
  computed images of Figures \ref{fig:Recos_omega45noise7_zero} (on
  the left) and Figures \ref{fig:Recos_omega50noise17_nonzero} (on the
  right), respectively.}
       \label{fig:Recos_omega45noise7_zero_mesh}
       \end{figure}

\begin{figure}
        \begin{center}
        
 \begin{tabular}{c c c}
 {\includegraphics[width =2.1cm, clip = true, trim = 9.0cm 4.0cm 6.0cm 4.0cm, angle = -90.0]{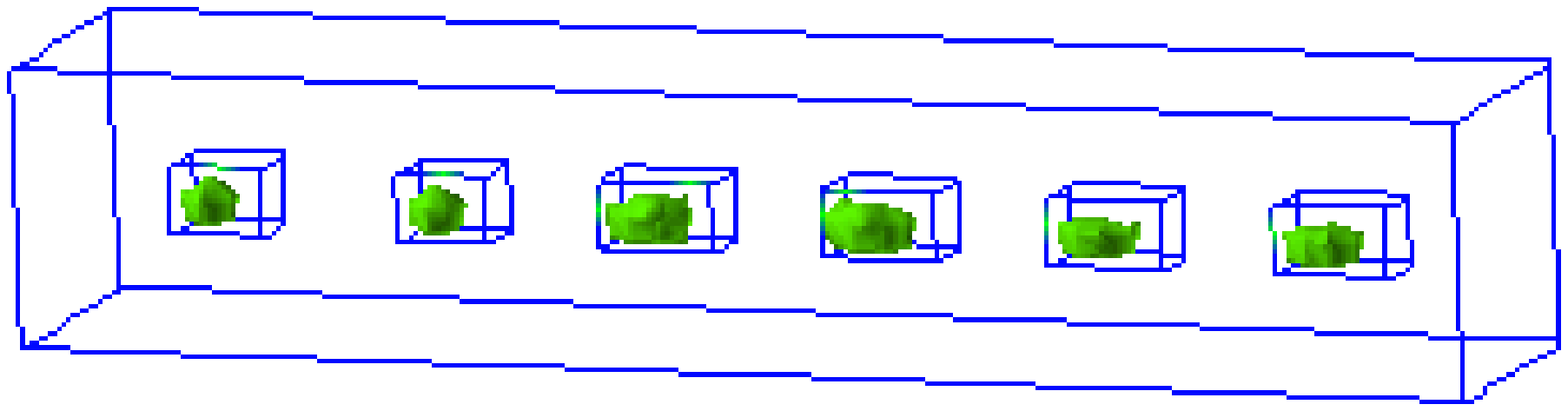}}
 & {\includegraphics[width = 2.1cm, clip = true, trim = 9.0cm 4.0cm 6.0cm 4.0cm, angle = -90.0]{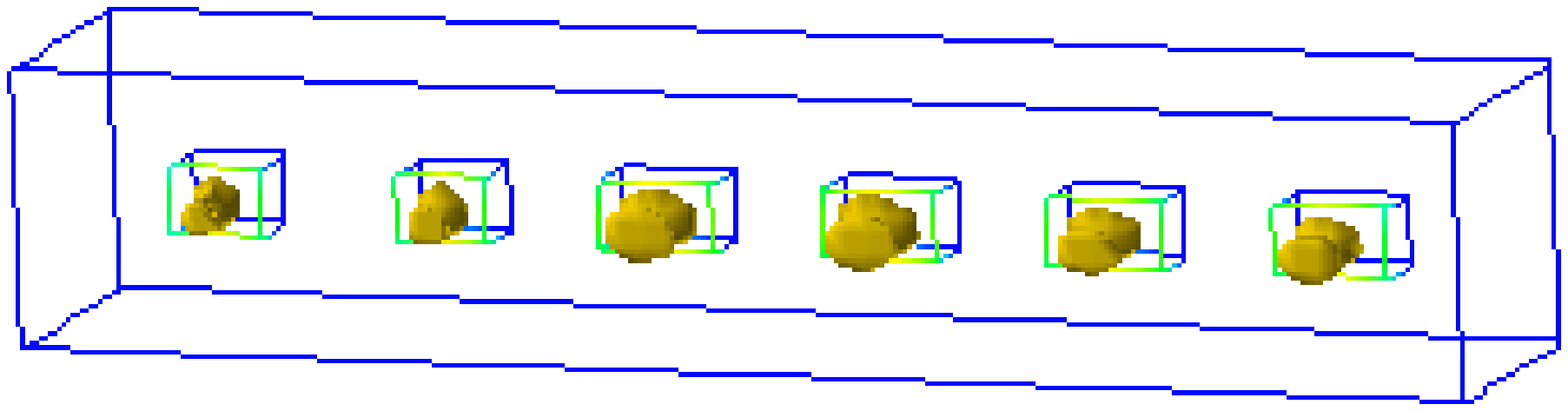}} & \\
c) $\max\limits_{\Omega_{FEM} }\varepsilon_r \approx 15$   & 
d)$\max\limits_{\Omega_{FEM} }\mu_r \approx 2.2$
  \end{tabular}
  
 \caption{Test 2. Computed images of reconstructed functions $\varepsilon_r(\mathbf{x}) \ \text{and} \ \mu_r(\mathbf{x})$ on a coarse mesh for $\omega = 50$, $\sigma =17\%$.}
 \label{fig:Recos_omega50noise17_nonzero_coarse}

\end{center}  
 \end{figure}

\begin{figure}
      
 \begin{tabular}{c c c c }
 {\includegraphics[width =2.1cm, clip = true, trim = 9.0cm 4.0cm 6.0cm 4.0cm, angle = -90.0]{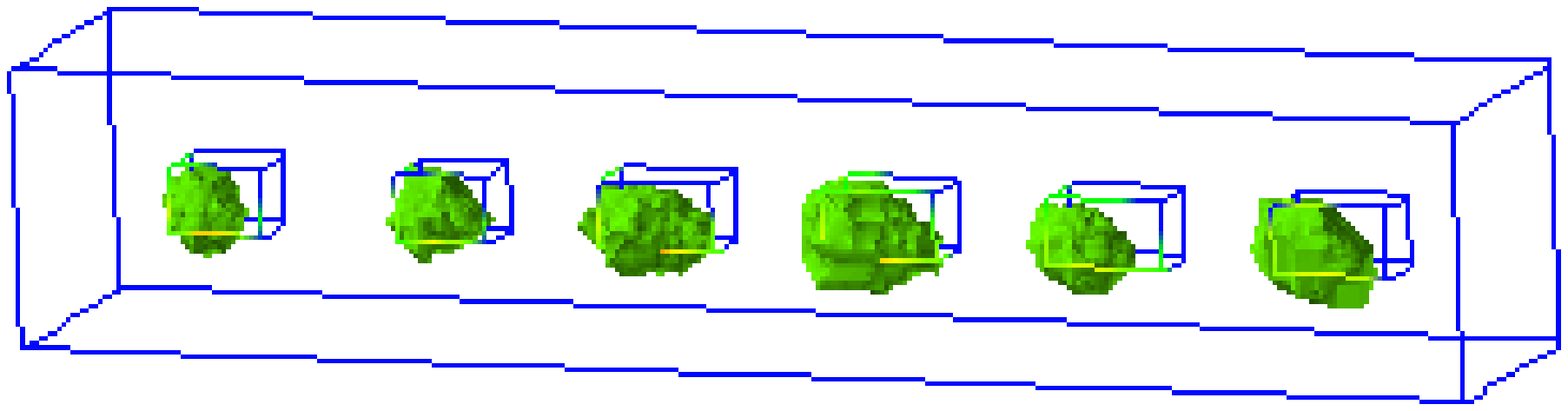}}
 & {\includegraphics[width =2.1cm, clip = true, trim = 9.0cm 4.0cm 6.0cm 4.0cm, angle = -90.0]{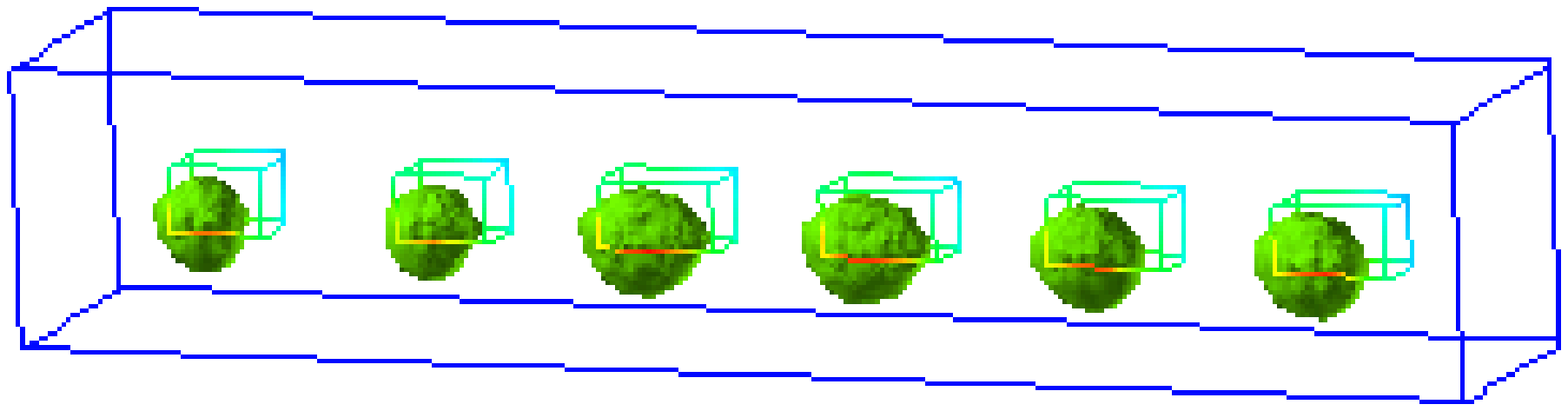}} &  \\
 (a) $\max\limits_{\Omega_{FEM} }\varepsilon_r \approx 14.4$  &
 (b) $\max\limits_{\Omega_{FEM} } \mu_r \approx 1.6$
  \end{tabular}
 \caption{Test 2. Computed images of reconstructed functions $\varepsilon_r(\mathbf{x}) \ \text{and} \ \mu_r(\mathbf{x})$ on a 5 times adaptively refined mesh. Computations are done with  $\omega = 50$, $\sigma = 17\%$.}
 \label{fig:Recos_omega50noise17_nonzero}

 \end{figure}

\begin{figure}
   \begin{center}
   \begin{tabular} {c c c }
   \includegraphics[width = 0.3\textwidth, clip = true, trim = 5.50cm 3.0cm 2.0cm 1.50cm, angle = -90]{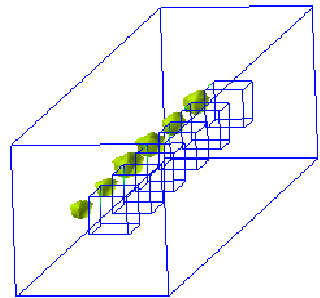}&
      \includegraphics[width = 0.3\textwidth, clip = true, trim = 5.50cm 6.0cm 2.0cm 1.50cm, angle = -90]{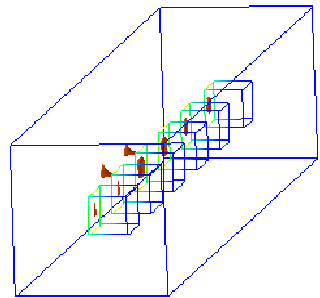} & \\
   \hspace{-1.5cm} (a)  $\max\limits_{\Omega_{FEM} } \varepsilon_r \approx 15$  &  \hspace{-1.5cm}(b) $\max\limits_{\Omega_{FEM} } \mu_r \approx 2.2$ \\   
      \includegraphics[width = 0.3\textwidth, clip = true, trim = 5.50cm 3.0cm 2.0cm 1.50cm, angle = -90]{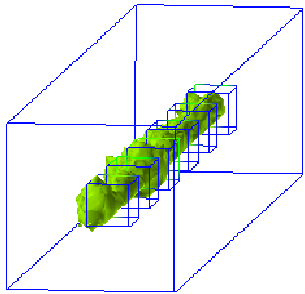}&
      \includegraphics[width = 0.3\textwidth, clip = true, trim = 5.50cm 6.0cm 2.0cm 1.50cm, angle = -90]{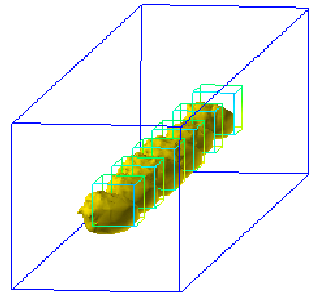} & \\
 \hspace{-1.5cm} (a)  $\max\limits_{\Omega_{FEM} } \varepsilon_r \approx 14.4$  &  \hspace{-1.5cm}(b) $\max\limits_{\Omega_{FEM} } \mu_r \approx 1.6$ \\

 \end{tabular}
 
     \caption{Test 2.  Computed images of reconstructed functions $\varepsilon_r(\mathbf{x}) \ \text{and} \ \mu_r(\mathbf{x})$ in $x_2 x_3$ view: a), b) on a coarse mesh, c), d) on a 5 times adaptively refined mesh. Computations are done  for  $\omega = 50$, $\sigma =17\%$.}
\label{fig:yz_reconstruction}

\end{center}
  \end{figure}

 Numerically we tested our algorithms with two different noise levels,
 $\sigma= 7 \%$ and $\sigma= 17\%$, on the frequency band $\omega \in
 [45, 60]$.  Main conclusion of our previous study of \cite{BCN} was
 that we could get the large contrast of the dielectric function
 $\varepsilon_r$ which allows us to reconstruct metallic targets, and
 that the contrast for $\mu_r$ was within limits of (\ref{admpar}).
 However, the size of $\mu_r$ in $x_1, x_2$ directions and location of all
 inclusions in $x_3$ direction should be improved. Using Figures
 \ref{fig:Recos_omega45noise7_zero_coarse},
 \ref{fig:Recos_omega50noise17_nonzero_coarse} and Tables 1-6 of this
 note we can conclude that on the coarse mesh we get similar results
 as were obtained in \cite{BCN}.  However, with mesh refinements, as
 was expected, quality of reconstruction is improved a lot,  see
 Figures \ref{fig:yz_reconstruction_zero},
 \ref{fig:Recos_omega50noise17_nonzero},
 \ref{fig:yz_reconstruction}.
 Using these Figures and Tables 1-6 we observe that now all
 inclusions have correct locations in $x_3$ direction as well as their
 contrasts and sizes in $x_1, x_2$ directions are also improved and
 reconstructed with a good accuracy.
We can conclude, that we have supported tests of our previous works
\cite{BJ, BMaxwell2, BTKB,BKK, KBB} and have shown that the adaptive
finite element method is very powerful tool for the reconstruction of
heterogeneous targets, their locations and shapes accurately.

Our adaptive algorithms can also be applied for the case when edge
elements are used for the numerical simulation of the solutions of
forward and adjoint problems, see \cite{CWZ99, CWZ14, FJZ10} for
finite element analysis in this case.  This as well as development of
iteratively regularized adaptive finite element method can be
considered as a challenge for the future research.




\newpage

\vspace{1cm}

\section*{Acknowledgments}

This research is supported by the
Swedish Research Council (VR).
 The computations were performed on
resources at Chalmers Centre for Computational Science and Engineering
(C3SE) provided by the Swedish National Infrastructure for Computing
(SNIC).

\medskip
\medskip


\begin{thebibliography}{99}

  \raggedright
\bibitem{Ass} 
\newblock F.~Assous, P.~Degond, E.~Heintze  and P.~Raviart,
\newblock   On a finite-element method for solving the   three-dimensional    Maxwell equations, 
\newblock \emph{Journal of Computational Physics}, {\bf 109} (1993), 222--237.


 \bibitem{BKS}  
 \newblock A. Bakushinsky, M. Y. Kokurin, A. Smirnova,  
 \newblock \emph{Iterative Methods for Ill-posed Problems},
 \newblock Inverse and Ill-Posed Problems Series 54, De Gruyter, 2011.

  \bibitem{BJ} 
 \newblock L. Beilina and C. Johnson, 
 \newblock A posteriori error estimation in computational inverse    scattering, 
 \newblock \emph{Mathematical Models in Applied Sciences}, {\bf 1} (2005), 23-35.


 \bibitem{BMaxwell} 
 \newblock   L. Beilina,  
 \newblock Energy estimates and numerical  verification of the
 stabilized Domain Decomposition Finite  Element/Finite Difference    approach for time-dependent Maxwell's  system, 
 \newblock \emph{Central European Journal of Mathematics},  {\bf 11} (2013), 702-733
 DOI: 10.2478/s11533-013-0202-3.


 \bibitem{BMaxwell2} 
 \newblock L. Beilina,  
 \newblock Adaptive Finite Element Method for a  coefficient    inverse problem for the Maxwell's system,  
 \newblock \emph{Applicable Analysis}, {\bf 90} (2011) 1461-1479.


\bibitem{BCN}  L. Beilina, M. Cristofol and K. Niinim\"aki, Optimization
  approach for the simultaneous reconstruction of the dielectric
  permittivity and magnetic permeability functions from limited
  observations, \emph{Inverse Problems and Imaging}, 9 (1), pp. 1-25, 2015



 \bibitem{BTKB} 
\newblock L. Beilina, N.~T. Th\`anh, M.~V. Klibanov,
   J. Bondestam-Malmberg, 
\newblock Reconstruction of shapes and refractive indices from blind
backscattering experimental  data using the adaptivity, 
\newblock \emph{Inverse Problems},  30, 105007, 2014.

 \bibitem{BTKF} 
 \newblock L. Beilina, N. T. Th\`anh, M. V. Klibanov, M. A. Fiddy, 
 \newblock Reconstruction from blind  experimental data for an    inverse problem for a hyperbolic equation,
 \newblock \emph{Inverse Problems}, {\bf 30} (2014), 025002.


  \bibitem{BOOK}  L. Beilina, M.V. Klibanov,
    \emph{Approximate global convergence and adaptivity for coefficient inverse problems},
    Springer, New-York, 2012.

\bibitem{BKK} L. Beilina, M.V. Klibanov and M.Yu. Kokurin, Adaptivity with relaxation
for ill-posed problems and global convergence for a coefficient inverse
problem, \emph{Journal of Mathematical Sciences}, 167, pp. 279-325, 2010.


\bibitem{BG} 
 \newblock  L.~Beilina, M. Grote, 
  \newblock Adaptive Hybrid Finite   Element/Difference Method for {M}axwell's equations, 
   \newblock \emph{TWMS Journal  of Pure and Applied Mathematics}, V.1(2), pp.176-197,  2010.



\bibitem{BCS}  M. Bellassoued, M. Cristofol, and E. Soccorsi,
\newblock  Inverse boundary value problem for the dynamical heterogeneous Maxwell's system,
\newblock \emph{Inverse Problems}, (28), 095009,  2012.
   

\bibitem{Bondestam1}  J. Bondestam Malmberg,
 \newblock A posteriori error estimate in the Lagrangian setting for an inverse problem based on a new
formulation of Maxwell’s system,
 \newblock \emph{Inverse Problems and Applications, Springer Proceedings in Mathematics and Statistics }, Vol. 120, 2015, pp. 42-53, 2015.


\bibitem{Bondestam2} J. Bondestam Malmberg, 
\newblock   A posteriori error estimation
in a finite element method for reconstruction of dielectric
permittivity,
\newblock  arXiv:1502.07658, 2015.



\bibitem{Brenner} S. C. Brenner  and L. R. Scott, \emph{The Mathematical Theory of Finite
  Element Methods} (2nd edn), Springer-Verlag, New York, 2002.

   
   
       \bibitem{CWZ99} 
 \newblock P. Ciarlet, Jr. and J. Zou, 
 \newblock Fully discrete finite element approaches for time-dependent
 Maxwell's equations,
 \newblock  \emph{Numerische Mathematik}  {\bf 82} (1999), 
193-219.

   
    \bibitem{CWZ14}
 \newblock P. Ciarlet, Jr., H. Wu and J. Zou, 
 \newblock Edge element methods for Maxwell's equations with strong convergence for Gauss' laws,
 \newblock  \emph{SIAM Journal on Numerical Analysis} {\bf 52} (2014), 
779-807.



  \bibitem{CFL67}
 \newblock R. Courant, K. Friedrichs and H. Lewy 
 \newblock  On the     partial differential equations of    mathematical physics, 
 \newblock \emph{IBM Journal of Research and Development}, {\bf 11} (1967), 215-234.

 \bibitem{Cohen} 
 \newblock G. C. Cohen, 
 \newblock \emph{Higher order numerical methods for  transient wave  equations}, 
 \newblock Springer-Verlag, 2002.

\bibitem{delta} 
\newblock A.~Elmkies and P.~Joly,  
\newblock Finite elements and mass lumping for Maxwell's equations:
the 2D case. 
\newblock\emph{Numerical Analysis}, {\bf 324} (1997), 1287--1293.



 \bibitem{Engl} 
 \newblock H. W. Engl,  M. Hanke and A. Neubauer,  
 \newblock \emph{Regularization of Inverse Problems} 
 \newblock Boston: Kluwer Academic Publishers, 2000.


 \bibitem{EM}  
 \newblock B. Engquist and A. Majda, 
 \newblock Absorbing boundary conditions for the numerical    simulation of waves,
 \newblock  \emph{Mathematics of Computation}, {\bf 31} (1977), 629-651.
 
  \bibitem{FJZ10}
 \newblock  H. Feng, D. Jiang and J. Zou, 
 \newblock Simultaneous identification of electric permittivity and magnetic permeability, 
 \newblock \emph{Inverse Problems}, {\bf 26} (2010), 095009.
 

\bibitem{WSKVH09}
\newblock D. W. Einters, J. D. Shea, P. Komas, B. D. Van Veen and
S. C. Hagness
\newblock Three-dimensional microwave breast imaging: Dispersive
dielectric properties estimation using patient-specific basis functions
\newblock \emph{IEEE Transactions on Medical Imaging} {\bf 28} (2009),
969--981

\bibitem{EEJ} K. Eriksson, D. Estep and C.\ Johnson,\emph{\ Calculus in Several
Dimensions}, Springer, Berlin, 2004.



    \bibitem{IJT11}
 \newblock  K~Ito, B~Jin, and T~Takeuchi
  \newblock Multi-parameter Tikhonov regularization,
  \newblock \emph{Methods and Applications of Analysis}, {\bf 18} (2011), 31--46.
  
\bibitem{JS} C. Johnson and A. Szepessy, Adaptive finite element methods for
conservation laws based on a posteriori error estimation, \emph{Comm. Pure
Appl. Math}., 48, 199--234, 1995.

\bibitem{Johnson} C. Johnson, \emph{Numerical Solution of Partial
  Differential Equations by the Finite Element Method}, Dover Books on
  Mathematics, 2009.

\bibitem{KBB} M.V. Klibanov, A.B. Bakushinsky and L. Beilina, Why a
minimizer of the Tikhonov functional is closer to the exact solution than
the first guess,\emph{\ J.\ Inverse and Ill-Posed Problems}, 19, 83-105,
2011.

\bibitem{lad} O. A. Ladyzhenskaya, \emph{Boundary Value Problems of
  Mathematical Physics}, Springer Verlag, Berlin, 1985.

 \bibitem{div_cor} 
 \newblock C.~D.~Munz, P.~Omnes, R.~Schneider, E.~Sonnendrucker
   and U.~Voss,  
 \newblock Divergence correction techniques for Maxwell   Solvers   based on a hyperbolic model, 
 \newblock \emph{Journal of Computational Physics}, {\bf 161} (2000),
   484--511.


 \bibitem{Peron} 
 \newblock O. Pironneau, 
 \newblock \emph{Optimal shape design for elliptic  systems},
 \newblock  Springer-Verlag, Berlin, 1984.

 \bibitem{NBKF} 
 \newblock N. T. Th\`anh, L. Beilina, M. V. Klibanov,
   M. A. Fiddy, 
 \newblock Reconstruction of the refractive index from experimental backscattering data using a globally convergent inverse method,
 \newblock \emph{SIAM J. Scientific Computing},36 (3)  (2014), 273-293.

\bibitem{samar} S. Hosseinzadegan, Iteratively regularized adaptive finite element method for reconstruction of coefficients in Maxwell's system, \emph{Master's Thesis in Applied Mathematics}, Department of Mathematical Sciences, Chalmers University of Technology and Gothenburg University, 2015.

\bibitem{SZ}  L. R. Scott, S. Zhang,  
\newblock Finite element interpolation of nonsmooth functions satisfying boundary conditions,
  \newblock\emph{Math.Comp.}, {\bf 54} (1990),
   483--493.

 \bibitem{SSMS} 
 \newblock D. R. Smith, S. Schultz, P. Markos and C. M. Soukoulis,
 \newblock  Determination of effective permittivity and permeability of  metamaterials from reflection and transmission coefficients,
 \newblock \emph{Physical Review B}, {\bf 65} (2002), DOI:10.1103/PhysRevB.65.195104.


 \bibitem{petsc}  
 \newblock PETSc, Portable, Extensible Toolkit for
 Scientific Computation, http://www.mcs.anl.gov/petsc/.


\bibitem{tikhonov} A.N. Tikhonov, A.V. Goncharsky, V.V. Stepanov and A.G. Yagola, 
\emph{Numerical Methods for the Solution of Ill-Posed Problems}, London:
Kluwer, London, 1995.

\bibitem{TGSK} Tikhonov A.N., Goncharskiy A.V., Stepanov V.V., Kochikov I.V. Ill-posed problems of the image processing, DAN USSR, 294, № 4. 832-837, 1987.


\bibitem{waves}  WavES, the software package, http://www.waves24.com


\bibitem{ZTRL11}
\newblock X. Zhang, H. Tortel, S. Ruy and A. Litman
\newblock Microwave imaging of soil water diffusion using the  linear sampling method,
\newblock \emph{IEEE Geoscience and Remote Sensing Letters} {\bf 8} (2011),421--425



\end{thebibliography}
\end{document}